\documentclass[10pt]{article}
\usepackage[textwidth=15cm, textheight=22cm]{geometry}
\usepackage{amsfonts}
\usepackage{amsmath}
\usepackage{amssymb}
\usepackage{amsthm}
\usepackage{enumitem}

\usepackage[dvipsnames]{xcolor}

\usepackage{tikz}
\usetikzlibrary{arrows.meta,shapes,calc}
\usetikzlibrary{positioning,shapes,arrows}

\usepackage{pgfplots}
\pgfplotsset{compat=1.12}
\usetikzlibrary{plotmarks}
\usepackage[ngerman, english]{babel}
\usepackage{babelbib}
\usepackage[utf8]{inputenc}

\usepackage{nicefrac}

\usepackage{subfig}
\usepackage{epstopdf}
\usepackage{bbm}

\usepackage{url}

\usepackage[colorinlistoftodos,prependcaption,textsize=tiny]{todonotes}

\usepackage{kpfonts}
\usepackage{adjustbox}

\usepackage{algorithmic}
\usepackage[section,Algorithm]{algorithm}

\usepackage[percent]{overpic}

\usepackage{graphicx}%
\graphicspath{{./img/}}

\newcommand{\AAA}{\mathcal{A}}

\newcommand{\MM}{\mathcal{M}}

\newcommand{\OO}{\mathcal{O}}
\newcommand{\PP}{\mathcal{P}}

\newcommand{\UU}{\mathcal{U}}

\newcommand{\E}{\mathbb{E}}

\newcommand{\N}{\mathbb{N}}

\newcommand{\R}{\mathbb{R}}

\newcommand{\dx}{\dot{x}}
\newcommand{\dz}{\dot{z}}
\newcommand{\dy}{\dot{y}}
\newcommand{\du}{\dot{u}}
\newcommand{\dv}{\dot{v}}
\newcommand{\dw}{\dot{w}}

\newcommand{\drho}{\dot{\rho}}

\newcommand{\Pol}{\Phi} 
\newcommand{\legPol}{P}
\newcommand{\g}{\tilde{g}}
\newcommand{\coeff}{c} 
\newcommand{\nmoms}{n_{\text{moms}}}
\newcommand{\npara}{n_{\text{para}}}

\newcommand{\ind}{\mathbbmss{1}}

\DeclareMathOperator{\Beta}{Beta}
\DeclareMathOperator{\Gamdist}{\text{Gamma}}
\DeclareMathOperator{\argmin}{arg\,min}

\theoremstyle{plain}
\newtheorem{thm}{Theorem}[section]

\newtheorem{prop}[thm]{Proposition}

\theoremstyle{definition}
\newtheorem{defin}{Definition}[section]
\newtheorem{ex}{Example}[section]

\theoremstyle{remark}
\newtheorem{rem}{Remark}

\begin{document}
\title{Uncertainty Quantification of Bifurcations\\ in Random Ordinary Differential Equations}
\author{Christian Kuehn, Kerstin Lux}
\footnotetext{Technical University of Munich, Department of Mathematics, (ckuehn@ma.tum.de, kerstin.lux@tum.de)}
\maketitle


We are concerned with random ordinary differential equations (RODEs). Our main question of interest is how uncertainties in system parameters propagate through the possibly highly nonlinear dynamical system and affect the system's bifurcation behavior. We come up with a methodology to determine the probability of the occurrence of different types of bifurcations (sub- vs super-critical) along a given bifurcation curve based on the probability distribution of the input parameters. In a first step, we reduce the system's behavior to the dynamics on its center manifold. We thereby still capture the major qualitative behavior of the RODEs. In a second step, we analyze the reduced RODEs and quantify the probability of the occurrence of different types of bifurcations based on the (nonlinear) functional appearance of uncertain parameters. To realize this major step, we present three approaches: an analytical one, where the probability can be calculated explicitly based on Mellin transformation and inversion, a semi-analytical one consisting of a combination of the analytical approach with a moment-based numerical estimation procedure, and a particular sampling-based approach using unscented transformation. We complement our new methodology with various numerical examples.

\bigskip

{\bf Keywords.} Uncertainty propagation, nonlinear random ordinary differential equations, Mellin transform, polynomial chaos expansion, method of moments, Gaussian mixture models, unscented transformation

\section{Introduction}

Nonlinear dynamics are omnipresent in various scientific disciplines such as climate science (see e.g.\ \cite{Dijkstra_bifurcations.2019,Palmer.1999}), ecology (see e.g.\ \cite{Lotka.1920}) and epidemiology (see e.g.\ \cite{Kermack.1927}). Numerous real-world phenomena can be described by models built upon nonlinear ordinary differential equations (ODEs). Key characteristics of the model are often encoded via deterministic \emph{parameters}. Yet, in most real-world situations, the parameters have to be estimated via data assimilation techniques~\cite{LawStuartZygalakis}. In particular, input parameters are prone to measurement errors, errors arising in statistical inference methods due to numerical inaccuracies or an insufficient size of the data base or even a lack of measurement data at all making the use of proxies necessary.

Parameter-related \emph{uncertainty} can crucially affect the system's dynamics. In the theory of dynamical systems, it is well-known that a system might undergo a \emph{bifurcation} as a parameter is varied (see \cite{Kuznetsov.2004,Strogatz.2015,Wiggins.2003}). That is the flow of the vector field induced by the ODE changes qualitatively as the parameter exceeds a certain value (cf.\ \cite[Def.\ 20.1.1]{Wiggins.2003}). For some bifurcations, a loss of stability of an equilibrium is linked to a rather smooth qualitative change, in the sense that small perturbations drive the system in a new equilibrium which is close by. Other bifurcation types cause drastic and sudden changes. In the latter case, a small perturbation can lead to a jump to a new equilibrium far off, which we then refer to as a \emph{critical transition}~\cite{Schefferetal}, \emph{tipping point}~\cite{AshwinWieczorekVitoloCox}, \emph{first-order (phase) transitions} or \emph{hard exchange of stability}. Typical bifurcation-induced critical transitions are subcritical pitchfork bifurcations or subcritical Hopf bifurcations~\cite{Kuehn.2013}. Yet, pitchfork and Hopf bifurcations can also occur in supercritical variants, where nearby stable steady states and limit cycles bifurcate respectively, so that no tipping or hard transition occurs in such cases. Hence, even if we know from modelling that a system does exhibit a certain bifurcation, it is crucial to further determine the \emph{precise type}. This situation generically~\cite{KuehnBick} occurs in many applications. It is highly relevant in areas such as climate science, ecology or epidemiology, where a single bifurcation-induced critical transition point in parameter space can have massive global impacts. Clearly, we may not hope to have exact estimates for all parameter values in climate science, ecology or epidemiology.\medskip

The first main insight of our work is to recognize that we desperately need better mathematical tools to rigorously estimate the occurrence of different bifurcation types under uncertainty. In this work, we look at a given curve of a certain bifurcation type. Under the assumption of being on this bifurcation curve, we investigate and develop tools to calculate the probability that the criticality of this bifurcation is such that it induces a critical transition, i.e., we are interested in quantifying when a bifurcation is super-critical or sub-critical. The magnitude of this probability of bifurcation type helps to come up with a quantitative assessment of the risk of a bifurcation induced tipping within the system.\medskip 

From a mathematical point of view, we can build upon bifurcation theory in a purely deterministic setting (see \cite{Guckenheimer.1983,Kuznetsov.2004,Strogatz.2015,Wiggins.2003}). Yet, the quantitative propagation of parametric uncertainty through bifurcation analysis is not covered by classical bifurcation theory. In our work, we present several ways to calculate the probability of having a certain bifurcation type. It is natural to start with bifurcations of codimension one, that is the variation of one parameter is sufficient to cause the bifurcation to happen. In particular, codimension-one phenomena, potentially with additional natural symmetries or constraints, are the most frequently encountered bifurcations in applications~\cite{Guckenheimer.1983,Kuznetsov.2004,Strogatz.2015,Wiggins.2003}.

We approach the problem of uncertain parameters by combining known statistical and probabilistic methods with classical analysis and bifurcation theory. The interplay of the diverse disciplines finally leads to promising first estimates of the bifurcation-type probability. In some cases, the latter will even be exact. We use a two-step procedure. The first step consists of using results from classical bifurcation theory. We know that the type of a given codimension-one bifurcation can be determined by the sign of a suitable normal form coefficient~\cite{Kuznetsov.2004}. This coefficient can be obtained by a sequence of reduction and transformation steps. One has to reduce the system, e.g., via center manifold techniques~\cite{Carr} or Lyapunov-Schmidt reduction~\cite{ChowHale}, or using other similar reduction invariant manifold reduction techniques~\cite{Veraszto.2020}. Then, one obtains a reduced system, which can be used to calculate normal form coefficients of the bifurcation. These coefficients generically involve nonlinear combinations of the input parameters, which are uncertain in our case. In the second step, we have to study the resulting nonlinear and problem-specific bifurcation normal form coefficients via efficient probabilistic techniques from uncertainty quantification. Here, we present three different approaches: an analytical, a semi-analytical, and a sampling-based one. Our main reasoning for considering three different approaches is that it is well-known in applied nonlinear dynamics that, depending on the situation, it can be helpful to have analytical, semi-analytical or formal asymptotic, as well as purely numerical methods.

\smallskip

For the first approach, the analytical one, we make excessive use of the \emph{Mellin transform}, an integral transform known from classical analysis. We exploit its probabilistic interpretation and obtain, for some combinations of uncertain input parameters, the exact analytical expression for the probability distribution of the bifurcation type. Moreover, these analytical expressions even allow for a perturbation analysis to assess the influence of small deviations in the input parameter distributions.

The second approach, the semi-analytical one, is motivated by the concrete example of the Lorenz system, where our purely analytical approach is, in general, not applicable. Aiming to keep as many benefits as possible of the availability of an exact probability distribution, we make use of properties of the Mellin transform to divide the bifurcation normal form coefficient into several parts. Some can still be treated by our analytical approach. A well-designed combination of a \emph{polynomial chaos expansion (PCE)}~\cite{Sullivan.2015} with the Mellin transform allows to reassemble the individual parts of the normal form coefficient of the bifurcation providing us closed-form expressions for its moments. To obtain an estimate of the probability distribution, we use \emph{moment-based approximations} via polynomial approximations and the \emph{method of moments} in the context of \emph{Gaussian mixture models}. This approach gives us a very flexible tool to address the bifurcation probabilities with minor assumptions on the input distributions only.

The first two approaches are beneficial for large, computation-intensive ODE systems as they do not rely on samples at all. Sometimes, it might no longer be clear how to separate the input parameters within the normal form coefficient of the bifurcation, especially for a high number of input parameters. For a first rough estimate, we therefore propose a third approach, which is sampling-based. Only a very low number of samples is needed due to the use of the \textit{unscented transformation}. Therefore, the probabilistic analysis of larger ODE systems is likely to still be feasible.

\bigskip

Our work is structured as follows. In Section \ref{sec:problemFormulation}, we introduce the mathematical notion of a random ODE (RODE), explain our two-step procedure to calculate the bifurcation-type probability, and illustrate the first part of the two-step procedure on the basis of the concrete example of the Lorenz system.

The second part of the two-step procedure is addressed in three different ways in Sections \ref{sec:anaAppUncertaintyProp} --\ref{sec:samplingBasedAppUncertaintyProp}. In Section \ref{sec:anaAppUncertaintyProp}, we present the details of our analytical approach and a perturbation analysis example is given. In Section \ref{sec:semiAnaAppUncertaintyProp}, we extend the analytical approach by the integration of the Mellin transform of a PCE. A numerical case study underlines the performance of our semi-analytical approach. In Section \ref{sec:samplingBasedAppUncertaintyProp}, we present an alternative way of estimating the bifurcation probability by using a sampling-based approach, namely the unscented transformation. We conclude our work with a summary and outlook.

\section{Problem formulation} \label{sec:problemFormulation}

We study systems of ordinary differential equations (ODEs)
\begin{align}
	\frac{dx}{dt} &= \dx = f(x,r), \quad x=x(t) \in \R^n, \label{eq:ODE}
\end{align}
where $f:\R^n\times\R^d \rightarrow \R^n$ is a sufficiently smooth vector field, and $r \in \R^d$ represents the $d \in \N$ given input parameters. There already exist sophisticated methods to analyze such systems in the case where the parameters are assumed to be given deterministic constants (see \cite{Guckenheimer.1983,Kuznetsov.2004,Strogatz.2015,Wiggins.2003}). Here, we extend the system of ODEs \eqref{eq:ODE} by uncertain input parameters leading to a system of \textit{random ordinary differential equations (RODEs)}
\begin{align}
	\dx &= f(x,r(\omega)), \quad x=x(t) \in \R^n, \label{eq:rODE}
\end{align}
where $r(\omega) \in \R^d$ is a random vector on a fixed probability space $(\Omega,\AAA,P)$ representing the $d$ uncertain parameters. Since most real-world systems are highly nonlinear, i.e. the function $f:\R^n\times\R^d \rightarrow \R^n$ is nonlinear. A key challenge is to study the uncertainty propagation through the nonlinear dynamics~\cite{Oberkampf.2004}.

\subsection{Bifurcation and instability} 
\label{sec:BifurcationAndInstab}

Our main interest consists in quantifying the probabilities of the occurrence of certain types of bifurcations in the dynamical system \eqref{eq:rODE} based on the uncertain nature of the input parameters. Many bifurcations are associated with a change in the stability of the system, yet only some of these bifurcations lead generically to large shifts of the system, so-called critical transitions or tipping pints. Some bifurcations always induce large shifts such as fold bifurcations, while for others, it crucially depends upon the precise type, e.g., sub-critical pitchfork bifurcations are tipping points while super-critical pitchfork bifurcations are not~\cite{Kuehn.2010}. For decision-makers, it is crucial to know the risk exposure to the presence of different types of bifurcations once it is known that a model can exhibit a certain bifurcation type such as a pitchfork bifurcation in which case one has to determine a certain normal form coefficient to check whether the bifurcation is sub- or super-critical.

If the system \eqref{eq:rODE} was purely deterministic, it would be natural to use classical bifurcation theory (see e.g.\ \cite{Kuznetsov.2004,Strogatz.2015,Wiggins.2003}). Yet, to the best of our knowledge, there exists no general methodology for calculating the probability of the presence of certain types of bifurcations within the RODE \eqref{eq:rODE} based on the probability distributions of the input parameters so far. In our work, we develop tools for the calculation of these \textit{bifurcation probabilities}. Our method can be roughly summarized in the following \textbf{two-step-procedure}:

\begin{enumerate}
	\item \textbf{Center manifold reduction and normal form of the dynamical system \eqref{eq:rODE}}
	
	For each realization of our random vector of input parameters, we can still use the classical theory of center manifold reduction (see e.g. \cite[Chapter 18]{Wiggins.2003}). This leaves us with a lower-dimensional system that locally still captures the main qualitative behavior of the original system. Of course, the resulting center manifold is in general a random object, and it generically induces a reduced RODE. For this RODE, one can then apply the theory of normal forms (see e.g. \cite[Chapter 3]{Kuznetsov.2004}), which allows us to bring the system into a canonical form, from which we can deduce the type of bifurcation present in the system. The calculation of normal forms has also been successfully used in a wide variety of contexts, see e.g.~\cite{Murdock1,Nayfeh3,DeVilleetal,Touze.2006}. In our case, the coefficients of the normal from are no longer deterministic but consist generically of \emph{nonlinear combinations of the random input parameters}.
	
	\item \textbf{Analysis of the reduced dynamics}
	
	For the analysis of the reduced system, we focus on the probability distributions of the certain key normal form coefficients. As these appear in terms of possibly nonlinear combinations of the random input parameters, oftentimes, they can no longer be identified as one of the standard probability distributions, where closed form analytical expressions for the probability distribution are known. However, the knowledge of this distribution would allow to determine the probability of the occurrence of different types of bifurcations.
	To take into account problem-specific challenges, we come up with three approaches: we analytically derive the bifurcation probability in Section \ref{sec:anaAppUncertaintyProp}. In Section \ref{sec:semiAnaAppUncertaintyProp}, we present a semi-analytical approach in case that a purely analytically solution cannot be found. A brief treatment of a particular sampling-based approach to estimate the bifurcation probability is presented in Section \ref{sec:samplingBasedAppUncertaintyProp}.
\end{enumerate}

We emphasize that our aim is to shed light on the problem of calculating the bifurcation probability type from various points of view and provide solution approaches under various different conditions. We do not claim to present a universally usable approach but rather provide a toolbox of variants for the second step of above two-step procedure that have their problem-specific advantages and drawbacks. Before addressing the main challenge in Sections \ref{sec:anaAppUncertaintyProp}-\ref{sec:samplingBasedAppUncertaintyProp}, we illustrate our two-step procedure in case of the Lorenz system~\cite{Lorenz} with uncertain parameters.

\subsection{Example of Lorenz system}

One of the tools in our methodology is center manifold reduction (see e.g.\ \cite[Chapter 18]{Wiggins.2003} performed in Step 1. We illustrate it here for the Lorenz equations, which are a mode-reduced model for fluid convection (see also \cite[Sec.\ 2.3]{Guckenheimer.1983}) and can be stated as
\begin{align}
	\left\{\begin{array}{ll}
		\dx &= \zeta(y-x) \\
		\dy &= (\rho+1)x - y -xz \\
		\dz &= xy - \theta z.
	\end{array}\right. . \label{eq:LorenzSyst}
\end{align}
Assume $\zeta$ and $\theta$ are $P$-almost surely ($P$-a.s.) positive and independent random variables on a probability space $(\Omega,\AAA,P)$, which possess corresponding probability density functions. Let $\rho$ be the bifurcation parameter, which remains deterministic to ensure being on the given bifurcation curve, and $(x,y,z)^\top \in \R^3$. The center manifold reduction of this problem or determining particular equilibria is standard and quite straightforward but we present it here for completeness.

For $\rho=0$, the eigenvalues of the Jacobian matrix of the linearized flow near the origin are $\lambda_1=0$, $\lambda_2=-(1+\zeta)$, $\lambda_3=-\theta$. By the center manifold theorem \cite[Thm. 18.1.2]{Wiggins.2003}, the system has a one-dimensional center manifold passing through the origin for each realization of $\zeta$ and $\theta$. Thus, we perform center manifold reduction for each realization of $\omega$.

\bigskip
\textbf{Step 1: Center manifold reduction and normal form of \eqref{eq:LorenzSyst}}

Consider the extended system
\begin{align*}
	\left\{\begin{array}{ll}
		\drho &=0 \\
		\begin{pmatrix}
			\dx \\ \dy \\ \dz
		\end{pmatrix} &= \underbrace{\begin{pmatrix}
				-\zeta & \zeta & 0 \\
				1 & -1 & 0 \\
				0 & 0 &-\theta
		\end{pmatrix}}_{=A} \begin{pmatrix}
			x \\ y \\ z
		\end{pmatrix} + \begin{pmatrix}
			0 \\ x(\rho-z) \\ xy
		\end{pmatrix}
	\end{array}\right. .
\end{align*}
The matrix $A$ has the eigenvalues and eigenvectors:
\begin{align*}
	\lambda_1&=0, \quad \lambda_2=-(1+\zeta), \quad \lambda_3=-\theta, \\
	e_1 &= \begin{pmatrix}
		1 \\ 1 \\ 0
	\end{pmatrix}, \quad
	e_2 = \begin{pmatrix}
		\zeta \\ -1 \\ 0
	\end{pmatrix}, \quad
	e_3 = \begin{pmatrix}
		0 \\ 0 \\ 1
	\end{pmatrix}.
\end{align*}
We perform the coordinate transformation $(x,y,z)^\top = T(u,v,w)^\top$, where
\begin{align*}
	T = (e_1 | e_2 | e_3) &= \begin{pmatrix}
		1 & \zeta & 0 \\
		1 & -1 & 0 \\
		0 & 0 & 1
	\end{pmatrix}
\end{align*}
Therefrom, we deduce
\begin{align*}
	\left\{\begin{array}{ll}
		x &= u+\zeta v \\
		y &= u-v \\
		z &= w.
	\end{array}\right., &\quad \begin{pmatrix}
		u \\ v \\ w
	\end{pmatrix} = T^{-1} \begin{pmatrix}
		x \\ y \\ z
	\end{pmatrix} \Rightarrow \begin{pmatrix}
		\du \\ \dv \\ \dw
	\end{pmatrix} = T^{-1} \begin{pmatrix}
		\dx \\ \dy \\ \dz
	\end{pmatrix}.
\end{align*}
Thus, we have
\begin{align}
	\begin{pmatrix}
		\du \\ \dv \\ \dw
	\end{pmatrix} &= \begin{pmatrix}
		\frac{\zeta}{1+\zeta}(\rho u - uw + \rho \zeta v - \zeta v w) \\ -(1+\zeta)v - \frac{1}{1+\zeta}(\rho u - uw + \rho \zeta v - \zeta v w) \\ -\theta w + (u^2 - uv + \zeta v u - \zeta v^2)
	\end{pmatrix}. \label{eq:transformedDynSyst}
\end{align}
The center manifold is given by $v=h_1(u,\rho)$ and $w=h_2(u,\rho)$ with the quadratic approximations
\begin{align}
	\left\{\begin{array}{ll}
		h_1(u,\rho) &= a_1u^2 + a_2u\rho + a_3\rho^2 + \OO(3)\\
		h_2(\rho,u) &= b_1u^2 + b_2u\rho + b_3\rho^2 + \OO(3).
	\end{array}\right. \label{eq:CMquadApprox}
\end{align}
By the invariance of the center manifold, we obtain:
\begin{align}
	\left\{\begin{array}{ll}
		\dv &= 2a_1u\du + a_2\rho\du + \OO(3)\\
		\dw &= 2b_1u\du + b_2\rho\du + \OO(3).
	\end{array}\right. \label{eq:CMderivative}
\end{align}
By substituting \eqref{eq:transformedDynSyst} and \eqref{eq:CMquadApprox} into \eqref{eq:CMderivative} and matching the exponents, we get
\begin{align*}
	a_1 &= 0, \quad a_2 = -\frac{1}{(1+\zeta)^2}, \quad a_3 = 0, \\
	b_1 &= \frac{1}{\theta}, \quad b_2 = 0, \quad b_3 = 0.
\end{align*}

Therefore, the graph of the center manifold reads as
\begin{align}
	\left\{\begin{array}{ll}
		v &= -\frac{1}{(1+\zeta)^2}\rho u + \OO(3)\\
		w &= \frac{1}{\theta}u^2 + \OO(3).
	\end{array}\right\} \quad \text{for} \ 0\leq |\rho|<<1. \label{eq:CMgraph}
\end{align}
By substituting \eqref{eq:CMgraph} in the first equation of \eqref{eq:transformedDynSyst}, we obtain the reduced equations on the CM as
\begin{align*}
	\du &= \frac{\zeta}{1+\zeta}\cdot u \cdot (\rho - \frac{1}{\theta}u^2)(1-\frac{\zeta\rho}{(1+\zeta)^2})
\end{align*}
For $\rho<0$, there is only one fixed point $u=0$ whereas for $\rho>0$, there are three fixed points $u=0$, and $u=\pm\sqrt{\rho\theta}$.

\bigskip
From now on, we consider the reduced equation for the Lorenz system \eqref{eq:LorenzSyst} on the center manifold, which reads as
\begin{align}
	\du &= \frac{\zeta}{1+\zeta} \rho u - \underbrace{\frac{\zeta}{\theta(1+\zeta)}}_{= X}u^3 -\frac{\zeta^2}{(1+\zeta)^3} \rho^2 u + \OO(4). \label{eq:redLorenzSyst}
\end{align}

\textbf{Step 2: Probabilistic analysis of the reduced system \eqref{eq:redLorenzSyst}}

The goal is now to derive the probability distribution of the random variable $X$ representing the \textit{normal form coefficient} of the reduced system. For example, in \eqref{eq:redLorenzSyst}, the sign of $X$ is decisive to determine whether a sub- or supercritical pitchfork bifurcation occurs. For details on sub- and supercritical pitchfork bifurcations, we refer the reader to \cite[Sec.\ 3.4]{Strogatz.2015}.

In general, these probability distributions are complicated to obtain as we often have to deal with nonlinear functions of the random input parameters that appear in the normal form of the reduced system. In case of the reduced Lorenz system \eqref{eq:redLorenzSyst}, the nonlinear transformation function $g:\R^2\rightarrow \R$ for $r_1=\zeta$ and $r_2=\theta$ is
\begin{align*}
	g(r_1,r_2) = \frac{r_1}{r_2(1+r_1)}.
\end{align*}
Our aim is to derive the probability distribution of $X=g(r_1,\dots,r_d)$ for general nonlinear functions $g:\R^d\rightarrow\R$ of the uncertain input parameters $r_1,\dots,r_d$ to determine probabilities for the precise type of bifurcation present in a given RODE of form \eqref{eq:rODE}. 
\begin{rem}
	From now on, we assume that the uncertain input parameters $r_1,\dots,r_d$ are absolutely continuous random variables, i.e.\ possess probability density functions $\rho_{r_1},\dots,\rho_{r_d}$, are independent of each other, and have finite moments. 
\end{rem}

\begin{figure}[h]
    \centering
	\tikzstyle{arrow} = [thick,->,>=stealth, line width=3pt]
	\begin{tikzpicture}
	%
	\node[draw=Apricot!90, fill=Apricot!90, rectangle, rounded corners, inner sep=12pt] (coeff) at (0,0) {Probabilistic bifurcation coefficient};
	%
	%
	\node[draw=Apricot!20, fill=Apricot!20, diamond, text width=2cm, align=center, inner sep=0pt] (Mellin) at (0,-2.5) {Analytical Mellin transform available?};
	\draw [arrow] ([yshift=0.25cm]coeff.south) -- ([yshift=-0.75cm]Mellin.north);
	\node[draw=NavyBlue!20, fill=NavyBlue!20, rectangle, rounded corners, inner sep=12pt] (ana) at (7,-2.5) {Use analytical approach from \textbf{Section \ref{sec:anaAppUncertaintyProp}}};
	\draw [arrow] (Mellin.east) -- (ana);
	\node at ([xshift=0.75cm,yshift=0.25cm]Mellin.east) {yes};
	\node at ([xshift=0.5cm,yshift=0cm]Mellin.south) {no};
	
	\node[draw=Apricot!20, fill=Apricot!20, diamond, text width=2.5cm, align=center, inner sep=0pt] (decomp) at (0,-6.25) {Proficient decomposition possible?};
	\draw [arrow] ([yshift=0.5cm]Mellin.south) -- ([yshift=-0.75cm]decomp.north);
	\node[draw=NavyBlue!20, fill=NavyBlue!20, rectangle, rounded corners, inner sep=12pt] (semiAna) at (7.4,-6.25) {Use semi-analytical approach from \textbf{Section \ref{sec:semiAnaAppUncertaintyProp}}};
	\draw [arrow] (decomp.east) -- (semiAna);
	\node at ([xshift=0.75cm,yshift=0.25cm]decomp.east) {yes};
	\node at ([xshift=0.5cm,yshift=0.2cm]decomp.south) {no};
	
	\node[draw=NavyBlue!20, fill=NavyBlue!20, rectangle, rounded corners, inner sep=12pt] (sampling) at (0,-8.75) {Use sampling-based procedure from \textbf{Section \ref{sec:samplingBasedAppUncertaintyProp}}};
	\draw [arrow] ([yshift=1cm]decomp.south) -- ([yshift=-0.25cm]sampling.north);
		%
	\end{tikzpicture}
    \caption{Guidance when to use which probabilistic method in Step 2.}
    \label{fig:bigPicture}
\end{figure}
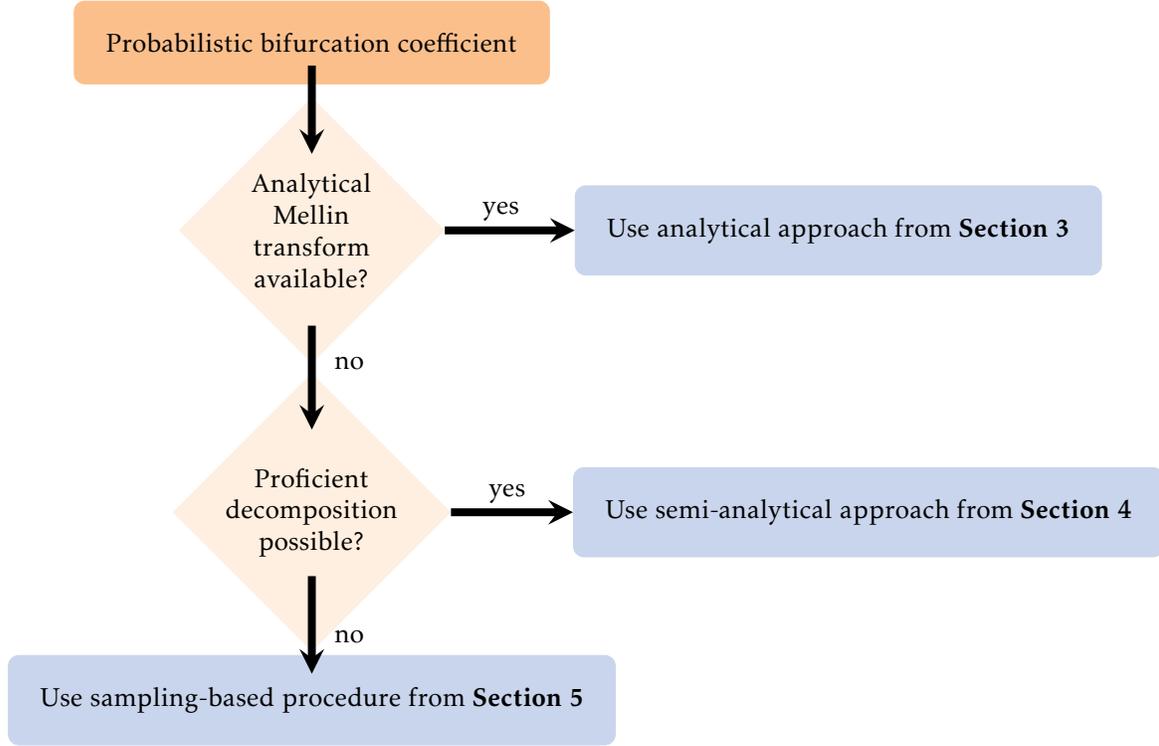

We come up with three different solution approaches in Section \ref{sec:anaAppUncertaintyProp}-\ref{sec:samplingBasedAppUncertaintyProp}. For a given RODE system, the question arises which method is most suited. Step 1 provides the corresponding bifurcation coefficient based on which the method in Step 2 is chosen. A guidance when each of the three probabilistic approaches is suitable is sketched in Figure \ref{fig:bigPicture}.

\bigskip
Before presenting the analytical and the semi-analytical approach in Sections \ref{sec:anaAppUncertaintyProp} and \ref{sec:semiAnaAppUncertaintyProp}, we introduce the \textit{Mellin transform} as a tool used therein for the probabilistic analysis of the bifurcation coefficient.
\subsection{Mellin transform} \label{subsec:MellinTrafo}
As emphasized in \cite{Springer.1979}, the Mellin transform is of crucial importance in the derivation of probability density functions of products, quotients and algebraic functions of independent random variables. It belongs to the class of integral transforms and is defined as follows.

\begin{defin}{(cf.\ \cite[p.\ 96]{Springer.1979})}
	The \textit{Mellin transform} is defined as
	\begin{align}
		\MM\left(f(x)\right)(s) &= \int_{0}^{\infty} x^{s-1}f(x) dx \label{eq:MellinTransform}
	\end{align}
	for a single- and real-valued function $f$ that is defined almost everywhere for $x\geq 0$.
\end{defin}

As usual in the theory of integral transforms, there is an inverse transformation linked to it.

\begin{defin}{(cf.\ [\cite[p.\ 96]{Springer.1979})}
	The \textit{inverse Mellin transform} is defined as
	\begin{align}
		f(x) &= \MM^{-1}\left(\MM\left(f(x)\right)(s)\right)(x) = \frac{1}{2\pi i}\int_{c-i\infty}^{c+i\infty} x^{-s} \MM\left(f(x)\right)(s) ~\textnormal{d}s \label{eq:invMellinTransform}
	\end{align}
	if the existence of the Mellin transform is ensured and the Mellin transform is an analytic function of the complex variable $s$ for $c_1 \leq \Re(s)=c \leq c_2$ for real $c_1$ and $c_2$.
\end{defin}

Under these conditions, the function $f$ is uniquely determined by its Mellin transform (see e.g.\ \cite{Springer.1979}). In the context of probability theory, this means that there is a one-to-one correspondence between probability density functions and their Mellin transforms (see \cite{Epstein.1948}).

As pointed out in \cite{Epstein.1948}, there is an obvious probabilistic interpretation of the Mellin transform:
\begin{align}
	\MM\left(\rho_\xi(x)\right)(s) = \E\left[\xi^{s-1}\right] =: \MM(\xi)(s), \label{eq:MellinProbabilistic}
\end{align}
where $\xi$ is an almost surely (a.s.) \textbf{positive} random variable on a probability space $\left(\Omega,\AAA,P\right)$ with corresponding probability density $\rho_\xi$. An extension of the Mellin transform to random variables also attaining negative values is described in \cite{Epstein.1948}.
For the remainder of this work, we shall define the Mellin transform of a random variable $\xi$ as the one of the corresponding probability density function $\rho_\xi$.
It has already been recognized in \cite{Rajan.2015} that the Mellin transform is a powerful tool in uncertainty evaluation. Therein, the authors analytically calculate the standard uncertainty, defined as the square-root of the second-order central moment of a given measurement distribution, by means of Mellin transforms. Note however that their input-output relationship is restricted to multivariate polynomials.
Also in our work, the interpretation \eqref{eq:MellinProbabilistic} of the Mellin transform $\MM(\xi)(s)$ as the $(s-1)$-th moment of the (a.s.) positive random variable $\xi$ will be of crucial importance for the semi-analytical approach to the probabilities of the presence of bifurcation types in Section \ref{sec:semiAnaAppUncertaintyProp}.

\bigskip

The Mellin transform has many operational properties (see e.g. \cite[Sec.\ 6.2.1]{Andrews.1999}). We recall here from \cite{Epstein.1948} the most appealing ones for our purpose.

\begin{prop} \label{prop:propertiesMellin}
	Let $\xi$ be a positive random variable on $(\Omega,\AAA,P)$ and the function $\rho_\xi$ its corresponding probability density function. Then, the following properties hold true:
	\begin{itemize}
		\item The positive random variable $\eta = a\xi$, where $a>0$, has the Mellin transform
		\begin{align}
			\MM(\eta)(s) &= a^{s-1}\MM(\xi)(s). \label{eq:scalePropMellin}
		\end{align}
		\item The Mellin transform of the positive random variable $\eta = \xi^a$ reads as
		\begin{align}
			\MM(\eta)(s) &=\MM(\xi)(as-a+1). \label{eq:powerRVMellin}
		\end{align}
		As an immediate consequence, we obtain the Mellin transform of the inverse $\eta = \frac{1}{\xi}$ of the positive random variable $\xi$ as
		\begin{align}
			\MM(\eta)(s) &=\MM(\xi)(-s+2). \label{eq:invRVMellin}
		\end{align}
		\item For positive independent random variables $\xi_1,\dots,\xi_n$, we can calculate the Mellin transform of the product $\eta=\prod_{i=1}^{n}\xi_i$ as
		\begin{align}
			\MM(\eta)(s) &= \prod_{i=1}^{n}\MM(\xi_i)(s). \label{eq:prodMellin}
		\end{align}
	\end{itemize}
\end{prop}

\begin{rem}
The product density of two non-negative independent random variables having probability density functions (PDFs) $\rho_{r_1}(x)$ and $\rho_{r_2}(x)$, i.e.\ the inversion of the Mellin transform in \eqref{eq:prodMellin}, can be expressed as a so called \textit{Mellin convolution} (see e.g.\ \cite[Ch.\ 4]{Springer.1979}), which reads as
\begin{align}
\rho_{r_1r_2}(y) &= \int_{0}^{\infty} \frac{1}{x}\rho_{r_1}\left(\frac{y}{x}\right)\rho_{r_2}(x) ~\textnormal{d}x. \label{eq:MellinConv}
\end{align}
\end{rem}

\section{Analytical approach for uncertainty propagation} 
\label{sec:anaAppUncertaintyProp}

In our first approach, we want to tackle the calculation of the probability distribution of the nonlinear transformation $X=g(r_1,\dots,r_d)$ of the uncertain input parameters analytically. For simple sums of random variables, we could employ the usual convolution formula for the resulting density, yet in many contexts, we cannot expect just sums but we also have to consider products. The Mellin transform, introduced in Subsection \ref{subsec:MellinTrafo}, provides an efficient tool for handling products of independent random variables. Hence, we are interested in analyzing the Mellin transform of the nonlinear transformation $X=g(r_1,\dots,r_d)$.

\subsection{Mellin transform of the bifurcation coefficient}

For some specific RODEs of type \eqref{eq:rODE}, the Mellin transform of the bifurcation coefficient and its inverse can be calculated analytically. Thereby, we directly obtain a closed-form analytical expression for the probability density function of the bifurcation coefficient (the random variable $X$ in case of the reduced Lorenz system \eqref{eq:redLorenzSyst}) and can analytically calculate the probability of the bifurcation type. We now illustrate the benefit of our procedure with an example.

\begin{ex} 
\label{ex:HuangManchesterSheet8}
To illustrate the main idea, we ask what is the probability to face a subcritical pitchfork bifurcation in the following reduced RODE system
\begin{align}
	\du &= -r_1r_2u^3 - au - a^2u + \OO(4). \label{eq:CMexHuangManchester}
\end{align}
Thereby, $a$ is the deterministic bifurcation parameter and $r_1,r_2 \in L^\infty(\Omega,\AAA,P)$ are random variables on a probability space $(\Omega,\AAA,P)$ with corresponding probability density functions $\rho_{r_1}(x)$ and $\rho_{r_2}(x)$.

We answer this question by applying our two-step procedure, introduced in Section \ref{sec:BifurcationAndInstab}. The dynamics \eqref{eq:CMexHuangManchester} are already in reduced normal form of a pitchfork bifurcation so that we can directly address Step 2 of our methodology.
Note that if $-r_1r_2(\omega)>0$, we have a subcritical pitchfork bifurcation in \eqref{eq:CMexHuangManchester} and if $-r_1r_2(\omega)<0$, we deal with a supercritical pitchfork bifurcation. Hence, the second step now consists in the probabilistic analysis of the reduced dynamics \eqref{eq:CMexHuangManchester}. To determine the probability of the presence of a sub- or supercritical bifurcation, we calculate the Mellin transform of $r_1r_2$. Note that not only positive random variables might appear in the bifurcation coefficient and, a priori, the Mellin transform is only defined for positive random variables. However, an extension to real-valued random variables is possible. In \cite{Epstein.1948}, the author makes use of the Mellin convolution \eqref{eq:MellinConv} to extend the Mellin transform procedure to the product of independent random variables which take both positive and negative values. The idea is to decompose the PDFs into their positive and negative parts and to consider each arising pair separately.

Assume now that $r_1\sim \UU(-1,3)$, i.e. $r_1$ has the probability density function
$$\rho_{r_1}(x) = \frac{1}{4}\cdot\ind_{[-1,3]}(x).$$
Assume further that $r_2 \sim \Gamdist(3,1)$, where $3$ is the shape parameter and $1$ denotes the rate, i.e.\ the probability density function of $r_2$ reads as
$$\rho_{r_2}(x) = \frac{x^{2}e^{-x}}{\Gamma(3)}\cdot\ind_{(0,\infty)}(x).$$
Moreover, $r_1$ and $r_2$ are assumed to be independent. As in \cite{Epstein.1948}, we decompose the PDF of $r_1$ as $\rho_{r_1}(x) = \rho_{\text{pos}}(x) + \rho_{\text{neg}}(x)$, where
\begin{align*}
	\rho_{\text{pos}}(x) &= \frac{1}{4}\cdot\ind_{[0,3]}(x), \quad \text{and} \quad \rho_{\text{neg}}(x) = \frac{1}{4}\cdot\ind_{[-1,0)}(x).
\end{align*}
Note that the negative part of $\rho_{r_2}(x)$ is zero and the positive part is the PDF itself. Hence, we consider the two pairs $P_1 = [\rho_{\text{pos}}(x), \rho_{r_2}(x)]$ and $P_2 = [\rho_{\text{neg}}(x),\rho_{r_2}(x)]$. For $P_1$, we can apply the Mellin convolution \eqref{eq:MellinConv} directly since both functions are positive.
\begin{align*}
	h_1(y) &= \int_{0}^{\infty} \frac{1}{x}\rho_{\text{pos}}\left(\frac{y}{x}\right)\rho_{r_2}(x) dx \\
	&= \int_{0}^{\infty} \left( \frac{1}{x}\frac{1}{4}\cdot\ind_{[0,3]}\left(\frac{y}{x}\right)\frac{x^{2}e^{-x}}{\Gamma(3)}\cdot\ind_{(0,\infty)}(x)\right) dx \\
	&= \frac{1}{4\Gamma(3)} \int_{\nicefrac{y}{3}}^{\infty} xe^{-x} dx \cdot \ind_{(0,\infty)}(y) \\
	&= \frac{1}{4\Gamma(3)} \left(1 + \frac{y}{3}\right)e^{-\nicefrac{y}{3}} \cdot \ind_{(0,\infty)}(y)
\end{align*}
For the second pair $P_2$, we first consider an auxiliary function.
\begin{align*}
	h_2^\ast(y) &= \int_{0}^{\infty} \frac{1}{x}\rho_{\text{neg}}\left(-\frac{y}{x}\right)\rho_{r_2}(x) dx \\
	&= \int_{0}^{\infty} \frac{1}{x}\frac{1}{4}\cdot\ind_{[-1,0)}\left(-\frac{y}{x}\right)\frac{x^2e^{-x}}{\Gamma(3)}\cdot\ind_{(0,\infty)}(x) dx \\
	&= \frac{1}{4\Gamma(3)}\int_{y}^{\infty} \left(xe^{-x}\right)dx \cdot \ind_{(0,\infty)}(y) \\
	&= \frac{1}{4\Gamma(3)}\left(1+y\right)e^{-y}\cdot\ind_{(0,\infty)}(y)
\end{align*}
By setting
\begin{align*}
h_2(y) &= h_2^\ast(-y) = \frac{1}{4\Gamma(3)}\left(1-y\right)e^{y}\cdot\ind_{(-\infty,0)}(y),
\end{align*}
we finally obtain the PDF of the product $r_1r_2$ as
\begin{align}
	\rho_{r_1r_2}(x) &= h_1(x) + h_2(x) = \frac{1}{4\Gamma(3)} \left(1 + \frac{x}{3}\right)e^{-\nicefrac{x}{3}} \cdot \ind_{(0,\infty)}(x) + \frac{1}{4\Gamma(3)}\left(1-x\right)e^{x}\cdot\ind_{(-\infty,0)}(x). \label{eq:PDFr1r2}
\end{align}

Based on the PDF \eqref{eq:PDFr1r2}, we can now calculate the cumulative distribution function (CDF) $\Phi_{r_1r_2}(x)$ of the product of $r_1$ and $r_2$ analytically and obtain
\begin{align}
	\Phi_{r_1r_2}(x) &= \int_{-\infty}^{x} \rho_{r_1r_2}(z) dz = \frac{1}{4\Gamma(3)} \left(2-x\right)e^{x} \cdot \ind_{(-\infty,0)}(x) + \frac{1}{4\Gamma(3)}\cdot \left( \left(-6-x\right)e^{-\nicefrac{x}{3}} +8\right) \cdot \ind_{[0,\infty)}(x). \label{eq:CDFr1r2}
\end{align}
With the CDF \eqref{eq:CDFr1r2}, we have a complete probabilistic characterization of the bifurcation coefficient $r_1r_2$. In Figure \ref{fig:probabDist_r1r2}, the PDF and the CDF are depicted together with the corresponding MC approximations based on $10^6$ samples. For both functions, we observe a very close match indicating the correctness of the transformation and inversion procedure.

\begin{figure}[h]
	\centering
	\subfloat[\label{fig:PDFr1r2}\ PDF of $r_1r_2$]{\begin{overpic}[width=0.45\textwidth]{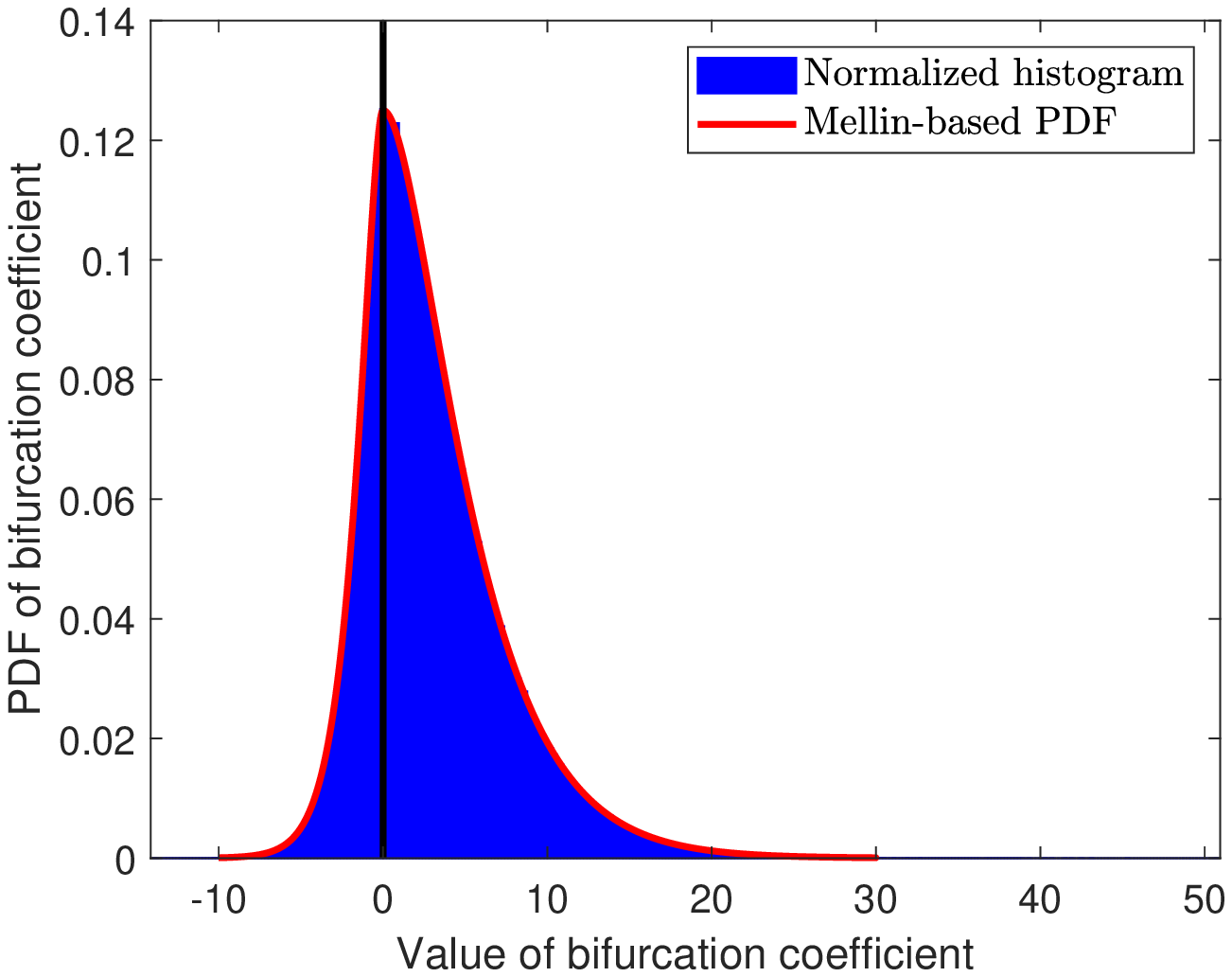}\end{overpic}}\hspace{0.2cm}
	\subfloat[\label{fig:CDFr1r2}\ CDF of $r_1r_2$]{\begin{overpic}[width=0.45\textwidth]{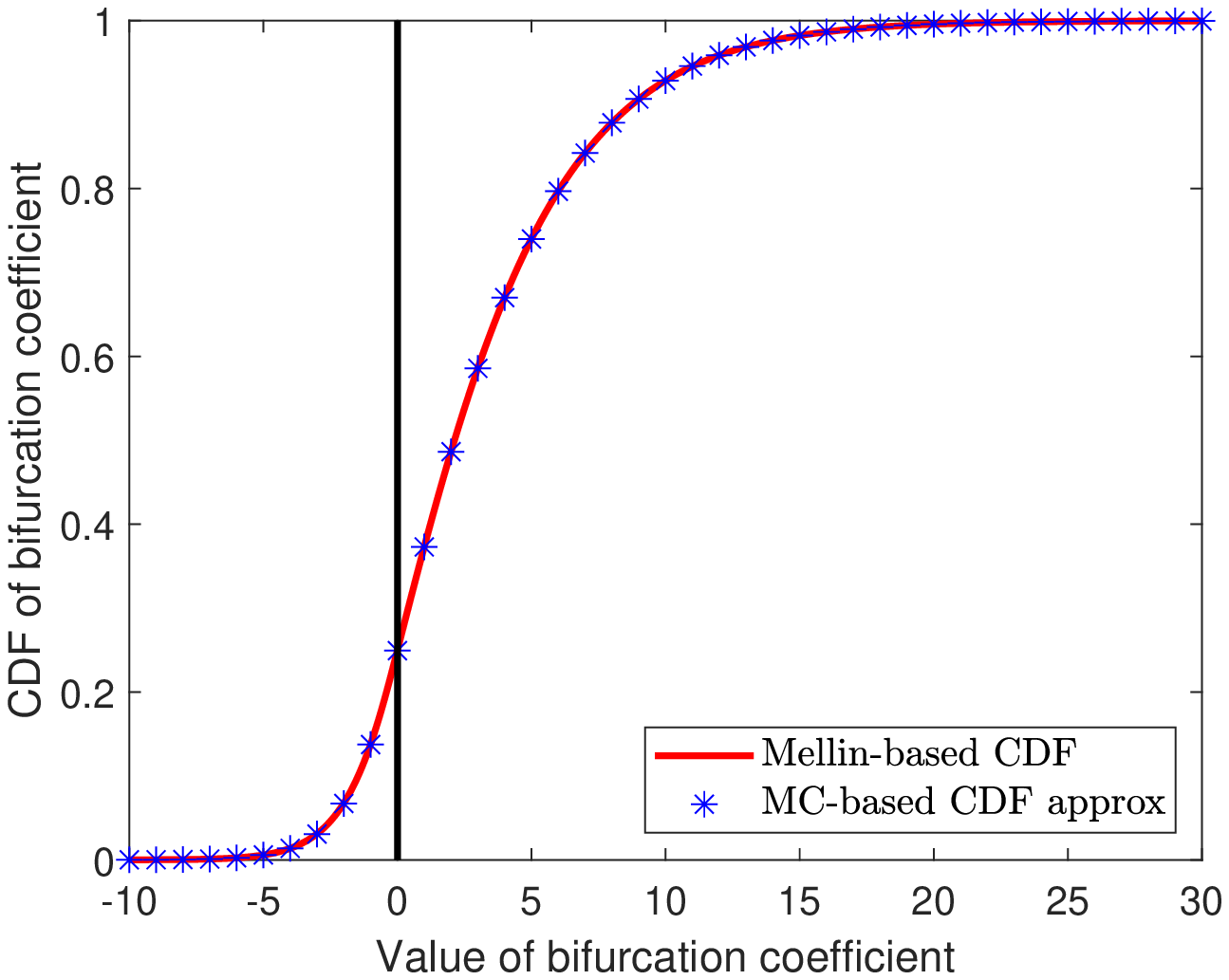}\end{overpic}}
	
	\caption{Analytical PDF and CDF of the bifurcation coefficient $r_1r_2$ of RODE \eqref{eq:CMexHuangManchester} based on the Mellin convolution; $r_1 \sim \UU(-1,3)$ and $r_2 \sim \Gamdist(3,1)$; corresponding MC approximations based on $M=10^6$ samples depicted in blue. Close agreement indicates correctness of analytical Mellin procedure.}
	
	\label{fig:probabDist_r1r2}
\end{figure}

As already mentioned, the sign of the bifurcation coefficient $r_1r_2$ is decisive to infer the type of bifurcation. Therefore, we calculate the bifurcation probability using \eqref{eq:CDFr1r2}, and obtain that the probability of the pitchfork bifurcation to be subcritical is
\begin{align*}
	P(-r_1r_2>0) = P(r_1r_2\leq 0) = \Phi_{r_1r_2}(0) = 0.25.
\end{align*}
The bifurcation probability of $0.25$ is validated by a MC simulation of $10^6$ samples, which lead to a value of $0.24936$. Accordingly, the probability of the pitchfork bifurcation to be supercritical amounts to $0.75$.
\end{ex}

\subsection{Quantification of the influence of small perturbations}

Do we obtain similar bifurcation probabilities in the presence of a small perturbation or do the probabilities differ drastically? Via our analytical solution approach, we can quantify how a perturbation $\varepsilon>0$ influences the Mellin transform and its inversion, and therefore how the PDF of the bifurcation normal form coefficient changes. We illustrate this in the setting stated below in Example \ref{ex:GammaGammawithGammaPerturb}.

\begin{ex} 
\label{ex:GammaGammawithGammaPerturb}
We start with an additive perturbation of a bifurcation coefficient $X=r_1r_2$. Assume that $r_1\sim\Gamdist(\alpha,1)$, $r_2\sim\Gamdist(\beta,1)$, $b_{\varepsilon}\sim\Gamdist(\varepsilon,1)$, where $\varepsilon>0$ and $r_1,r_2$, and $b$ are independent of each other.
	The Mellin transform of a random variable $\xi\sim\Gamdist(\alpha,1)$ with shape parameter $\alpha>0$ and  rate parameter $1$ is given by (see \cite[Table 1]{Epstein.1948})
	\begin{align}
		\MM\left(\xi\right)(s) &= \frac{\Gamma(\alpha-1+s)}{\Gamma(\alpha)}, \quad \text{for} \ \Re(s)>-\alpha+1, \label{eq:MellinGamma}
	\end{align}
	where $\Gamma(\cdot)$ denotes the Gamma function.
	Hence, without perturbation by $b_\varepsilon$, we obtain the probability density function $\rho_{r_1r_2}(x)$ via Mellin transformation and subsequent inversion as follows.
	\begin{align}
		\MM(r_1r_2)(s) &= \MM(r_1)(s)\MM(r_2)(s) \notag\\
		&= \frac{\Gamma(\alpha-1+s)}{\Gamma(\alpha)}\cdot \frac{\Gamma(\beta-1+s)}{\Gamma(\beta)}, \label{eq:MellinGammaGamma}
	\end{align}
	Note that we used property \eqref{eq:prodMellin} of the Mellin transform for the product of independent random variables.
	The probability density function $\rho_{r_1r_2}$ is then obtained by inverting the Mellin transform \eqref{eq:MellinGammaGamma} as
	\begin{align}
		\rho_{r_1r_2}(x) &= \MM^{-1}\left(\frac{\Gamma(\alpha-1+s)}{\Gamma(\alpha)}\cdot \frac{\Gamma(\beta-1+s)}{\Gamma(\beta)}\right)(x) = \frac{2x^{\nicefrac{1}{2}(-2+\alpha+\beta)}K_{-\alpha+\beta}\left(2\sqrt{x}\right)}{\Gamma(\alpha)\Gamma(\beta)}, \label{eq:invMellinGammaGamma}
	\end{align}
	where $K_n(x)$ is the modified Bessel function of the second kind and the inversion was obtained by using the software \textit{Mathematica} \cite{Mathematica.2019}; see \cite{Malik.1968} for an analytical proof of this fact.
	
	If we now consider an additive gamma perturbation $b_{\varepsilon}\sim\Gamdist(\varepsilon,1)$ with mean $\varepsilon$ of the random variable $r_2$, we can again perform the Mellin transformation of $X=r_1(r_2+b_{\varepsilon})$ and its inversion to obtain the probability density function $\rho_{r_1(r_2+b_{\varepsilon})}$.
	\begin{align*}
		\MM(r_1(r_2+b_{\varepsilon}))(s) &= \MM(r_1)(s)\MM(r_2+b_{\varepsilon})(s) \\
		&= \MM(r_1)(s)\MM(\tilde{r_2})(s), 
	\end{align*}
	where $\tilde{r_2}\sim\Gamdist(\beta+\varepsilon,1)$ (see e.g.\ \cite{Dufresne.2010}).
	The perturbed Mellin transform thus reads as
	\begin{align}
		\MM(r_1(r_2+b_{\varepsilon}))(s) &= \frac{\Gamma(\alpha-1+s)}{\Gamma(\alpha)}\cdot\frac{\Gamma(\beta+\varepsilon-1+s)}{\Gamma(\beta+\varepsilon)}. \label{eq:MellinGammaGammaperturbed}
	\end{align}
	By again using \textit{Mathematica} \cite{Mathematica.2019}, we obtain the perturbed probability density function as
	\begin{align}
		\rho_{r_1(r_2+b_{\varepsilon})}(x) &=\MM^{-1}\left(\frac{\Gamma(\alpha-1+s)}{\Gamma(\alpha)}\cdot\frac{\Gamma(\beta+\varepsilon-1+s)}{\Gamma(\beta+\varepsilon)}\right)(x) = \frac{2x^{\nicefrac{1}{2}(-2+\alpha+\beta+\varepsilon)}K_{-\alpha+\beta+\varepsilon}\left(2\sqrt{x}\right)}{\Gamma(\alpha)\Gamma(\beta+\varepsilon)} \label{eq:PDFGammaperturbed}
	\end{align}
	with corresponding series expansion
	\begin{align*}
		\rho_{r_1(r_2+b_{\varepsilon})}(x) =& \rho_{r_1r_2}(x) + \varepsilon\cdot\left[\frac{x^{\nicefrac{1}{2}(-2+\alpha+\beta)}}{\Gamma(\alpha)\Gamma(\beta)} \cdot \left(K_{-\alpha+\beta}\left(2\sqrt{x}\right)\left(\log(x)-2\Psi_0(\beta)\right) + 2K_{-\alpha+\beta}^{(1,0)}\left(2\sqrt{x}\right) \right)\right] + \OO(\varepsilon^2),
	\end{align*}
	where $\Psi_0(x)$ is the digamma function and $K_n^{(1,0)}(x)$ is the first derivative with respect to $n$ of the modified Bessel function of the second kind. This reveals how the perturbation influences the PDF of the bifurcation coefficient.
\end{ex}

For arbitrary RODEs of type \eqref{eq:rODE}, such explicit analytical calculations of the Mellin transform, its inverse, and a perturbed version are not always possible as one often only obtains implicit formulas, particularly for complicated combined sums and products of random variables. A possible remedy consists of the extension of the above analytical approach by suitable approximations of random variables and a numerical estimation procedure.

\section{Semi-analytical approach for uncertainty propagation} \label{sec:semiAnaAppUncertaintyProp}

In this section, we extend the analytical approach of Section \ref{sec:anaAppUncertaintyProp} to address nonlinear combinations of input parameters where a direct Mellin transformation and inversion are not readily performed. To do so, we take the bifurcation coefficient $X=g(r_1,r_2) = \nicefrac{r_1}{r_2(1+r_1)}$ of the reduced Lorenz system \eqref{eq:redLorenzSyst} as inspiration. By using property \eqref{eq:invRVMellin}, the Mellin transform of the factor $\nicefrac{1}{r_2}$ is readily calculated. It remains to calculate the Mellin transform of $\nicefrac{r_1}{1+r_1}$.
We will illustrate the working principle of our semi-analytical approach by answering the following questions for the reduced Lorenz system \eqref{eq:redLorenzSyst} in Example \ref{ex:MellinTransformPCEpartLorenzRV}:
\begin{itemize}
	\item What is the Mellin transform of the part containing $r_1$?
	\item What is the Mellin transform of the bifurcation coefficient?
	\item What is the probability for the pitchfork bifurcation to be subcritical?
\end{itemize}
Our semi-analytical approach is not limited to the particular RODE \eqref{eq:LorenzSyst}. We do not use the particular type of nonlinearity nor do we use special characteristics of the system.

\subsection{Combining polynomial chaos expansion with the Mellin transformation}

A nice way to handle measurable functions of random variables is to approximate them by using a generalized \textit{polynomial chaos expansion (PCE)} (see e.g.\ \cite{Ghanem.1991}). As emphasized in \cite[Sec.\ 2.4]{Muehlpfordt.2020}, the PCE of a random variable can be compared to the Fourier series of a periodic signal: both allow to trace back an infinite-dimensional mathematical construct to an (in)finite number of coefficients. Under suitable assumptions, our nonlinear transformation $X=g(r)$ can be represented by the PCE
\begin{align}
	g &= \sum_{n \in \N_0^d} g_n\Pol_n, \label{eq:PCEnonlinTrafo_complete}
\end{align}
where $\left(\Pol_n\right)_{n\in \N^d}$ represents an appropriate orthogonal polynomial basis and $g_n$ are the corresponding coefficients. The generalized PCE is also used in \cite{Scott.2011} to model the propagation of a combination of uncertainties, amongst others in model parameters. It is well-known that mean and variance of $X$ can be easily obtained based on the coefficients $g_n$ in \eqref{eq:PCEnonlinTrafo_complete} (see e.g. \cite[Sec.\ 11.3]{Sullivan.2015}). There are many other advantages of using a PCE (see e.g.\ \cite{Ghanem.2017}). Aiming at using as much analytical information as possible, we come up with a combined use of the Mellin transformation method and the PCE. Thus, our aim is to simplify the Mellin transformation as far as possible by making use of known results for standard distributions (see e.g. \cite[Table 1]{Epstein.1948} or \cite[Table D.2]{Springer.1979}) and the properties stated in Proposition \ref{prop:propertiesMellin}. The remaining combinations of input parameters that could not be separated are then taken into account via PCE.

In the example of the reduced Lorenz system \eqref{eq:redLorenzSyst}, this would mean that we first calculate the Mellin transform of $\nicefrac{1}{r_2}$ by using \eqref{eq:invRVMellin}. Then, we use a PCE for $\nicefrac{r_1}{1+r_1}$ and calculate the Mellin transform thereof. 

\subsubsection{Mellin transform of a polynomial chaos expansion} \label{subsubsec:MellinPCE}

As we aim at separating the individual uncertain input parameters via Mellin transformation, from now on, we assume to work with one-dimensional PCEs of the form
\begin{align}
	\g(r_i) &= \sum_{n=0}^{\infty} \g_n\Pol_n(\xi), \label{eq:PCEnonlinTrafo_individual}
\end{align}
where $\xi$ is the \textit{stochastic germ}, i.e.\ a random variable with probability distribution corresponding to the chosen polynomial basis. The latter choice of the polynomial basis depends on the probability distribution of $r_i$ that suggests a suitable stochastic germ e.g.\ in terms of a common support. There exists a unique family of orthogonal polynomials $\left(\Pol_n\right)_{n\in \N}$ such that for given positive normalization constants $\left(h_n\right)_{n\in \N}$, the orthogonality condition with respect to the weight $\rho_{\xi}$ of the stochastic germ, i.e.
\begin{align*}
	<\Pol_{n_1},\Pol_{n_2}> &:= \int_{\R} \Pol_{n_1}(y)\Pol_{n_2}(y)\rho_{\xi}(y) dy = h_{n_1}\delta_{n_1n_2}, \quad \forall \ n_1,n_2 \in \N,
\end{align*}
holds true (see e.g.\ \cite{Kuehn.2020}). For many standard probability distributions, the corresponding family of orthogonal polynomials $\left(\Pol_n\right)_{n\in \N}$ is listed in the Wiener--Askey scheme (see e.g.\ \cite[Table 4.1]{Xiu.2002}).

To include the PCE \eqref{eq:PCEnonlinTrafo_individual} into our numerical estimation routine, we work with a truncated version of the latter, i.e.
\begin{align}
\g(r_i) &\approx \sum_{n=0}^{N} \g_n\Pol_n(\xi), \label{eq:truncPCEnonlinTrafo_individual}
\end{align}
where $N$ is an appropriately chosen truncation index. For a discussion on arising truncation errors, we refer the reader to \cite{Blatman.2009,Sullivan.2015,Xiu.2002}. We already emphasize that this will not be our focus here and a very low truncation index will be enough to reproduce the Mellin transforms reasonably well as we will see later.

\bigskip
Let us now assume that we fixed an orthogonal polynomial basis $\left(\Pol_n\right)_{n\in \N}$ with corresponding stochastic germ $\xi$  and a truncation index $N$. The idea is to use the Mellin transform of the stochastic germ $\xi$, which is known explicitly in many cases (see e.g. \cite[Table 1]{Epstein.1948} or \cite[Table D.2]{Springer.1979}). The arising powers of $\xi$ can easily be handled by using property \eqref{eq:powerRVMellin} of the Mellin transform. The difficulty consists in handling the sums of the powers of $\xi$, which are not independent of each other, and in adding up the polynomials of $\xi$ themselves. Note that for two a.s. positive absolutely continuous random variables $X_1$ and $X_2$, the Mellin transform of their sum $\MM(X_1+X_2)(s)$ does not equal the sum of the corresponding Mellin transforms $\MM(X_1)(s)$ and $\MM(X_2)(s)$ in general.

To circumvent this problem, in the calculation of the Mellin transform of the PCE \eqref{eq:truncPCEnonlinTrafo_individual}, we make use of an iterated application of the binomial formula. Note that the Mellin transform of the sum of random variables and the use of a binomial formula has also been considered in \cite{Dufresne.2010}. We are then able to trace back the Mellin transform of the PCE \eqref{eq:truncPCEnonlinTrafo_individual} to the individual Mellin transform of $\xi$. 

\bigskip
In a first step, we address the issue that the stochastic germ might not be an a.s.\ positive random variable. For example, a common choice is a uniform random variable $U \in \UU(-1,1)$ with corresponding familiy of orthogonal polynomials given by the Legendre polynomials. Note that different normalizations are used in the literature\footnote{Be aware of the different normalization of the \textit{MATLAB} function \textit{legendreP}; see also \url{https://de.mathworks.com/help/symbolic/legendrep.html#buei9e5-6}, last checked: November, 3, 2020.}. As we will use the MATLAB-based software framework \textit{UQLab} \cite{Sudret.2014} later within our numerical estimation procedure, we use the definition specified in \cite{Sudret.2019}, which corresponds to
\begin{align*}
	<\legPol_{n_1},\legPol_{n_2}> &= \int_{-1}^{1}\legPol_{n_1}(y)\legPol_{n_2}(y)\cdot\frac{1}{2} = h_{n_1}\delta_{n_1n_2}, \quad \text{with}\ h_{n_1}=(2n_1+1)\cdot\frac{1}{2}\cdot \frac{2}{2n_1+1}=1, \quad \forall \ n_1,n_2 \in \N.
\end{align*}
To work with an a.s.\ positive uniform random variable $\tilde{U} \in \UU(0,1)$, we use the shifted Legendre polynomials, that is $\tilde{\legPol}_n(x) = \legPol_n(2x-1)$.

\bigskip

In a second step, we collect the powers of $\tilde{U}$. This leads to the representation of the PCE \eqref{eq:truncPCEnonlinTrafo_individual} for a uniform stochastic germ $\xi = U \sim \UU(-1,1)$ with the Legendre polynomials $\left(\legPol_n\right)_{n\in \N}$ being the corresponding orthogonal polynomials as
\begin{align*}
	\g(r_i) &\approx \sum_{n=0}^{N} \g_n\legPol_n(U) = \sum_{n=0}^{N} \g_n\tilde{\legPol}_n(\tilde{U}) = \sum_{n=0}^{N} \coeff_n \tilde{U}^n,
\end{align*}
where $\left(\coeff_n\right)_{n\in \N}$ are the coefficients calculated by collecting the powers of $\tilde{U}$ arising in $\tilde{\legPol}_n(\tilde{U})$ for all $n\in \{0,\dots,N\}$. 
Note that a representation in terms of collected powers can also be achieved for other stochastic germs $\xi$ such as a Beta distributed one as well. We thus work with the reformulation of \eqref{eq:truncPCEnonlinTrafo_individual}, which reads as
\begin{align}
\g(r_i) &\approx \sum_{n=0}^{N} \coeff_n \xi^n, \label{eq:truncPCE_collectedPowers}
\end{align}
where $\left(\coeff_n\right)_{n\in \N}$ are the coefficients calculated based on the collected powers of $\xi$ in the corresponding familiy of orthogonal polynomials $\left(\Pol_n\right)_{n\in \N}$.
One of the advantages of representation \eqref{eq:truncPCE_collectedPowers} is that the Mellin transform of $\g(r_i)$ can be calculated in the same way for an a.s.\ positive random variable $\tilde{r_i}$ and a random variable $r_i$ which also takes negative values. The Mellin transform of the reformulated PCE \eqref{eq:truncPCE_collectedPowers} finally reads as
\begin{align}
	\MM\left(\sum_{n=0}^{N} \coeff_n \xi^n\right)(s) &= \sum_{i=0}^{N\cdot(s-1)} \hat{c}_i(s) \cdot \MM(\xi)(i+1), \label{eq:MellinTransform_PCE}
\end{align}
where $\hat{c}_i(s)$ denote the cumulated coefficients for the powers of $\xi$ which arise due to a repeated use of the binomial formula and depend on the chosen evaluation $s$ of the Mellin transform. Details of our calculation can be found in Appendix \ref{app:MellinPCE}.

Note that we limit our calculation to integer values of $s\geq 1$.
An alternative representation of the Mellin transform in \eqref{eq:MellinTransform_PCE} for non-integer values of $s$ by means of complex integrals can be found in \cite[Thm.\ 3.1]{Dufresne.2010}. Therein, the author provides an explicit formula for the Mellin transform of sums of random variables involving multiple integrations along complex paths. It might be that summing residues enables an even more concrete representation of the Mellin transform of the PCE \eqref{eq:MellinTransform_PCE}. This more concrete representation might also pave the way for a perturbation analysis of the Mellin transform of the PCE \eqref{eq:truncPCEnonlinTrafo_individual}.

\bigskip

However, by again taking a closer look at the probabilistic interpretation of the Mellin transform \eqref{eq:MellinProbabilistic}, we recognize that the Mellin transform evaluated at integer values of $s$ corresponds to the moments of the random variable. Hence, even with this restriction, we still get valuable information about the probability distribution of the nonlinear transformation of the uncertain input in terms of an analytical formula for an arbitrary high number of its moments, provided they exist. Note that moments of polynomial expressions can also be computed based on Mellin transforms in the automated framework for uncertainty evaluation from \cite{Rajan.2015}.
For now, we content ourselves with the representation \eqref{eq:MellinTransform_PCE} and continue our analysis of the probability distribution of the bifurcation coefficient.

\bigskip

The stochastic germ $\xi$ usually has a standard distribution, which means that its Mellin transform is known in analytical form in many cases (see e.g.\ \cite[Table 1]{Epstein.1948}). The great advantage is that \eqref{eq:MellinTransform_PCE} can be further simplified once the distribution of $\xi$ is specified. We consider two common choices of stochastic germs: uniform random variables and beta random variables.

\bigskip

\textbf{Expansion based on uniform random variables}

The Mellin transform of a uniformly distributed stochastic germ $\xi_U\sim\UU(0,1)$ reads as (see e.g.\ \cite[Table 1]{Epstein.1948})
\begin{align}
	\MM\left(\xi_U\right)(s) &= \frac{1}{s}, \quad \text{for} \ \Re(s)>0. \label{eq:MellinUniform01}
\end{align}
The expression \eqref{eq:MellinUniform01} can now be plugged into \eqref{eq:MellinTransform_PCE} to obtain
\begin{align}
\MM\left(\sum_{n=0}^{N} \coeff_n \xi_U^n\right)(s) &= \sum_{i=0}^{N\cdot(s-1)} \frac{\hat{c}_i(s)}{i+1}. \label{eq:MellinTransform_PCE_Uniform}
\end{align}

\bigskip

\textbf{Expansion based on beta random variables}

Similarly, for a beta distributed stochastic germ $\xi_B\sim\text{Beta}(\alpha,\beta)$, the Mellin transform is also known in analytical closed form. According to \cite[Table 1]{Epstein.1948}, the Mellin transform of a random variable $X\sim\text{Beta}(\alpha,\beta)$ with $\alpha,\beta>0$ is 
\begin{align}
	\MM\left(X\right)(s) &= \frac{\Gamma(\alpha+\beta)\cdot \Gamma(\alpha-1+s)}{\Gamma(\alpha)\cdot\Gamma(\alpha+\beta-1+s)}, \quad \text{for} \ \Re(s)>-\alpha+1. \label{eq:MellinBeta}
\end{align}
We plug the Mellin transform \eqref{eq:MellinBeta} into \eqref{eq:MellinTransform_PCE} and end up with the Mellin transform of the PCE
\begin{align}
	\MM\left(\sum_{n=0}^{N} \coeff_n \xi_B^n\right)(s) &= \sum_{i=0}^{N\cdot(s-1)} \hat{c}_i(s) \frac{\Gamma(\alpha+\beta)\cdot \Gamma(\alpha+i)}{\Gamma(\alpha)\cdot\Gamma(\alpha+\beta+i)}. \label{eq:MellinTransform_PCE_Beta}
\end{align}

We come back to the task of calculating the probability distribution of the normal form coefficient in the reduced Lorenz system \eqref{eq:redLorenzSyst}. Remember that the problem revealed at the end of Section \ref{sec:anaAppUncertaintyProp} consisted in calculating the Mellin transform of $\nicefrac{r_1}{1+r_1}$. Approximating the random variable $\g(r_1) = \nicefrac{r_1}{1+r_1}$ first by a PCE and then calculating the Mellin transform thereof solves the problem.

\begin{ex} 
\label{ex:MellinTransformPCEpartLorenzRV} We are interested in the Mellin transform of the part containing $r_1$ in the normal form coefficient of the reduced Lorenz system \eqref{eq:redLorenzSyst}. To illustrate our procedure of calculating the Mellin transform of a PCE more concretely, we perform the steps for $\g(r_1) = \nicefrac{r_1}{1+r_1}$ with $r_1 \sim \Beta(2,2)$, a uniform stochastic germ $\xi = \tilde{U}\sim\UU(0,1)$ and a chosen truncation of $N=2$. We use the MATLAB-based software framework \textit{UQLab} \cite{Sudret.2014}, to obtain the PCE \eqref{eq:truncPCEnonlinTrafo_individual} (rounded to 4 digits) for $\g(r_1)$ as
	\begin{align}
		\g(r_1) &\approx 0.3188 + 0.1002 \cdot \tilde{P}_1(\tilde{U}) -0.0130 \cdot \tilde{P}_2(\tilde{U}). \label{eq:truncPCEnonlinTrafo_Beta22RV}
	\end{align}
	The reformulation \eqref{eq:truncPCE_collectedPowers} is obtained by collecting the powers in \eqref{eq:truncPCEnonlinTrafo_Beta22RV} and reads as
	\begin{align}
		\tilde{g}(r_1) &\approx 0.1163 + 0.5210 \cdot \tilde{U} -0.1739 \cdot \tilde{U}^2. \label{eq:truncPCEnonlinTrafo_Beta22RV_collected}
	\end{align}
	In our notation in formula \eqref{eq:truncPCE_collectedPowers}, this corresponds to
	\begin{align*}
		\coeff_0 &= 0.1163, \quad \coeff_1 = 0.5210, \quad \text{and} \ \coeff_2 = -0.1739.
	\end{align*}
	To get a better uncerstanding of the $\hat{c}_i(s)$ in \eqref{eq:MellinTransform_PCE_Uniform}, we note that, here, the repeated use of the binomial formula to obtain the Mellin transform of \eqref{eq:truncPCEnonlinTrafo_Beta22RV_collected} (see the appendix for the general derivation) leads to
	\begin{align}
		\MM(\coeff_0 + \coeff_1 \cdot \tilde{U} + \coeff_2 \cdot \tilde{U}^2)(s) &= \int_{0}^{1} \left(\sum_{k=0}^{s-1} \binom{s-1}{k}\left(\sum_{l=0}^{s-1-k}\binom{s-1-k}{l} \coeff_0^{s-1-k-l} \coeff_1^l x^l \right) \coeff_2^k x^{2k} \right)~\textnormal{d}x.
	\end{align}
	The coefficients $\hat{c}_i(s)$ result from the collection of all powers of the integration variable $x$ in the integrand. The approximated Mellin transform of $\g(r_1)$ is then obtained by plugging the $\hat{c}_i(s)$ into \eqref{eq:MellinTransform_PCE_Uniform}.
	Table \ref{tab:coefficients_hat_c_i} shows explicitly how to calculate the coefficients $\hat{c}_i(s)$ based on the coefficients $\coeff_n$ for $s\in\{2,3,4\}$ here in our example.
	\begin{table}[htbp]
		\centering
		\begin{tabular}{c|c|c|c|c|c|c|c}
			& $\hat{c}_0(s)$ & $\hat{c}_1(s)$ & $\hat{c}_2(s)$ & $\hat{c}_3(s)$ & $\hat{c}_4(s)$ & $\hat{c}_5(s)$ & $\hat{c}_6(s)$ \\
			\hline
			$s=2$ & $\coeff_0$ & $\coeff_1$ & $\coeff_2$ & - & - & - &  \\
			\hline
			$s=3$ & $\coeff_0^2$ & $2\coeff_0\coeff_1$ & $2\coeff_0\coeff_2 + \coeff_1^2$ & $2\coeff_1\coeff_2$ & $\coeff_2^2$ & - & - \\
			\hline
			$s=4$ & $\coeff_0^3$ & $3\coeff_0^2\coeff_1$ & $3\coeff_0\coeff_1^2 + 3\coeff_0^2\coeff_2$ & $\coeff_1^3 + 6\coeff_0\coeff_1\coeff_2$ & $3\coeff_1^2\coeff_2 + 3\coeff_0\coeff_2^2$ & $3\coeff_1\coeff_2^2$ & $\coeff_2^3$ \\
		\end{tabular}
		\caption{Relation between $\coeff_n$ and $\hat{c}_i(s)$ for truncation $N=2$ and $s\in\{2,3,4\}$}
		\label{tab:coefficients_hat_c_i}
	\end{table}
\end{ex}
We now have a tool at hand that allows us to deal with parts that cannot be further decomposed by using the properties of the Mellin transform stated in Proposition \ref{prop:propertiesMellin}. We first approximate these parts of the bifurcation coefficient via a PCE and then calculate the Mellin transform of the PCE by using \eqref{eq:MellinTransform_PCE}.

\subsubsection{Mellin transform of bifurcation coefficient}
We now come back to the combination of the analytical Mellin transform (wherever possible) with the Mellin transform of the PCE approximation, derived in Subsection \ref{subsubsec:MellinPCE}. It remains to build the Mellin transform of the normal form coefficient based on the individual components.
As announced at the beginning of the section, we will test our semi-analytical approach for the bifurcation coefficient of the reduced Lorenz system \eqref{eq:redLorenzSyst} and continue our Example \ref{ex:MellinTransformPCEpartLorenzRV}.

\addtocounter{ex}{-1}
\begin{ex}[continued] \label{ex:bifurcationCoeff_redLorenzSyst}
Assume that $r_1 \sim \Beta(2,2)$ and $r_2 \sim \Gamdist(8,1)$, where $8$ is the shape parameter and $1$ denotes the rate. Let further $r_1$ and $r_2$ be independent as required. In Example \ref{ex:MellinTransformPCEpartLorenzRV}, we already derived  the Mellin transform \eqref{eq:MellinTransform_PCE_Uniform} of $\nicefrac{r_1}{1+r_1}$. By using the property \eqref{eq:invRVMellin} of the Mellin transform and the known Mellin transform \eqref{eq:MellinGamma} of a Gamma random variable, we obtain
\begin{align}
	\MM\left(\frac{1}{r_2}\right)(s) &= \MM\left(r_2\right)(-s+2) = \frac{\Gamma(8-s+1)}{\Gamma(8)}. \label{eq:MellinTransformInvPartLorenzRV}
\end{align}
Hence, to obtain the Mellin transform of the normal form coefficient, we combine these results to obtain
\begin{align}
	\MM\left(\nicefrac{r_1}{r_2(1+r_1)}\right)(s) &= \frac{\Gamma(8-s+1)}{\Gamma(8)}\cdot \sum_{i=0}^{2\cdot(s-1)} \frac{\hat{c}_i(s)}{i+1}. \label{eq:MellinLorenzRV}
\end{align}
The coefficients $\hat{c}_i(s)$ are explicitly stated for $s \in \{2,3,4\}$ in Table \ref{tab:coefficients_hat_c_i}.
\end{ex}

The Mellin inversion of \eqref{eq:MellinLorenzRV} would provide us the desired probability density function of the bifurcation coefficient. Though, it is unclear whether the latter can be obtained analytically, at least not directly based on \eqref{eq:MellinLorenzRV} as our calculation was restricted to integer values of $s$. Note however that by means of \eqref{eq:MellinLorenzRV}, we can state \emph{all} moments of the normal form coefficient (provided they exist) in terms of an analytical expression.

\subsection{From Mellin transforms to the probabilities of bifurcation types} \label{subsec:probBifurcationTypes}

We can now profit from the probabilistic interpretation \eqref{eq:MellinProbabilistic} of the Mellin transform as moments of the corresponding probability density function. The reconstruction of the PDF from the series of moments of a distribution with compact support is known under the name \textit{Hausdorff moment problem} (see e.g.\ \cite{Mnatsakanov.2008b,Mnatsakanov.2008a}) and as \textit{Stieltjes moment problem} for a support on the positive half-line (\cite{Stieltjes.1894}). There exist several conditions under which the distribution is uniquely determined by its moments (see e.g.\ \cite{Mnatsakanov.2008b} and references therein).

We emphasize that the Mellin transform \eqref{eq:MellinTransform_PCE} of a general PCE, possibly combined with other analytically available Mellin transforms as in Example \ref{ex:bifurcationCoeff_redLorenzSyst}, makes moments available in closed form without sampling. Therefore, we use a moment-based estimation procedure of the probability distribution of the bifurcation coefficient based on which we can then derive the probabilities of bifurcation types. We analyze the applicability of nonparametric polynomial approximations of the probability density function (see Section \ref{subsubsec:nonparaEst_PolynomialEstPDF}) and a semi-parametric estimation procedure via a method of moments \cite{Hansen.1982} in Section \ref{subsubsec:semiParaEst_MethodOfMoments}, both having their advantages and drawbacks. However, it will turn out that the semi-parametric approach much better suits our purpose of estimating bifurcation type probabilities.

\subsubsection{Nonparametric estimation: polynomial approximation of probability distribution functions} \label{subsubsec:nonparaEst_PolynomialEstPDF}
A nonparametric estimation of the PDF has the advantage that no assumptions on the shape of the distribution are needed and the reconstruction is very flexible. It is often possible to approximate the PDF of the quantity of interest via summation of orthogonal polynomials based on the given moment sequence $(\mu_j)_{j\in \N}$ (cf.\ \cite{Provost.2005}). The number of moments needed to calculate a meaningful approximation depends on the irregularity of the PDF to approximate (cf.\ \cite{Provost.2005}). In many instances where PDFs are unimodal or show a similar behavior to standard univariate PDFs, a few moments are already sufficient (see e.g.\ the case of $\nmoms=4$ in our numerical results in Section \ref{sec:numRes}).
We now present different techniques for the reconstruction of the PDF. The latter contains the probabilistic information about the bifurcation coefficient based on which we can derive the probabilities of bifurcation types.

\bigskip
\textbf{Method based on Legendre polynomials}

According to \cite{Provost.2005}, the $\nmoms$-th degree polynomial approximation of the density function of a random variable defined on the interval $[a,b]$ given its first $\nmoms+1$ moments $\mu_0,\mu_1,\dots,\mu_N$, where $\mu_0=1$, reads as
\begin{align}
\rho_{N, \text{legPol}}(y) &= \sum_{k=0}^{N} \frac{2k+1}{b-a}\left(\sum_{j=0}^{k}c_{j}^{\text{leg},k}\cdot\mu_j\right)\cdot P_k\left(\frac{2y -(a+b)}{b-a}\right), \label{eq:PDF_approx_legPol}
\end{align}
where $c_{j}^{\text{leg},k}$ are the coefficients corresponding to the $k$-th Legendre polynomial $P_k$ with evaluation in $\frac{2y -(a+b)}{b-a}$.

\bigskip
\textbf{Method based on monic orthogonal polynomials}

In \cite{Provost.2015}, given the moments of a random variable, the PDF is approximated via orthogonal polynomials generated based on a normalized weight function, which serves as an initial approximation of the PDF. Therefore, the support of this weight function and the one of the bifurcation coefficient should coincide.

The $N$-th degree orthogonal polynomial approximant of the PDF $\rho_X(x)\ind_{(a,b)}(x)$ with respect to the weight function $w(y)\ind_{(a,b)}(y)$ reads as
\begin{align}
	\rho_{N,monic}(y) = w(y)\left(c_w+\sum_{i=n_p+1}^{N}\lambda_i \pi_i(y)\right)\ind_{(a,b)}(y), \label{eq:PDF_approx_monic}
\end{align}
where the monic orthogonal polynomials are $\pi_i(y) = \sum_{k=0}^{i}d_{i,k}y^k$, and $\lambda_i = \nicefrac{\sum_{k=0}^{i}d_{i,k}\mu_k}{\int_{a}^{b}w(y) \pi_i(y)^2dy}$.

\bigskip
\textbf{Method based on transformed moments}

The density approximant in \cite{Mnatsakanov.2008a} for a given truncation value $N \in \N$ of the moment sequence $\mu_0,\mu_1,\dots,\mu_N$, where $\mu_0=1$, is derived for $a=0$ and $b>0$ as
\begin{align}
\rho_{N, \text{trafoMom}}(y) &= \left(\frac{\Gamma(N+2)}{\Gamma(\lfloor Ny/b\rfloor+1)} \cdot \frac{1}{b^{\lfloor Ny/b\rfloor+1}} \sum_{m=0}^{N-\lfloor Ny/b\rfloor} \frac{(-\nicefrac{1}{b})^m \mu_{m+\lfloor Ny/b\rfloor}}{m!(N-\lfloor Ny/b\rfloor-m)!}\right)\ind_{(0,b)}(y). \label{eq:PDF_approx_trafoMom}
\end{align}
Although choosing one of the above approximations seems natural, we will see in the numerical Section \ref{sec:numRes} that problems with instabilities arise in our setting and the a priori unkown support of the bifurcation coefficient complicates the estimation procedure.

\subsubsection{Semi-parametric estimation: Gaussian mixture models and method of moments} \label{subsubsec:semiParaEst_MethodOfMoments}

Therefore, besides the polynomial approximation techniques introduced above, we also suggest another direction and analyze the applicability of a semi-parametric estimation procedure. We focus on the estimation of the probability distribution by using the method of moments in the context of Gaussian mixture models. However, before explaining our use of the method of moments in the context of Gaussian mixture models in detail, we complete our example of the reduced Lorenz system \eqref{eq:redLorenzSyst} by stating the final estimate of the probability of the bifurcation type present in the system and showing that the method of moments for Gaussian mixture models captures well the probability distribution of the bifurcation coefficient.
	
\begin{figure}[h!]
		\centering
		\subfloat[\label{fig:PDF_genBeta25}\ Uncertain input $r_1 \sim \text{gen}\Beta_{-0.5}^{0.5}(2,5)$]{\begin{overpic}[width=0.48\textwidth]{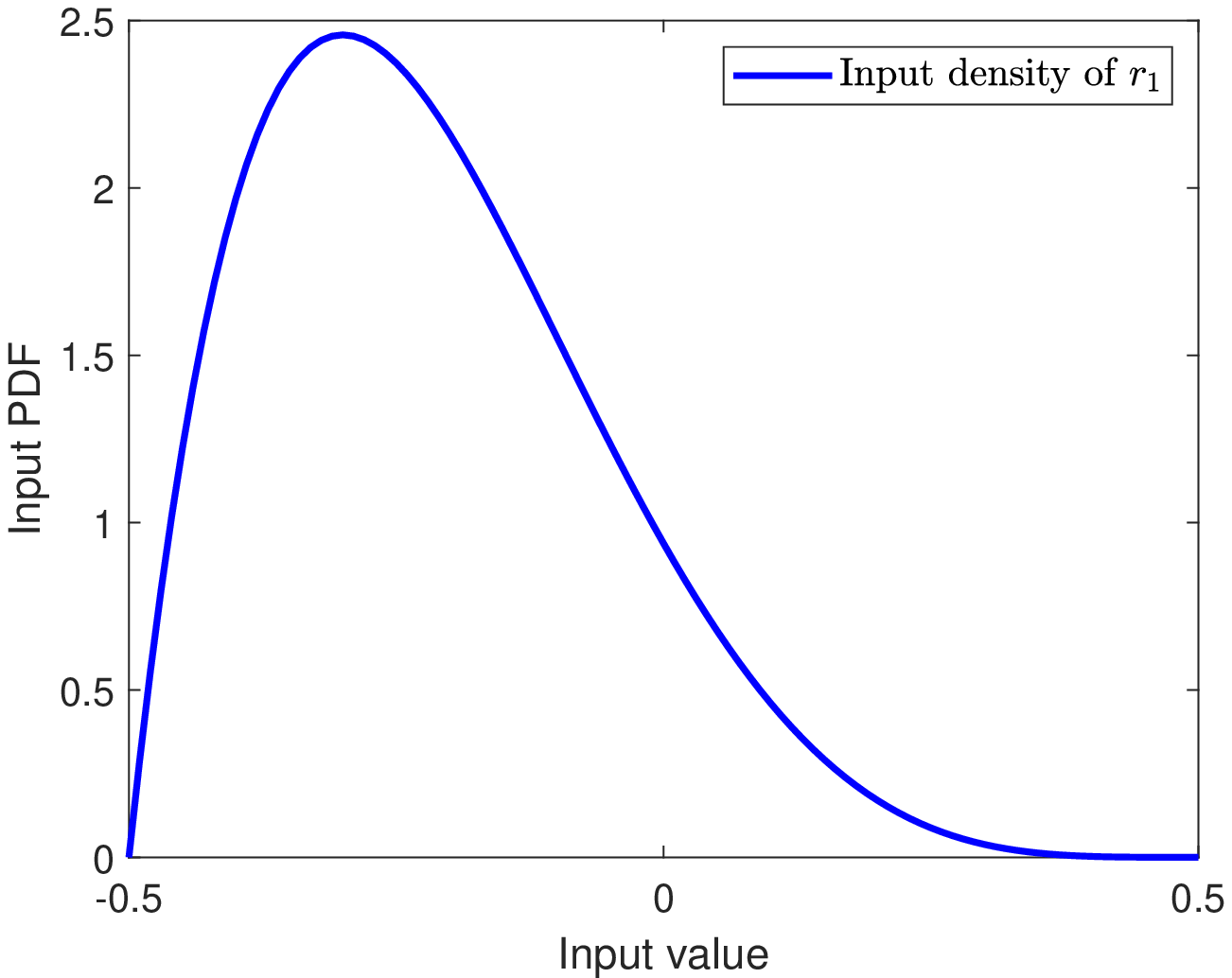}\end{overpic}}
		\subfloat[\label{fig:PDF_gamma81}\ Uncertain input $r_2 \sim \Gamdist(8,1)$]{\begin{overpic}[width=0.48\textwidth]{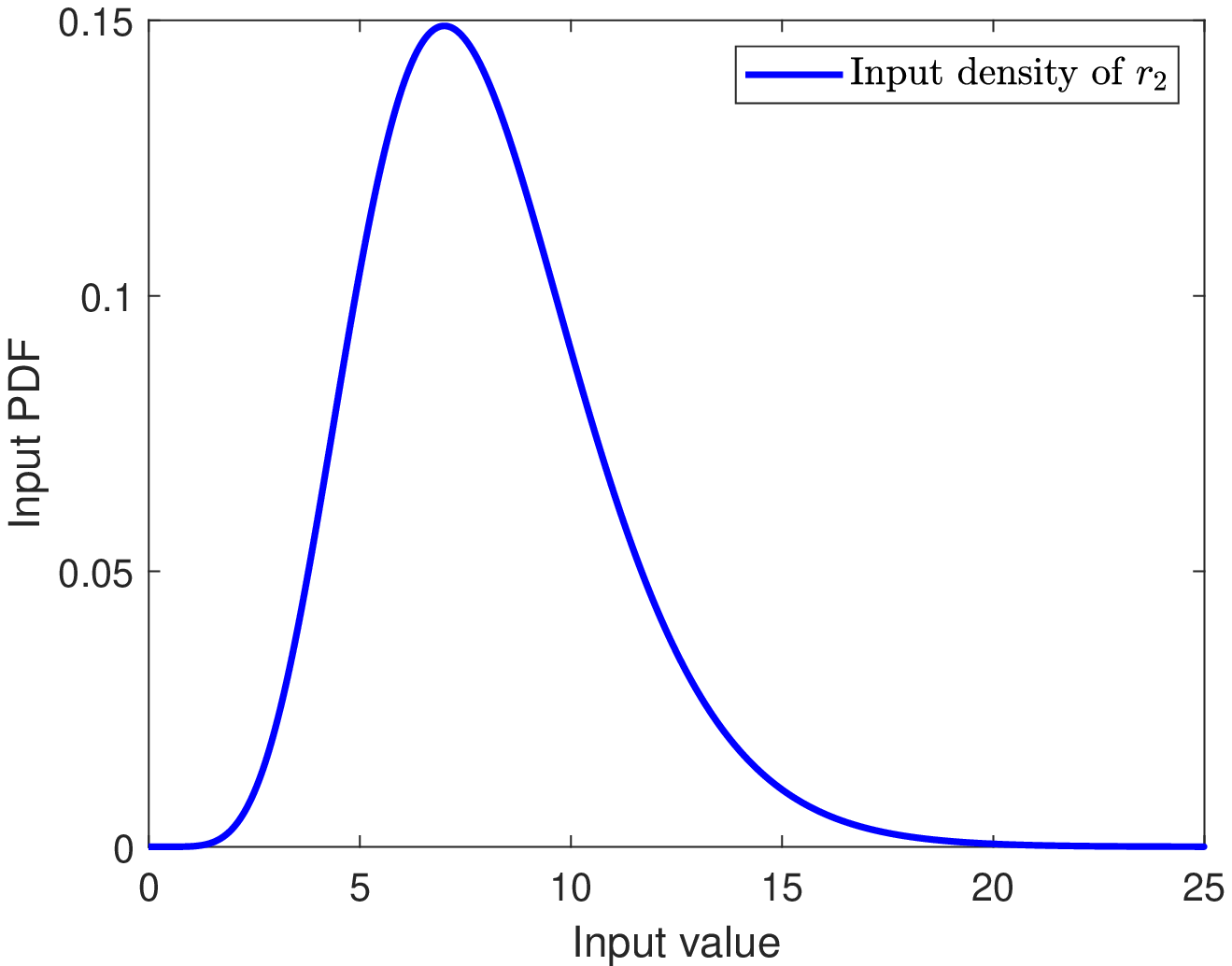}\end{overpic}}
		
		\subfloat[\label{fig:PDF_approx_LorenzRV}\ Approximated PDF of bifurcation coefficient $X$ in \eqref{eq:redLorenzSyst}]{\begin{overpic}[width=0.48\textwidth]{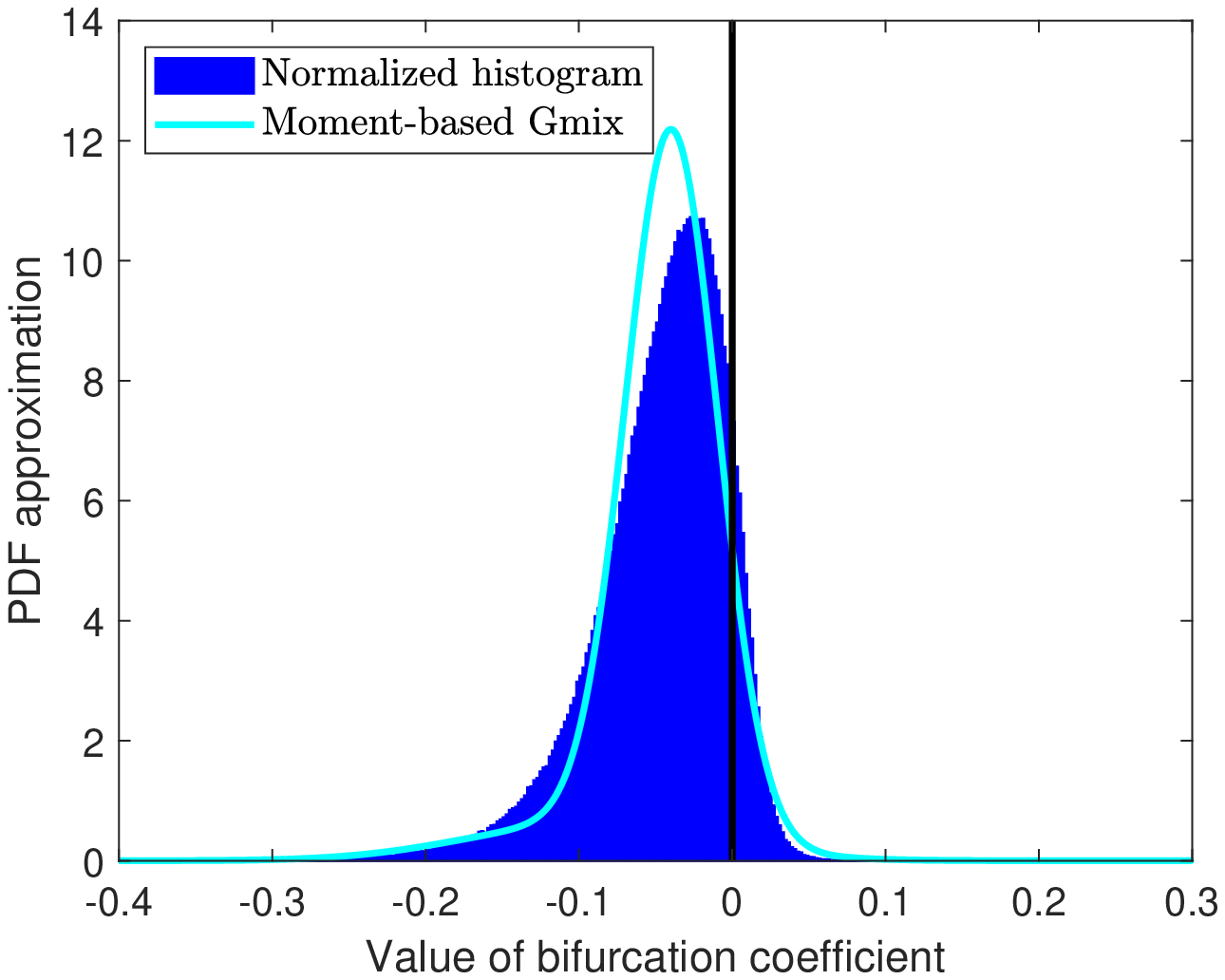}\end{overpic}}

		\caption{PDF approximation of the bifurcation coefficient $X$ in the reduced Lorenz system \eqref{eq:redLorenzSyst}; Upper figures show PDFs of the input parameters $r_1$ and $r_2$; well-designed combination of PCE with Mellin transform gives moments of $X$; PDF of Gaussian mixture (Gmix) model is obtained via a method of moments; good agreement of the PDF with the normalized histogram over $M=10^6$.}
		
		\label{fig:PDF_approx_teaserLorenzSyst}
	\end{figure}

\addtocounter{ex}{-1}
\begin{ex}[continued]
Assume the uncertain input $r_1$ follows a generalized Beta distribution supported on $[-0.5,0.5]$ with parameters $\alpha=2$ and $\beta=5$, i.e.\ $r_1\sim\Beta_{-0.5}^{0.5}(2,5)$, and $r_2$ is Gamma distributed with shape $8$ and rate $1$, i.e. \ $r_2\sim\Gamma(8,1)$. Note that we relax the positivity assumption on the parameter $r_1$ here to illustrate the methodological concept. Figures \ref{fig:PDF_genBeta25} and \ref{fig:PDF_gamma81} show the corresponding input densities. Then, the final estimate for the presence of a subcritical pitchfork bifurcation in the Lorenz system \eqref{eq:LorenzSyst} is $0.9049$. The corresponding estimate obtained based on $M=10^6$ Monte Carlo (MC) samples amounts to $0.8903$, so our estimation is in close agreement. The approximated PDF of the bifurcation coefficient together with the sample-based normalized histogram is depicted in Figure \ref{fig:PDF_approx_LorenzRV}.
\end{ex}

With this example of the approach in mind, we now explain the method from semi-parametric estimation in detail. More specifically, in our approach, we make use of univariate \textit{Gaussian mixture models (GMMs)}.
\begin{defin}[Finite mixture distribution (cf.\ {\cite[p.\ 357]{McLachlan.2019}})]
	A \textit{finite mixture distribution} of a $p$-dimensional [absolutely continuous] random vector $Y$ has the probability density function
	\begin{align}
		\rho_Y(y) = \sum_{i=1}^{k} \pi_i \rho_i(y), \label{eq:mixtModel}
	\end{align}
	where $\pi_i$ denote the \textit{mixing proportions} that are nonnegative and sum up to one. The $\rho_i$ represent the \textit{component densities} and $k$ is the number of components.
\end{defin}

In the case of GMMs, the $\rho_i$ in \eqref{eq:mixtModel} are taken to be normal densities.
For details on finite mixture models, we refer the reader to \cite{Fruehwirth-Schnatter.2006,McLachlan.2000}. GMMs have a long-standing history and estimating their parameters is one of the oldest estimation problems in the statistical community (cf.\ \cite{Everitt.1981}). Mixture models are able to represent unknown distributions arising from random interactions in a computationally convenient way (cf. \cite[p.\ 356]{McLachlan.2019}). As they are very flexible, they are not only used for cluster analysis but also in situations where a single component distribution might not be capable of capturing the underlying heterogeneity (see e.g.\ \cite[p.\ 359]{McLachlan.2019}). Thereby, the Gaussian mixture densities are particularly appealing as they lie dense in the set of all density functions (see \cite[Thm.\ 3.1]{Li.1999}). When it comes to density estimation, Gaussian mixtures take the role of sigmoidal neural
networks in the context of function fitting and approximation (cf. \cite[p.\ 279]{LiBarron.1999}). Here, we use Gaussian mixture models to estimate bifurcation probabilities. It is known that a finite number of moments is sufficient to identify a finite Gaussian mixture model (except for permuting components) (see \cite{Wu.2020}).

\bigskip
To estimate the parameters of the Gaussian mixture model, we will use the \textit{generalized method of moments}, which appeared first in the econometrics literature in \cite{Hansen.1982} (see e.g.\ \cite{Hall.2005}). As the name indicates, the method is a generalization of the classical method of moments, which can be traced back to K.\ Pearson \cite{Pearson.1894} (see \cite[p.\ 5]{Hall.2005} as well as \cite[p.\ 30]{Everitt.1981} and~\cite{Anandkumar.2012,Byrnes.2008,Gandhi.2016,Lindsay.1989,Wu.2020}).

The generalized method of moments allows to take into account more moments than there are parameters to fit. Thus, the system of moment conditions is overdetermined and one minimizes the difference between the estimated and the given moments to match. This is beneficial as it circumvents the problem that, in the classical method of moments, the system might not be solvable meaning that the combination of estimated moments does not describe the moments of any probability distribution (see \cite{Wu.2020}). For example, conditions for a sequence of numbers to be identified as the moments of a mixture distribution are given in terms of conditions on moment matrices in \cite{Lindsay.1989}.



Oftentimes, the \textit{Maximum Likelihood} estimation technique is preferred over the method of moments (see e.g.\ \cite{Lindsay.1989}) or \textit{Bayesian} estimation techniques are in the focus (see e.g.\ \cite[p.\ 41]{Fruehwirth-Schnatter.2006}). 
However, one striking difference of our setting to the classical setting in the literature is that we are given the exact moments or the ones calculated based on the PCE approximation \eqref{eq:MellinTransform_PCE} in closed form instead of samples/data from which empirical moments need to be estimated. Therefore, we do not consider \textit{Maximum Likelihood} or \textit{Bayesian} estimation in this work.

\subsection{Numerical results} \label{sec:numRes}
In this section, we test the applicability of the two types of estimation procedures described in Sections \ref{subsubsec:nonparaEst_PolynomialEstPDF} and \ref{subsubsec:semiParaEst_MethodOfMoments} for estimating the probabilities of bifurcation types such as sub- or supercritical ones determined by the sign of the bifurcation coefficient inherent to the RODE \eqref{eq:rODE}.
In a preparatory step, we carry out numerical experiments to reproduce a known standard probability density function and analyze the advantages and drawbacks of the polynomial approximations in Section \ref{subsubsec:appApproaches_poly}. Thereby, we keep in mind that the overall goal is the estimation of probability distributions of bifurcation coefficients such as $X$ in the reduced Lorenz system \eqref{eq:redLorenzSyst}. Furthermore, in Section \ref{subsubsec:appApproaches_GMMs}, we give some details of our numerical implementation of the method of moments in the context of GMMs.

The main step finally consists in testing the applicability of above methods to estimate the PDF of the bifurcation coefficient.
From the approximated PDF, the probabilities of the bifurcation types can then be derived.

\subsubsection{Applicability of polynomial approximations} \label{subsubsec:appApproaches_poly}

To get an intuition whether the polynomial approximations \eqref{eq:PDF_approx_legPol}, \eqref{eq:PDF_approx_monic}, and \eqref{eq:PDF_approx_trafoMom} might be suitable, we test the PDF approximation methods by trying to recover a Beta distribution based on the knowledge of their first $\nmoms+1$ moments (calculated using the known analytical expressions for the standard distributions and not by means of their Mellin transforms).

\begin{figure}[h!]
	\centering
	\subfloat[\label{fig:PDF_approx_BetaRV_nmoms4}\ PDF $\nmoms=4$]{\includegraphics[width=0.5\textwidth]{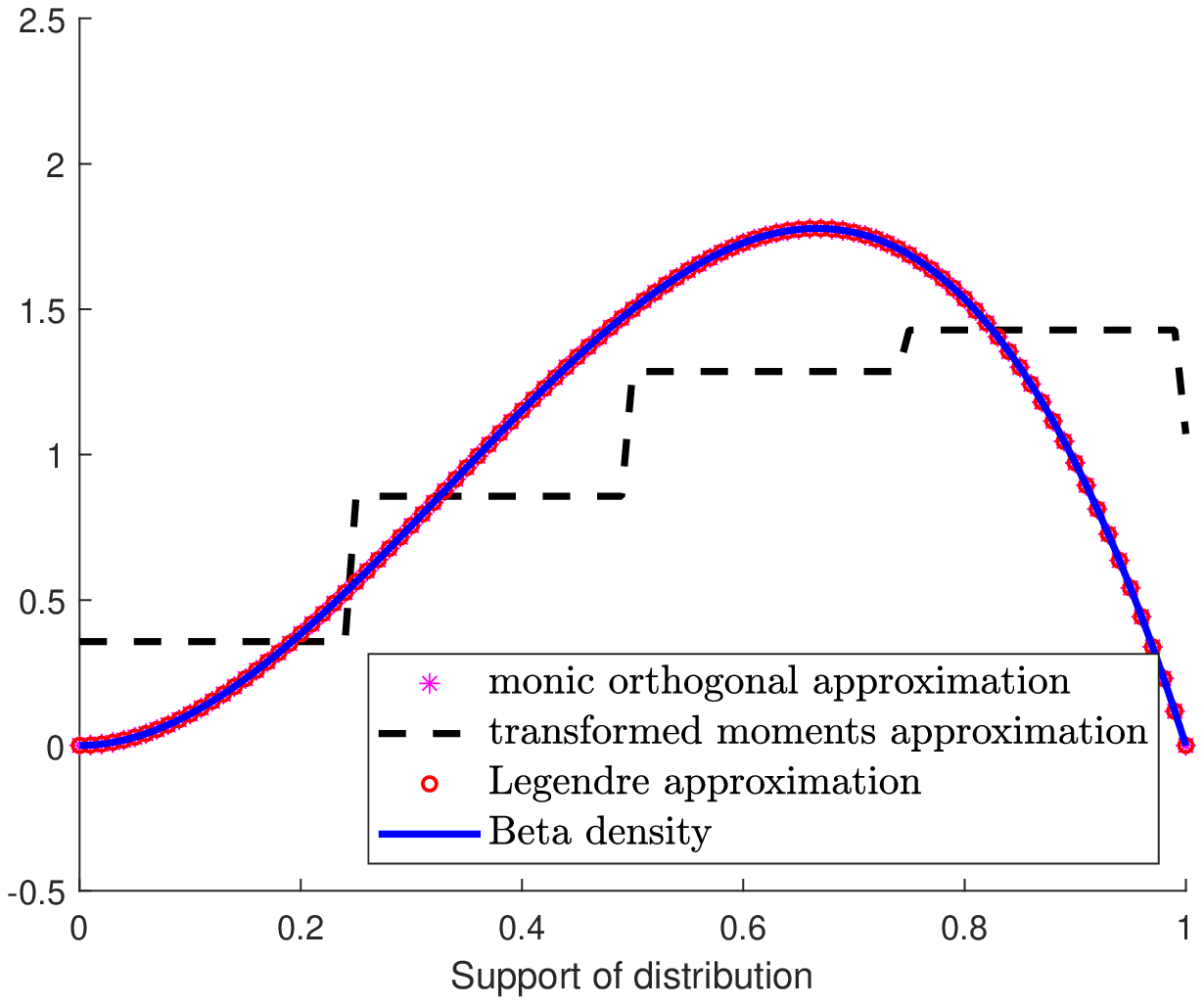}}
	\subfloat[\label{fig:PDF_approx_BetaRV_nmoms10}\ PDF $\nmoms=10$]{\includegraphics[width=0.5\textwidth]{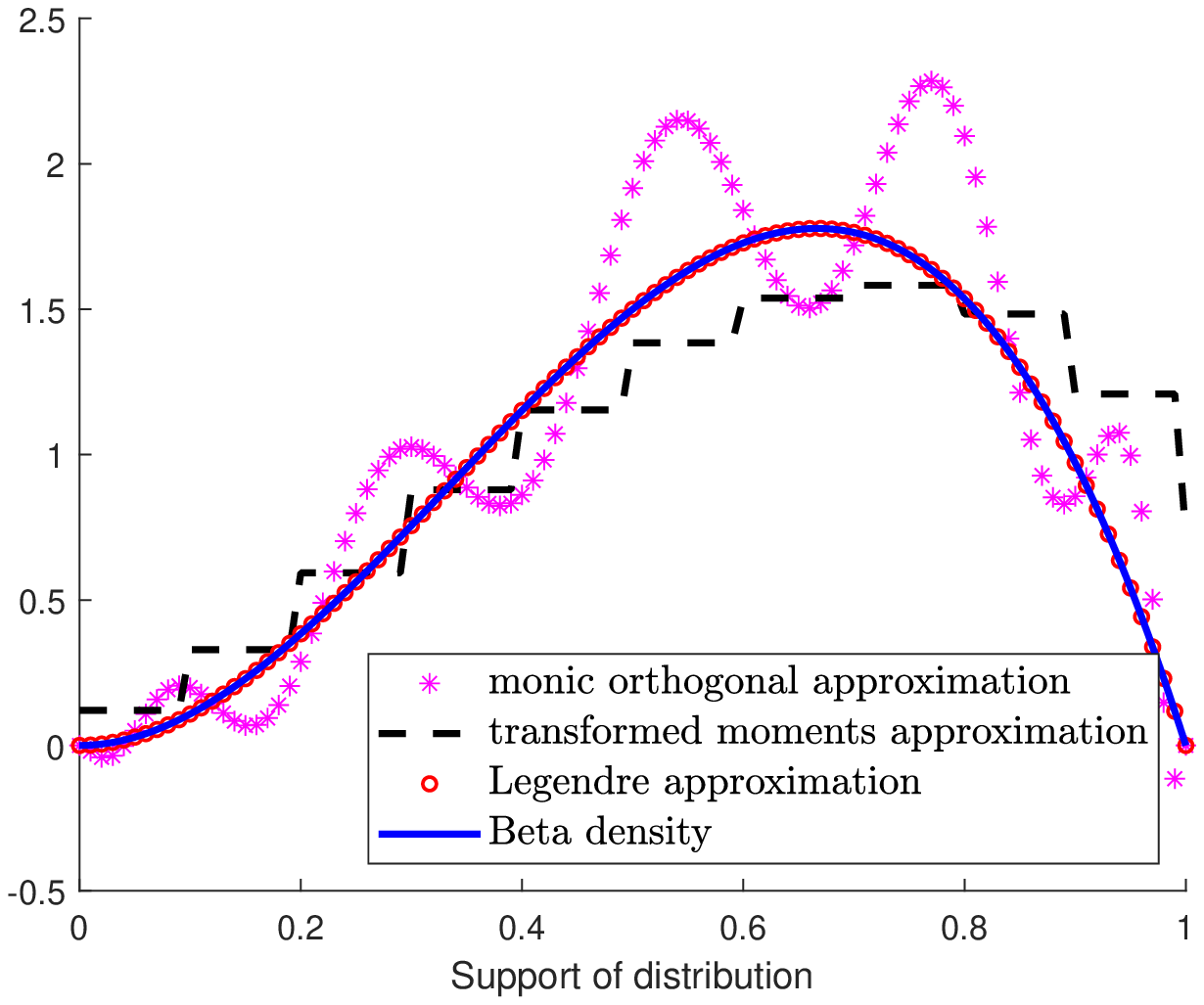}}
	
	\caption{Reconstruction of the PDF of the Beta RV $X\sim\Beta(3,2)$; left figure shows the excellent performance of approximations \eqref{eq:PDF_approx_legPol} and \eqref{eq:PDF_approx_monic} for a small number of moments used; right figure shows that performance for \eqref{eq:PDF_approx_trafoMom} improves significantly with the higher number of moments while \eqref{eq:PDF_approx_monic} starts oscillating.}
	
	\label{fig:PDF_approx_BetaRV}
\end{figure}

In Figure \ref{fig:PDF_approx_BetaRV}, we observe that the PDF approximation based on Legendre polynomials \eqref{eq:PDF_approx_legPol} provides excellent results provided that the number of moments used is not too high, which is beneficial. A magnitude of around $\nmoms=10$ turned out to be appropriate here. 

The PDF approximation based on monic orthogonal polynomials \eqref{eq:PDF_approx_monic} approximates the true PDF well already for a very small number of moments (see Figure \ref{fig:PDF_approx_BetaRV_nmoms4}). However, it suffers from oscillations that arise in our numerical simulations as can be seen in Figure \ref{fig:PDF_approx_BetaRV_nmoms10}. The PDF approximation \eqref{eq:PDF_approx_trafoMom} based on transformed moments does not lead to satisfactory results for a number of moments $\nmoms \in \{4,10\}$ (see Figure \ref{fig:PDF_approx_BetaRV}).

\bigskip

Above results suggest that the number of moments needed to produce a reasonable approximation of the probability distribution is highly dependent on the method of reconstruction used.

\begin{figure}[h!]
	\subfloat[\ $N=2$ ]{\begin{overpic}[width=0.5\textwidth]{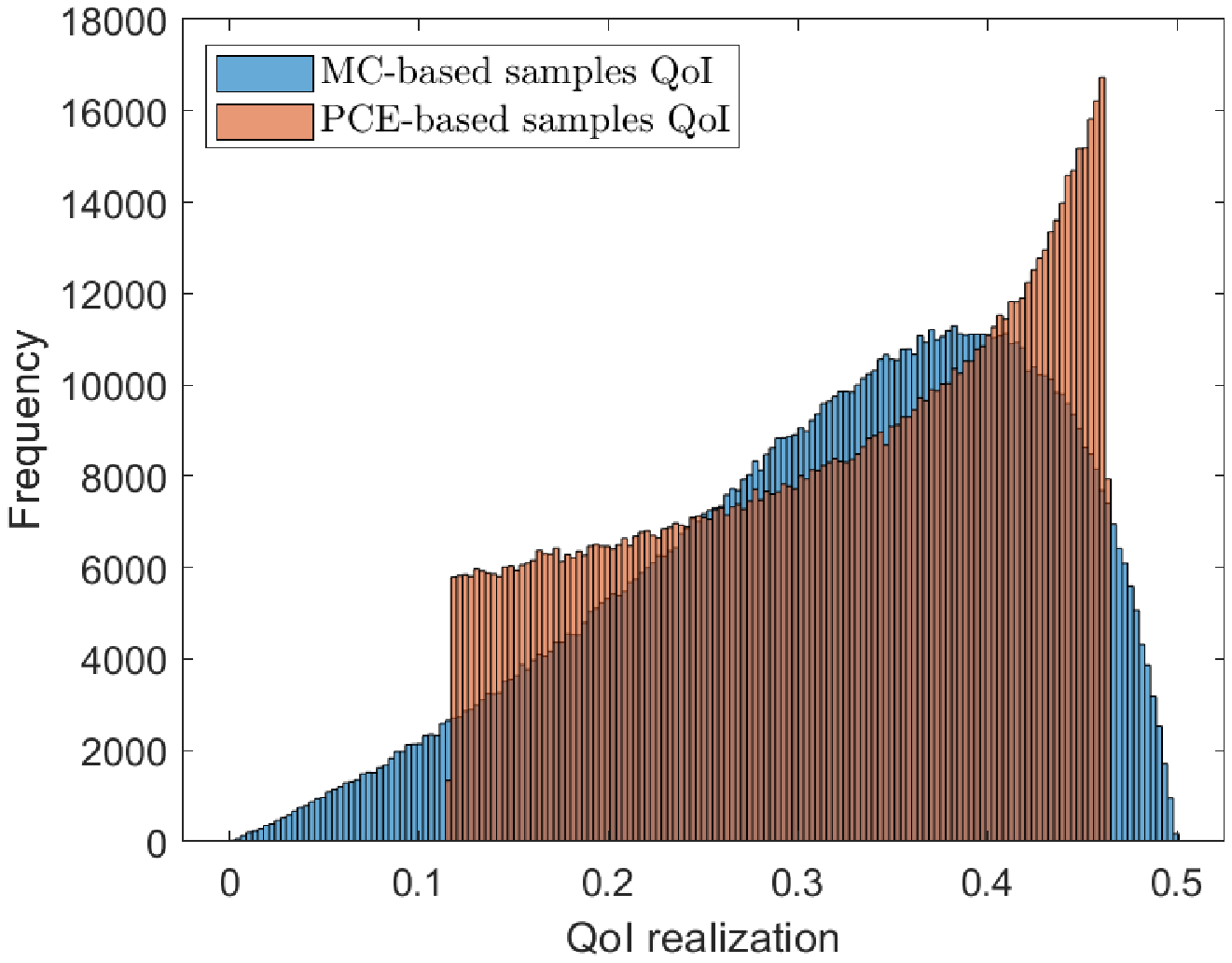}\end{overpic}}
	\subfloat[\ $N=14$ ]{\begin{overpic}[width=0.5\textwidth]{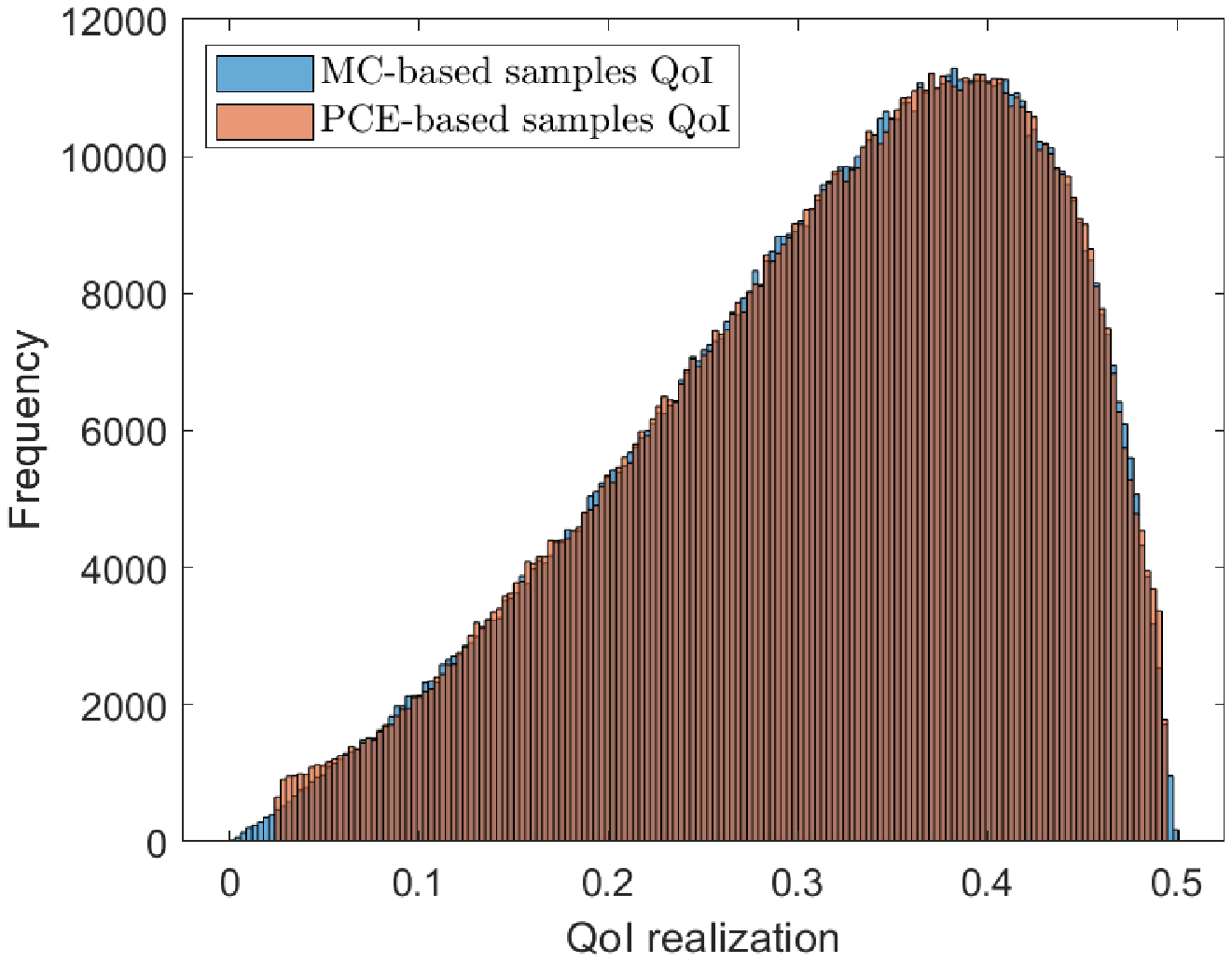}\end{overpic}}
	
	\caption{Approximation quality of the PCE (Legendre polynomial basis) of the quantity of interest (QoI) $\frac{\zeta}{1+\zeta}$ from \eqref{eq:redLorenzSyst} for $\zeta \sim \Beta(2,2)$ calculated via \textit{UQLab} \cite{Sudret.2014}; the approximation quality itself is limited in the left figure but we will see that the moments are still well captured (see Table \ref{tab:moms_LorenzRV_Beta22_Gamma81_LegBasis_N2_numMoms5_numComp2}); the right figure shows that the quality can be enhanced by choosing a higher truncation $N$ possibly causing numerical issues in the calculation of \eqref{eq:MellinTransform_PCE}.}
	
	\label{fig:approxQuali_X_div_1plusX_Beta22}
\end{figure}

We will see later that, although theoretically provided by formula \eqref{eq:MellinTransform_PCE}, the number of moments that can be reasonably calculated numerically is somewhat limited due to a high number of coefficients summed up and the possibly badly scaled product of coefficient $\hat{c}_i(s)$ and Mellin transform evaluation $\MM(\xi)(i+1)$. In the light of this observation, the method of transformed moments does not seem particularly suited for the purpose of reconstructing the probability distribution of the bifurcation coefficient.
Note however that the numerical issues in the calculation of the moments based on equation \eqref{eq:MellinTransform_PCE_Uniform} can be significantly reduced by choosing a very low precision $N$ in the PCE \eqref{eq:truncPCEnonlinTrafo_individual}. The choice of a low truncation index should not harm the method too much as it is an intermediate step in the calculation of the moments of the bifurcation coefficient $X$. A possible indication on a suitable truncation of the PCE is given by the Mellin-based moments of $X$. If the latter are sufficiently close to the Monte Carlo moments, the cut off of the PCE seems acceptable. An example illustrates the distinction between capturing well the overall distribution and providing decent moment estimates: in Figure \ref{fig:approxQuali_X_div_1plusX_Beta22}, we see PCE approximations for $N=2$ and $N=14$ calculated using the MATLAB-based software framework \textit{UQLab} \cite{Sudret.2014}. The PCE approximation itself for $N=2$ would not serve as a good approximation of the transformed random variable $\g(r_1) = \frac{r_1}{1+r_1}$. Being combined to the bifurcation coefficient $X$ from \eqref{eq:redLorenzSyst}, where the moments of $\frac{1}{r_2}$ are calculated analytically, the overall moments of $X$ are very well captured (see Table \ref{tab:moms_LorenzRV_Beta22_Gamma81_LegBasis_N2_numMoms5_numComp2}), which is our primary purpose of the PCE here.

\begin{table}[h!]
	\centering
	\begin{tabular}{lrrrrr}
		   & \multicolumn{1}{l}{$\mu_1$} & \multicolumn{1}{l}{$\mu_2$} & \multicolumn{1}{l}{$\mu_3$} & \multicolumn{1}{l}{$\mu_4$} & \multicolumn{1}{l}{$\mu_5$} \\
		\textbf{Mellin-based} & 4.55E-02 & 2.66E-03 & 1.99E-04 & 1.94E-05 & 2.60E-06 \\
		\textbf{MC}    & 4.54E-02 & 2.67E-03 & 2.00E-04 & 1.93E-05 & 2.46E-06 \\
	\end{tabular}%
	\caption{Comparison of moments for bifurcation coefficient $X$ from \eqref{eq:redLorenzSyst} for $\zeta \sim \Beta(2,2)$ and $\theta\sim \Gamdist(8,1)$ with truncation $N=2$.}
	\label{tab:moms_LorenzRV_Beta22_Gamma81_LegBasis_N2_numMoms5_numComp2}%
\end{table}%

\subsubsection{Applicability of GMMs} \label{subsubsec:appApproaches_GMMs}

The term \textit{method of moments} is rather generic and a specification of the way it is used is in order. Here, we use the moment conditions and minimize the corresponding quadratic form as explained in \cite{Hansen.2008}. To state the moment conditions, we first need some notation. We denote the $n$-th moment of the Gaussian mixture distribution by $m_n(\eta)$, where $\eta \in \PP\subset\R^{\npara}$ denotes the parameters defining the mixture distribution belonging to the parameter space $\PP\subset\R^{\npara}$. In case of a $k$-component Gaussian mixture distribution, the parameters are given by $\eta = (\pi_1,\dots,\pi_{k-1},\mu_1,\dots,\mu_k,\sigma_1,\dots,\sigma_k)$ with $\pi_i, i=1,\dots,k-1,$ being the mixing proportions, $\mu_i, i=1,\dots,k,$ the means and $\sigma_i, i=1,\dots,k,$ the standard deviations of the corresponding normal component distributions.

A suitable choice of $f$ in the moment conditions in \cite{Hansen.2008} for our setting leads to the moment conditions
\begin{align}
	E\left[f(X,\eta)\right] &= E\begin{bmatrix}
		X-m_1(\eta) \\ X^2 - m_2(\eta) \\ \vdots \\ X^{\nmoms} - m_{\nmoms}(\eta)
	\end{bmatrix} = E\begin{bmatrix}
	X \\ X^2 \\ \vdots \\ X^{\nmoms}
\end{bmatrix} - \begin{bmatrix}
m_1(\eta) \\ m_2(\eta) \\ \vdots \\ m_{\nmoms}(\eta)
\end{bmatrix} = \begin{bmatrix}
\MM(X)(2) \\ \MM(X)(3) \\ \vdots \\ \MM(X)(\nmoms+1)
\end{bmatrix} - \begin{bmatrix}
m_1(\eta) \\ m_2(\eta) \\ \vdots \\ m_{\nmoms}(\eta)
\end{bmatrix}, \label{eq:momentCond}
\end{align}
where $X$ denotes the bifurcation coefficient. Here, we benefit from the knowledge of the Mellin transform of the bifurcation coefficient in closed form. A big advantage is that we do not need to rely on sample estimates of the expected value in \eqref{eq:momentCond} as they are usually used in the generalized method of moments (see e.g. \cite{Hansen.2008,Wu.2020}). This is particularly beneficial for higher order moments that often need a high number of samples. Here, we do not need any samples at all. Note however, that the Mellin transform of the bifurcation coefficient in \eqref{eq:momentCond} might in fact be replaced by a version, where parts of the bifurcation coefficient are approximated via a PCE and thus, inducing an approximation error in \eqref{eq:momentCond} as well.  A rigorous error analysis is out of the scope of the present work. However, as we have already seen in Example \ref{ex:MellinTransformPCEpartLorenzRV} and Table \ref{tab:moms_LorenzRV_Beta22_Gamma81_LegBasis_N2_numMoms5_numComp2}, the approximation procedure still works reasonably well.

As in \cite{Hansen.2008}, the minimization problem via the quadratic form can be stated as
\begin{align}
	\argmin_{\eta \in \PP} E\left[f(X,\eta)\right]^T W E\left[f(X,\eta)\right],  \label{eq:minProbGMM}
\end{align}
where $W$ is a positive definite weighting matrix. Different choices of weighting matrices might lead to different estimators at least by using sample counterparts (see \cite{Hansen.2008}). However, our focus here is not on finding the most performant weighting matrix or optimization algorithm. To avoid facing a badly scaled optimization problem due to possibly huge differences of the magnitudes of the moments, we choose $W$ as the positive definite diagonal matrix with entries equal to the inverse of the absolute value of the moments of the bifurcation coefficient, i.e.
\begin{align*}
	W &= \begin{pmatrix}
		\lvert\frac{1}{\MM(X)(2)}\rvert & & 0 \\
		& \ddots & \\
		0 & & \lvert\frac{1}{\MM(X)(\nmoms+1)}\rvert
	\end{pmatrix}.
\end{align*}
To ensure that $W$ is well defined the Mellin transforms must not be zero. This is in particular an issue for symmetric distributions for which all odd moments are equal to zero and if the mean of the bifurcation coefficient itself is zero. An easy remedy would be to choose $W$ to be the identity matrix $W = I$. For more sophisticated choices of the weighting matrix, we refer the reader to \cite{Wolf.2016} and references therein. For our numerical implementation of the minimization problem \eqref{eq:minProbGMM}, we use the MATLAB solver \textit{fmincon}\footnotemark[1]\footnotetext[1]{\url{https://de.mathworks.com/help/optim/ug/fmincon.html}, last checked: November 17, 2020}.

According to \cite{Wu.2020}, the optimization problem arising from the generalized method of moments for Gaussian mixture models is not efficiently solvable as the moment conditions are non-convex in that case. However, as we aim at obtaining a first estimate of the probability of the bifurcation types based on the knowledge of the first $\nmoms+1$ moments only, using the generalized method of moments in the estimation of the Gaussian mixture models serves our purpose. Note that the solution found by the solver depends on the initialization, which is not surprising in the context of non-convex optimization. However, the different solutions found gave reasonable estimates of the probabilities. In particular, when dealing with a two-component mixture, setting the initial means nearby the boundaries of the support lead to very satisfying results.

\bigskip

Due to the drawbacks of the polynomial approximation mentioned in \ref{subsubsec:appApproaches_poly}, in our bifurcation context, we prefer the approach via the method of moments to fit a Gaussian mixture distribution, for which we already showed a first promising result in Figure~\ref{fig:PDF_approx_teaserLorenzSyst}.
Accordingly, in a next step, we use the approximation procedure of Section \ref{subsubsec:semiParaEst_MethodOfMoments} to estimate the probability density function of the bifurcation coefficient. Further details on the performance of the polynomial approximations can be found in Appendix \ref{app:probBifGMM_poly}.

We start by introducing some numerical test settings in Table \ref{tab:numericalSets} that will be considered in the remainder of this section. The abbreviation $\text{gen}\Beta_{a}^{b}$ stands for the generalized Beta distribution with support $[a,b]$. Although $r_1=\zeta$ has assumed to be an $P$-a.s.\ positive random variable in the Lorenz system \eqref{eq:LorenzSyst}, it is here allowed to be negative for conceptual purposes. Note again that our methodology is not restricted to a particular RODE. An example from computational neuroscience that allows for negative parameter values, can be found in \cite{Kuehn.2009}.

\begin{table}[h]
	\centering
	\begin{adjustbox}{width=\textwidth}
	\begin{tabular}{c|c|c|c|c}
		& \textbf{PS A} & \textbf{PS B} & \textbf{PS C} & \textbf{PS D} \\
		\hline
		Input parameter $r_1$ & $\zeta \sim \UU(4,6)$ & $\zeta \sim \Beta(2,2)$ & $\zeta \sim \Beta(2,5)$ & $\zeta \sim \text{gen}\Beta_{-0.5}^{0.5}(2,5)$\\
		\hline
		Input parameter $r_2$ & $\theta \sim \Gamdist(8,1)$ & PS A & PS A & PS A \\
		\hline
		PCE basis & Legendre polynomials & PS A & PS A & PS A \\
		\hline
		Truncation of PCE & $N=2$ & PS A & PS A & PS A\\
		\hline
		Number of moments & $\nmoms = 7$ & $\nmoms = 5$ & PS B & PS B \\
	\end{tabular}
	\end{adjustbox}
	\caption{Numerical setting for analysis of bifurcation coefficient $X$ from \eqref{eq:redLorenzSyst}.}
	\label{tab:numericalSets}%
\end{table}%

A modification of PS D including a truncated power law and an application of our methodology to the Hindmarsh-Rose model \cite{Liu.2012}, where positive parameters can lead to super- as well as subcritical bifurcations can be found in the supplementary material at \url{https://github.com/kerstinLux/UQbifurcation}.

\subsubsection{GMMs for the probability distribution of the bifurcation coefficient} \label{subsubsec:probBifGMM_GMMs}

Also carrying over the GMM approach to the estimation of non-standard PDFs, such as the one of the bifurcation coefficient $X$ in the reduced Lorenz system \eqref{eq:redLorenzSyst}, comes with several challenges.

Again, the number of moments that can be reasonably calculated is an important aspect. As mentioned above, for our purpose, it is beneficial to reduce the precision in the PCE to a low number $N$ to face less numerical issues in the calculation of the moments in equation \eqref{eq:MellinTransform_PCE_Uniform}.

We start our numerical analysis of the performance of GMMs in estimating the probability distribution of the bifurcation coefficient again with the numerical setting PS A from Table \ref{tab:numericalSets}.

\bigskip

\begin{figure}[h!]
	\subfloat[\label{fig:histVsPDFapprox_LorenzRV_Uni46_Gamma81_LegBasis_N2_numMoms7_numComp2_Maxfeval1e+4_maxIter5mal1e+3_StepTol1e-12_stopNondecr}\ {GMM estimation of PDF}]{\begin{overpic}[width=0.5\textwidth]{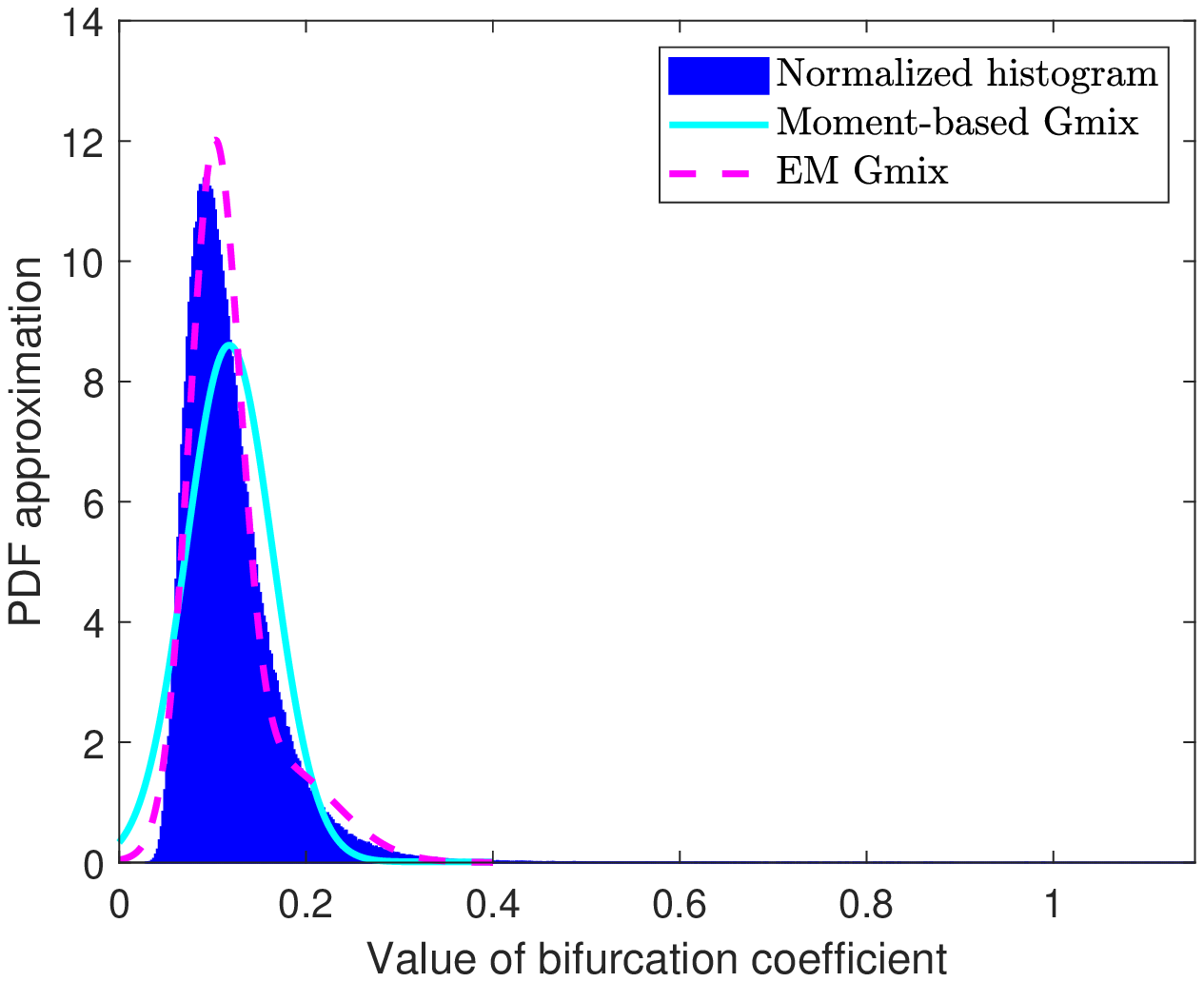}
		\end{overpic}
	}
	\subfloat[\ {GMM estimation of CDF}]{\begin{overpic}[width=0.5\textwidth]{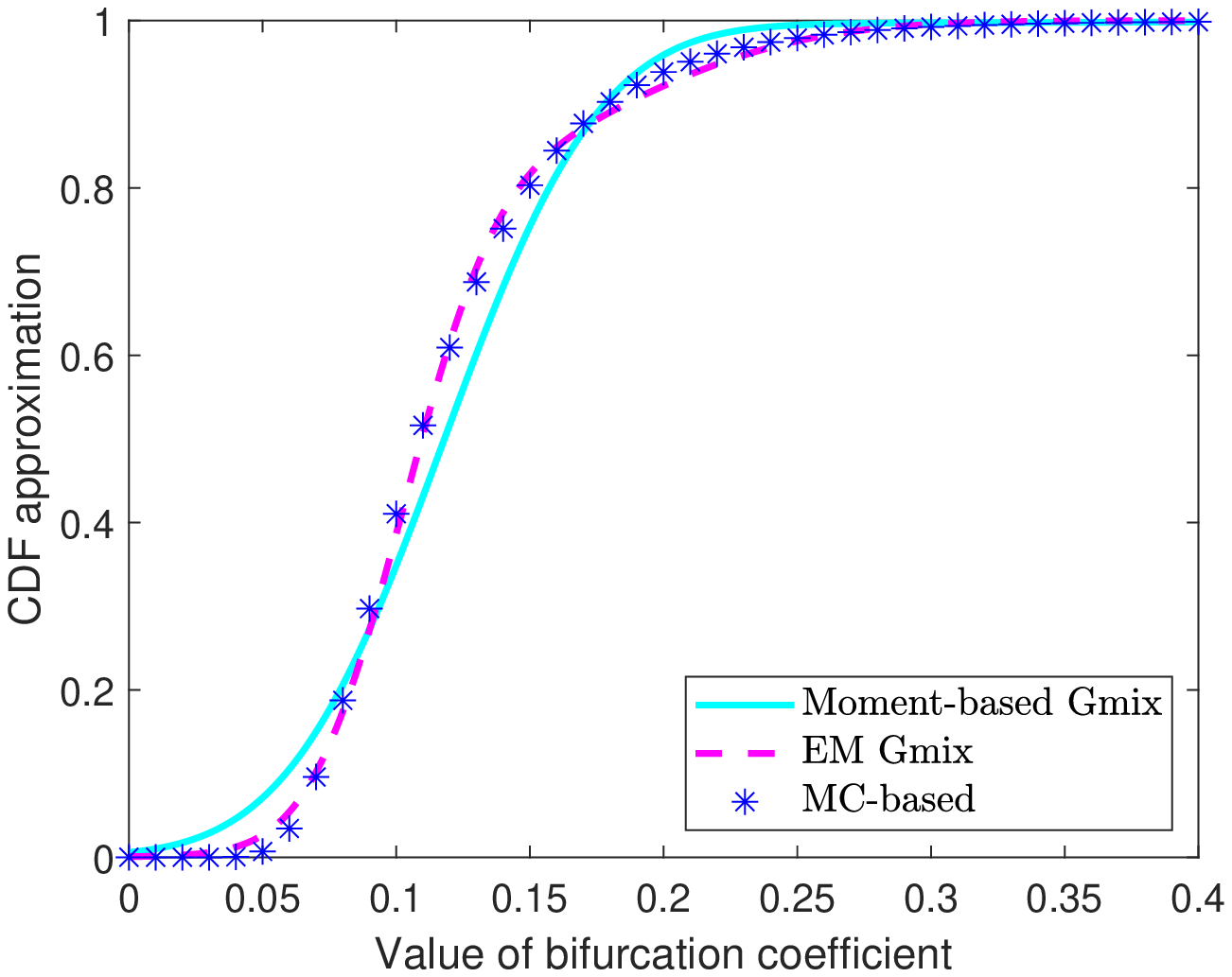}\end{overpic}}
	
	\caption{Reconstruction of PDF and CDF of $X$ in reduced Lorenz system \eqref{eq:redLorenzSyst} via Gaussian mixture (Gmix) model for PS A; our Mellin-moments-based GMM approximation of the PDF and CDF in turqoise is compared to the sample-based MATLAB solution (using the EM algorithm) in magenta and the blue normalized histogram/MC estimates: results are in good qualitative agreement.}
	
	\label{fig:GMMapprox_LorenzRV_Uni46_Gamma81_LegBasis_N2_numMoms7_numComp2_Maxfeval1e+4_maxIter5mal1e+3_StepTol1e-12_stopNondecr}
\end{figure}

In Figure \ref{fig:GMMapprox_LorenzRV_Uni46_Gamma81_LegBasis_N2_numMoms7_numComp2_Maxfeval1e+4_maxIter5mal1e+3_StepTol1e-12_stopNondecr}, we compare the GMM estimate of the PDF and the CDF of the bifurcation coefficient of the reduced Lorenz system with a MC sampling-based solution for $M=10^6$ samples and the GMM estimation obtained by using the MATLAB routine \textit{fitgmdist}\footnote{source: \url{https://de.mathworks.com/help/stats/fitgmdist.html} (based on the \textit{EM-algorithm}), last checked: November 21, 2020}, which gets all the $M=10^6$ samples as input instead of the $\nmoms = 7$ moments used for the moment-based GMM estimation.

\begin{table}[htbp]
	\centering
	\begin{adjustbox}{width=\textwidth}
	\begin{tabular}{lrrrrrrr}
		& \multicolumn{1}{l}{$\mu_1$} & \multicolumn{1}{l}{$\mu_2$} & \multicolumn{1}{l}{$\mu_3$} & \multicolumn{1}{l}{$\mu_4$} & \multicolumn{1}{l}{$\mu_5$} & \multicolumn{1}{l}{$\mu_6$} & \multicolumn{1}{l}{$\mu_7$} \\
		\textbf{MC}    & 1.1880E-01 & 1.6470E-02 & 2.7368E-03 & 5.6489E-04 & 1.5111E-04 & 5.4213E-05 & 2.6066E-05 \\
		\textbf{Mellin-based} & 1.1882E-01 & 1.6479E-02 & 2.7434E-03 & 5.7112E-04 & 1.5859E-04 & 6.6080E-05 & 5.5089E-05 \\
		\textbf{GMM estimates} & 1.1874E-01 & 1.6594E-02 & 2.6960E-03 & 5.5522E-04 & 1.6678E-04 & 7.5774E-05 & 4.5018E-05 \\
	\end{tabular}%
	\end{adjustbox}
	\caption{Moments for PS A: Mellin-based moments are in close agreement with the MC moments based on $M=10^6$ samples and GMM estimates match well the Mellin-based moments that appear in the underlying moment conditions \eqref{eq:momentCond}.}
	\label{tab:moms_LorenzRV_Uni46_Gamma81_LegBasis_N2_numMoms7_numComp2_Maxfeval1e+4_maxIter5mal1e+3_StepTol1e-12_stopNondecr}%
\end{table}%

\begin{table}[htbp]
	\centering
	\begin{tabular}{lrrrr}
		\textbf{y} & \textbf{0.00} & \textbf{0.05} & \textbf{0.10} & \textbf{0.35} \\
		\textbf{Moment-based Gmix} & 0.0057 & 0.0711 & 0.3489 & 0.9979 \\
		\textbf{MC} & 0.0000 & 0.0069 & 0.4106 & 0.9969 \\
	\end{tabular}%
	\caption{Exemplary values of CDF for PS A: probability for a subcritical bifurcation $P(X<0)$ is closely matched.}
	\label{tab:prob_LorenzRV_Uni46_Gamma81_LegBasis_N2_numMoms7_numComp2_Maxfeval1e+4_maxIter5mal1e+3_StepTol1e-12_stopNondecr}%
\end{table}%

Although we made assumptions on the shape of the bifurcation coefficient (in terms of being representable as a mixture of Gaussian distributions) and accepted that the supports of the distributions might not coincide, the results are in good qualitative agreement. We can also quantify this in numbers. In Table \ref{tab:moms_LorenzRV_Uni46_Gamma81_LegBasis_N2_numMoms7_numComp2_Maxfeval1e+4_maxIter5mal1e+3_StepTol1e-12_stopNondecr}, we show that all the $\nmoms=7$ moments calculated based on the Mellin transform are well matched by our GMM estimates and that they are very close to the MC moments.

In Table \ref{tab:prob_LorenzRV_Uni46_Gamma81_LegBasis_N2_numMoms7_numComp2_Maxfeval1e+4_maxIter5mal1e+3_StepTol1e-12_stopNondecr}, we see evaluations of the CDF in $y \in \{0, 0.05, 0.1, 0.35\}$.  These are the numbers we are interested in for estimating the probabilities of bifurcation types. In this case, the probability of observing a subcritical bifurcation ($P(X<0)$) would be zero, which is closely matched by our approximation.

\bigskip

We have already seen the very close match of our moment-based GMM approximation of the PDF for PS D from Table \ref{tab:numericalSets} in Figure \ref{fig:PDF_approx_LorenzRV}. The quantification for PS D is provided additionally in Tables \ref{tab:moms_LorenzRV_genBeta25_suppMin0K5_0K5_Gamma81_LegBasis_N2_numMoms5_numComp2_Maxfeval1e+4_maxIter5mal1e+3_StepTol1e-12_stopnondecr} and \ref{tab:prob_LorenzRV_genBeta25_suppMin0K5_0K5_Gamma81_LegBasis_N2_numMoms5_numComp2_Maxfeval1e+4_maxIter5mal1e+3_StepTol1e-12_stopnondecr}. Note the excellent match of all $\nmoms=5$ moments calculated based on the Mellin transform with the ones of the estimated GMM both being close to the MC moments. Moreover, the evaluations of the CDF in $y \in \{0, 0.05, 0.1, 0.35\}$ are matched to the precision of at least $10^{-2}$. In particular note again that the probability of observing a subcritical pitchfork bifurcation ($P(X<0)$) is very well matched. Results suggest that the symmetry in the support of $r_1$ even further reduces numerical inaccuracies.

\begin{table}[h]
	\centering
	\begin{tabular}{lrrrrr}
		& \multicolumn{1}{l}{$\mu_1$} & \multicolumn{1}{l}{$\mu_2$} & \multicolumn{1}{l}{$\mu_3$} & \multicolumn{1}{l}{$\mu_4$} & \multicolumn{1}{l}{$\mu_5$} \\
		\textbf{MC} & -4.6251E-02 & 4.0976E-03 & -4.8920E-04 & 8.0369E-05 & -1.7875E-05 \\
		\textbf{Mellin-based} & -4.6247E-02 & 3.9638E-03 & -4.5777E-04 & 7.1032E-05 & -1.5358E-05 \\
		\textbf{GMM estimates} & -4.6246E-02 & 3.9646E-03 & -4.5596E-04 & 7.3300E-05 & -1.4255E-05 \\
	\end{tabular}%
	\caption{Moments for PS D: Mellin-based moments are in close agreement with the MC moments based on $M=10^6$ samples and GMM estimates give an excellent match of the Mellin-based moments that appear in the underlying moment conditions \eqref{eq:momentCond}.}
	\label{tab:moms_LorenzRV_genBeta25_suppMin0K5_0K5_Gamma81_LegBasis_N2_numMoms5_numComp2_Maxfeval1e+4_maxIter5mal1e+3_StepTol1e-12_stopnondecr}%
\end{table}%

\begin{table}[h!]
	\centering
	\begin{tabular}{lrrrr}
		\textbf{y} & \textbf{-0.15} & \textbf{-0.05} & \textbf{0} & \textbf{0.05} \\
		\textbf{Moment-based Gmix} & 0.0286 & 0.4069 & 0.9049 & 0.9954 \\
		\textbf{EM Gmix} & 0.0332 & 0.3928 & 0.8883 & 0.9971 \\
		\textbf{MC} & 0.0264 & 0.3950 & 0.8903 & 0.9992 \\
	\end{tabular}%
	\caption{Exemplary values of CDF for PS D: probability for a subcritical bifurcation $P(X<0)$ is very well matched to the precision of at least $10^{-2}$.}
	\label{tab:prob_LorenzRV_genBeta25_suppMin0K5_0K5_Gamma81_LegBasis_N2_numMoms5_numComp2_Maxfeval1e+4_maxIter5mal1e+3_StepTol1e-12_stopnondecr}%
\end{table}%

\begin{figure}[h!]
	\centering
	\subfloat[\label{fig:PDF_approx_Beta22}\ PDF PS B]{\begin{overpic}[width=0.5\textwidth]{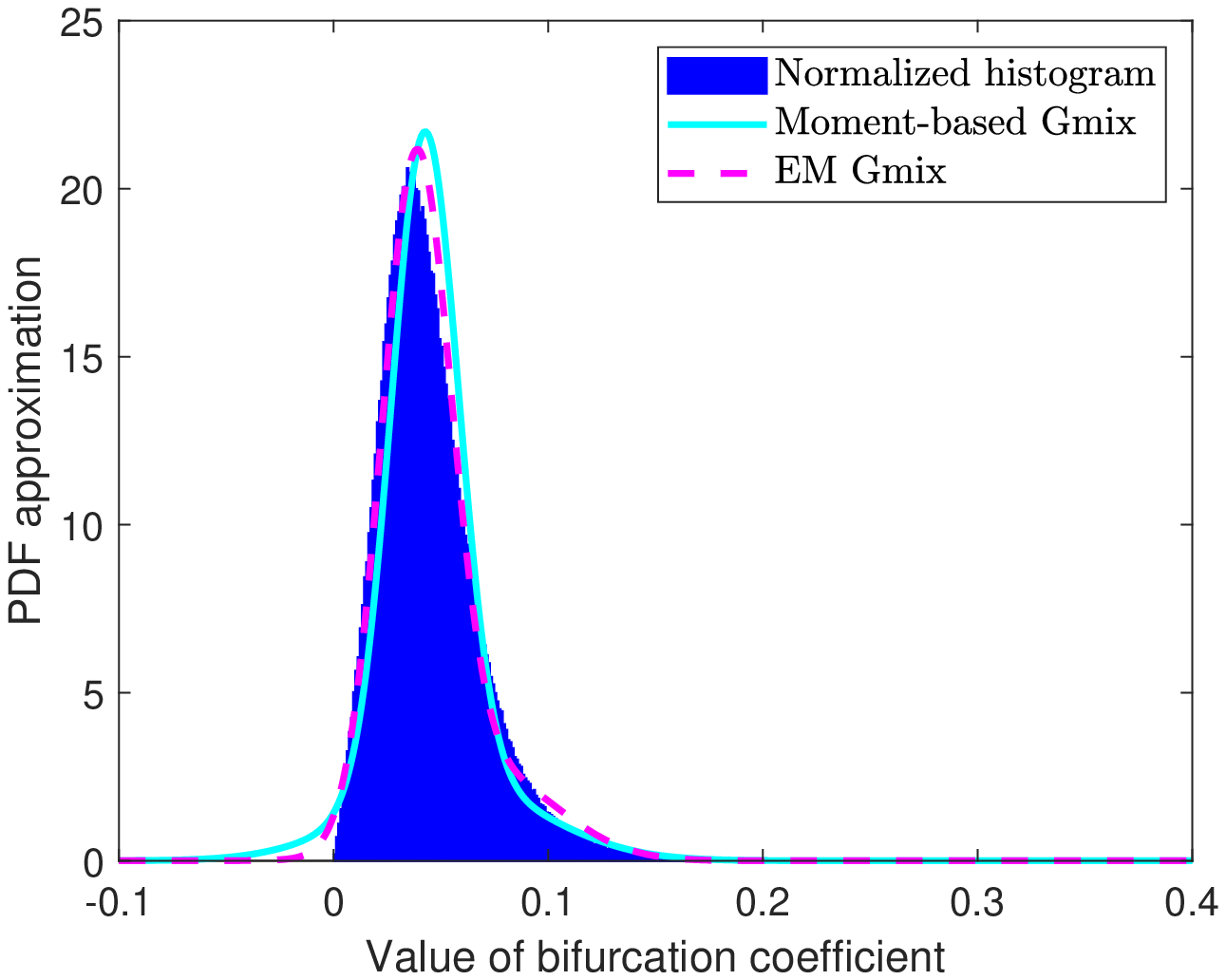}\end{overpic}}
	\subfloat[\label{fig:CDF_approx_Beta22}\ CDF PS B]{\begin{overpic}[width=0.5\textwidth]{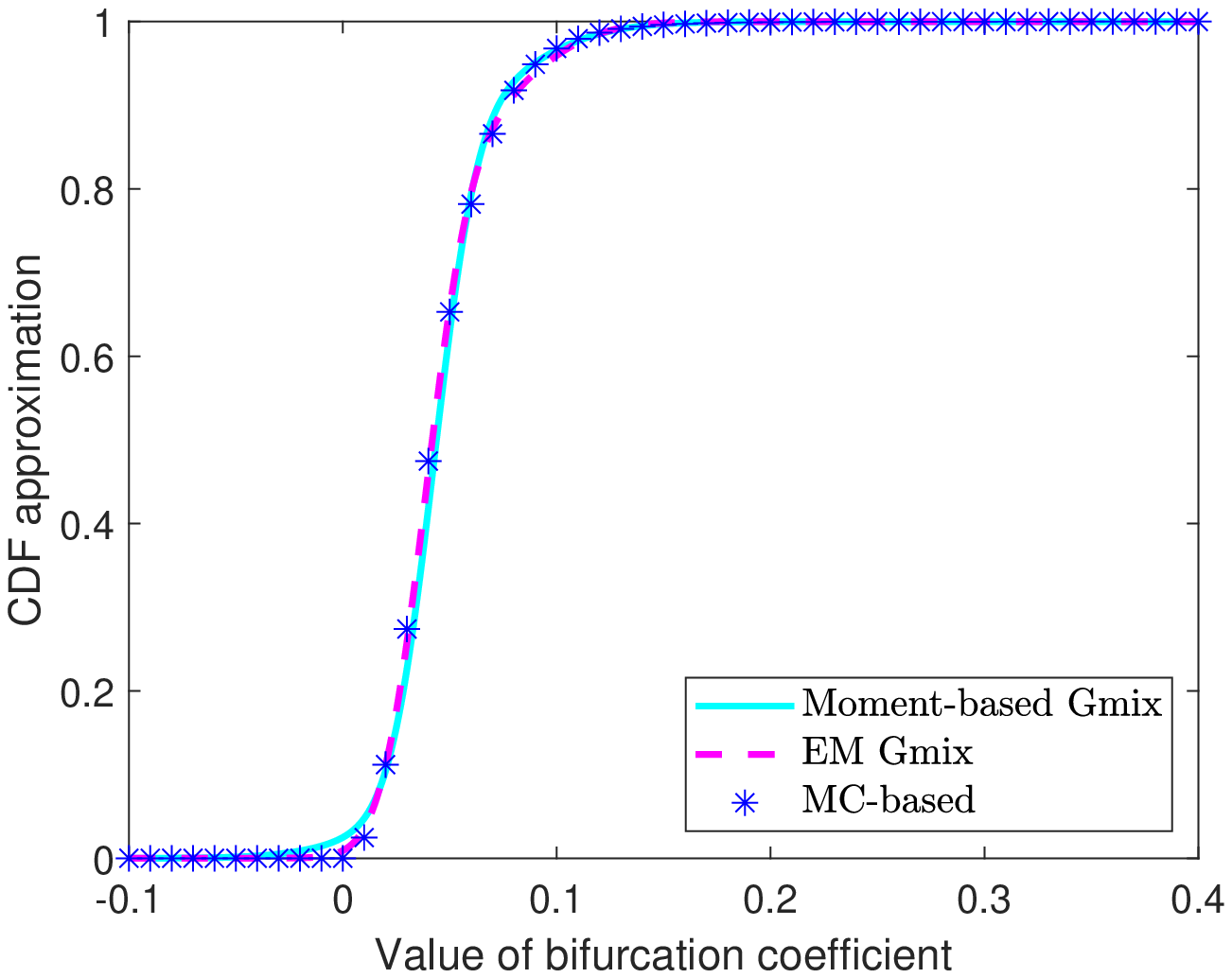}\end{overpic}}
	
	\subfloat[\label{fig:PDF_approx_genBeta25}\ PDF PS C]{\begin{overpic}[width=0.5\textwidth]{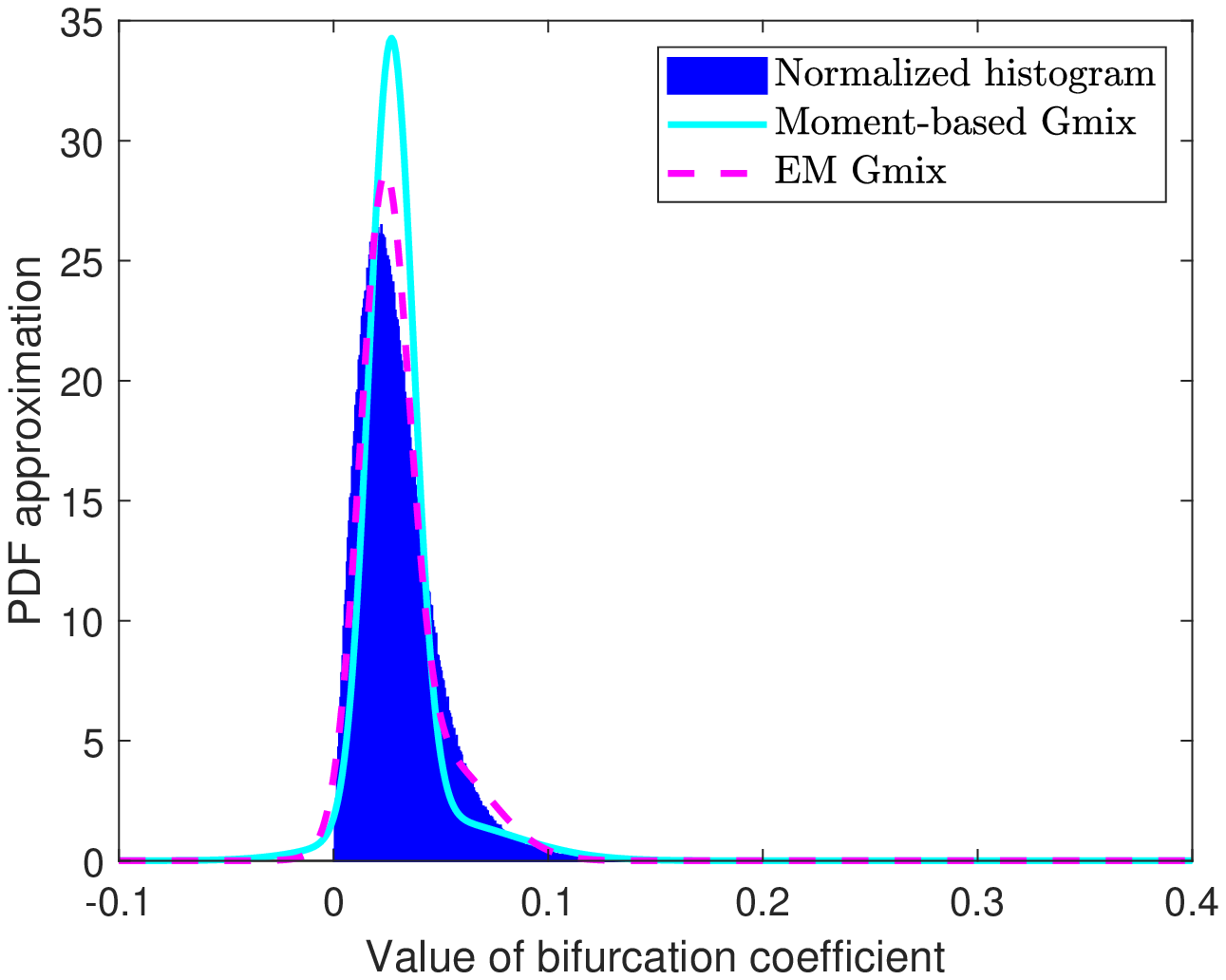}\end{overpic}}
	\subfloat[\label{fig:CDF_approx_genBeta25}\ CDF PS C]{\begin{overpic}[width=0.5\textwidth]{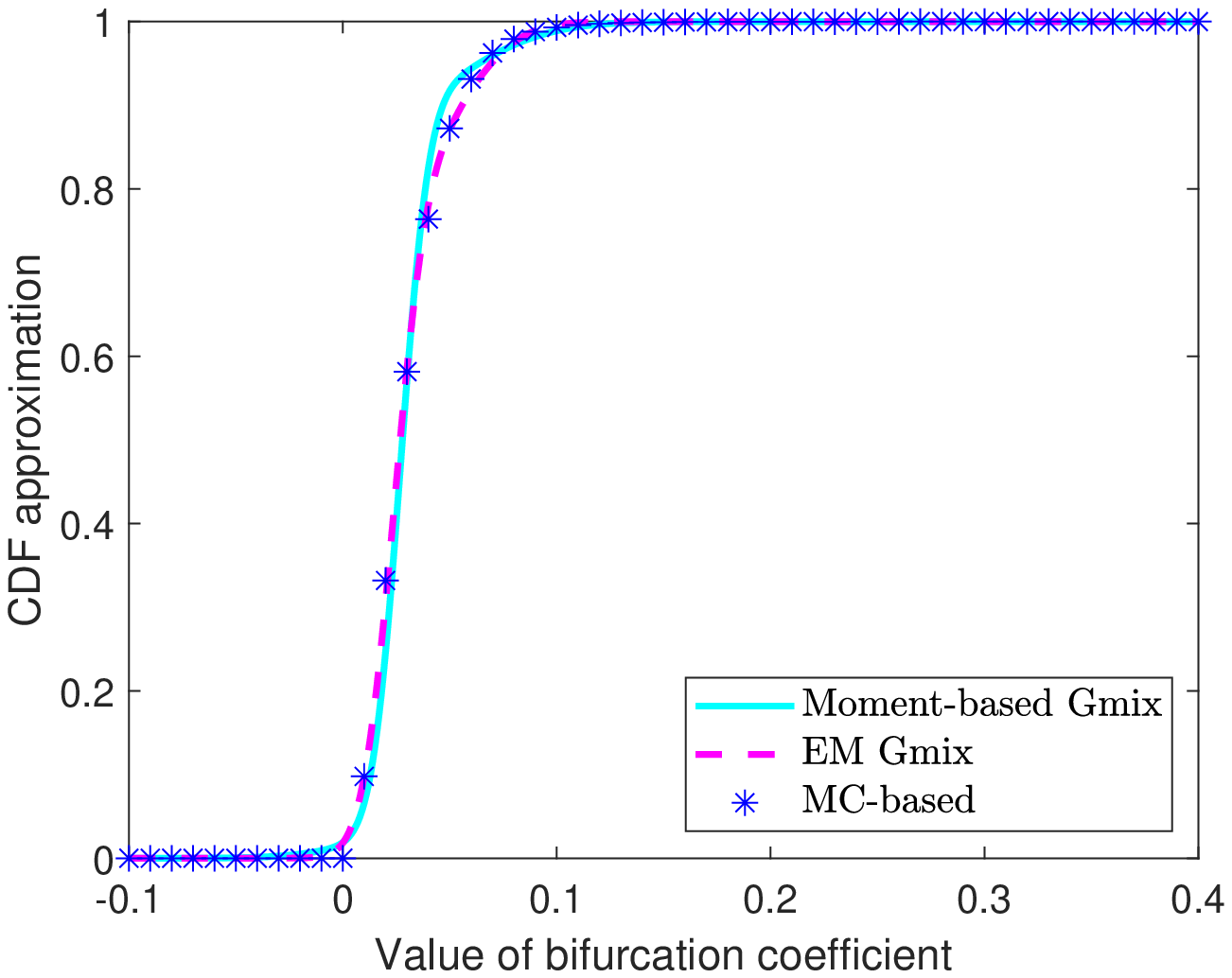}\end{overpic}}
	
	\caption{Reconstruction of PDF and CDF of $X$ in reduced Lorenz system \eqref{eq:redLorenzSyst} via GMM for PS B and PS C from Table \ref{tab:numericalSets}: our Mellin-moments-based GMM approximations of the PDF and CDF in turqoise are compared to the sample-based MATLAB solutions (using the EM algorithm) in magenta and the blue normalized histogram/MC estimates: results are in very close agreement.}
	
	\label{fig:GMM_approx_Betainput}
\end{figure}

\bigskip
The corresponding results for the numerical settings PS B and PS C underline the excellent performance of our method for the approximation of the probability distribution of the bifurcation coefficient and are given in Figure \ref{fig:GMM_approx_Betainput}. We see an even more precise match both in terms of PDF and CDF approximation than for PS A. This might be due to the absence of discontinuities in the Beta PDF of $r_1$ in PS B and PS C compared to the uniform input PDF of $r_1$ in PS A.


\section{Sampling-based approach for uncertainty propagation: unscented Transformation} \label{sec:samplingBasedAppUncertaintyProp}
Here, we make a cut and present a purely sampling-based alternative if one wants to obtain a very rough first estimate of the bifurcation probability.
The unscented Kalman filter idea has been proposed by Uhlmann in his PhD thesis \cite{Uhlmann.1995} in 1995. Julier, Uhlmann and Durrant-Whyte then came up with the unscented Kalman filter in \cite{Julier.1995} in 1995. It explicitly addresses drawbacks of the extended Kalman filter. It is based on a so-called \textit{unscented transformation}. The basic idea of this unscented transformation is to determine a set of so-called \textit{sigma points} by a deterministic algorithm that captures a predefined number of moments of the input distribution. This number depends on whether the input follows a Gaussian distribution or not, but a non-Gaussian input is possible. This sigma point set can be interpreted as a discrete approximation of the input probability distribution. The unscented Kalman filter belongs to the more general class of sigma-point methods.
\begin{figure}[t]
	\centering
	\subfloat[\ $a_1=0, b_1=1, a_2=0, b_2=1$]{\includegraphics[width=0.46\textwidth]{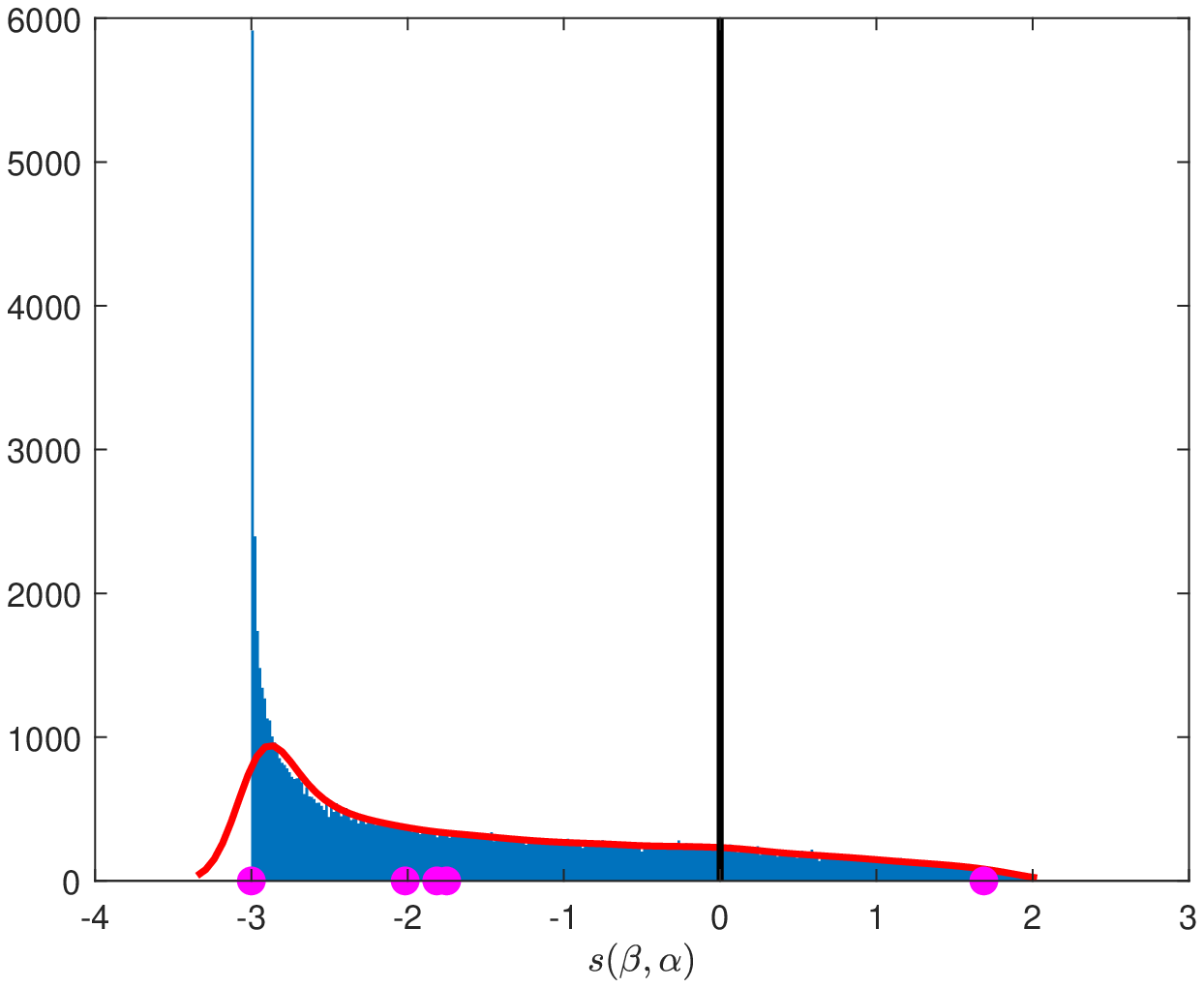}} \hspace{0.2cm}
	\subfloat[\ $a_1=0, b_1=1, a_2=0.6, b_2=1$]{\includegraphics[width=0.46\textwidth]{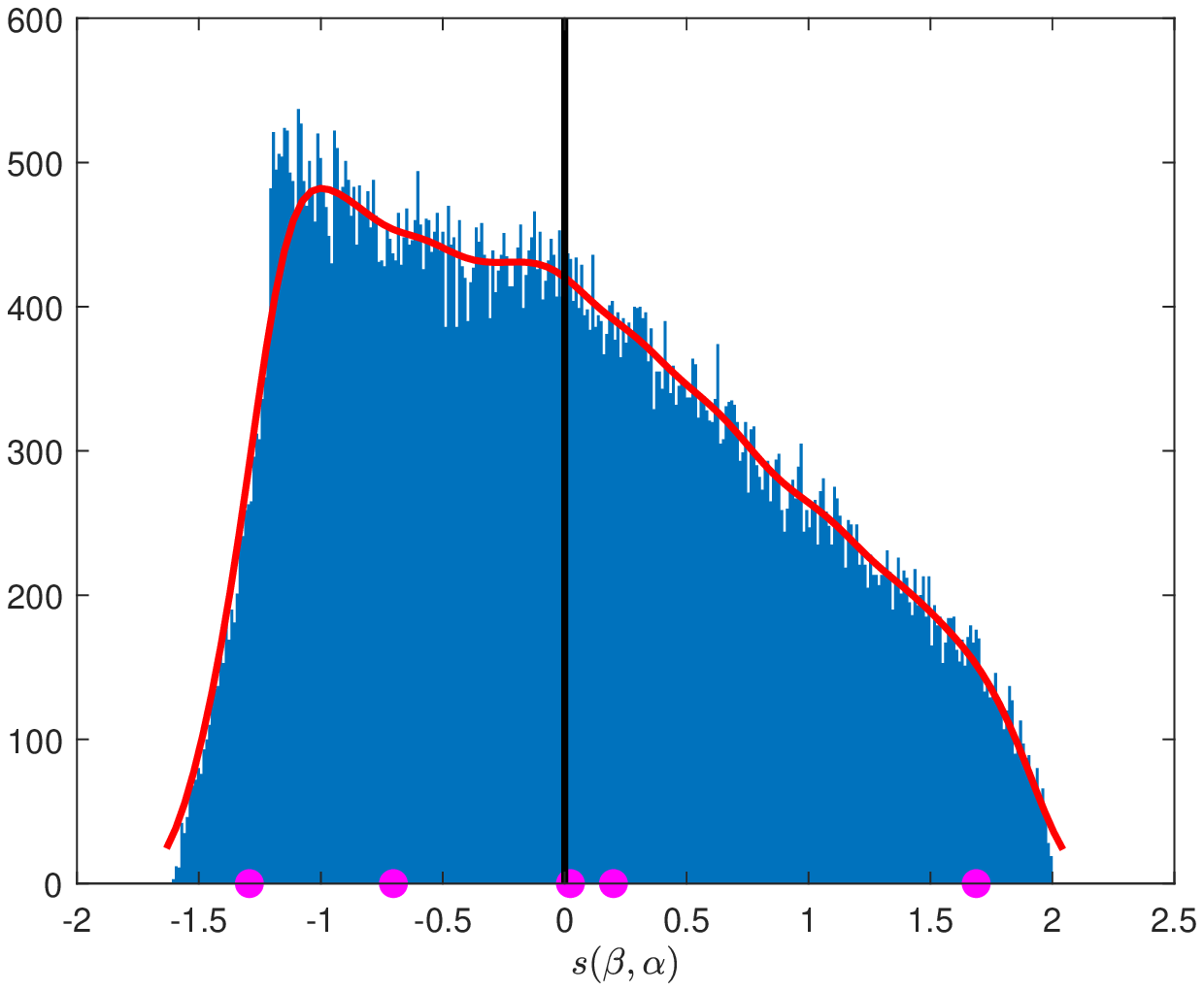}}
	
	\subfloat[\ $a_1=0, b_1=1, a_2=0.9, b_2=1$]{\includegraphics[width=0.46\textwidth]{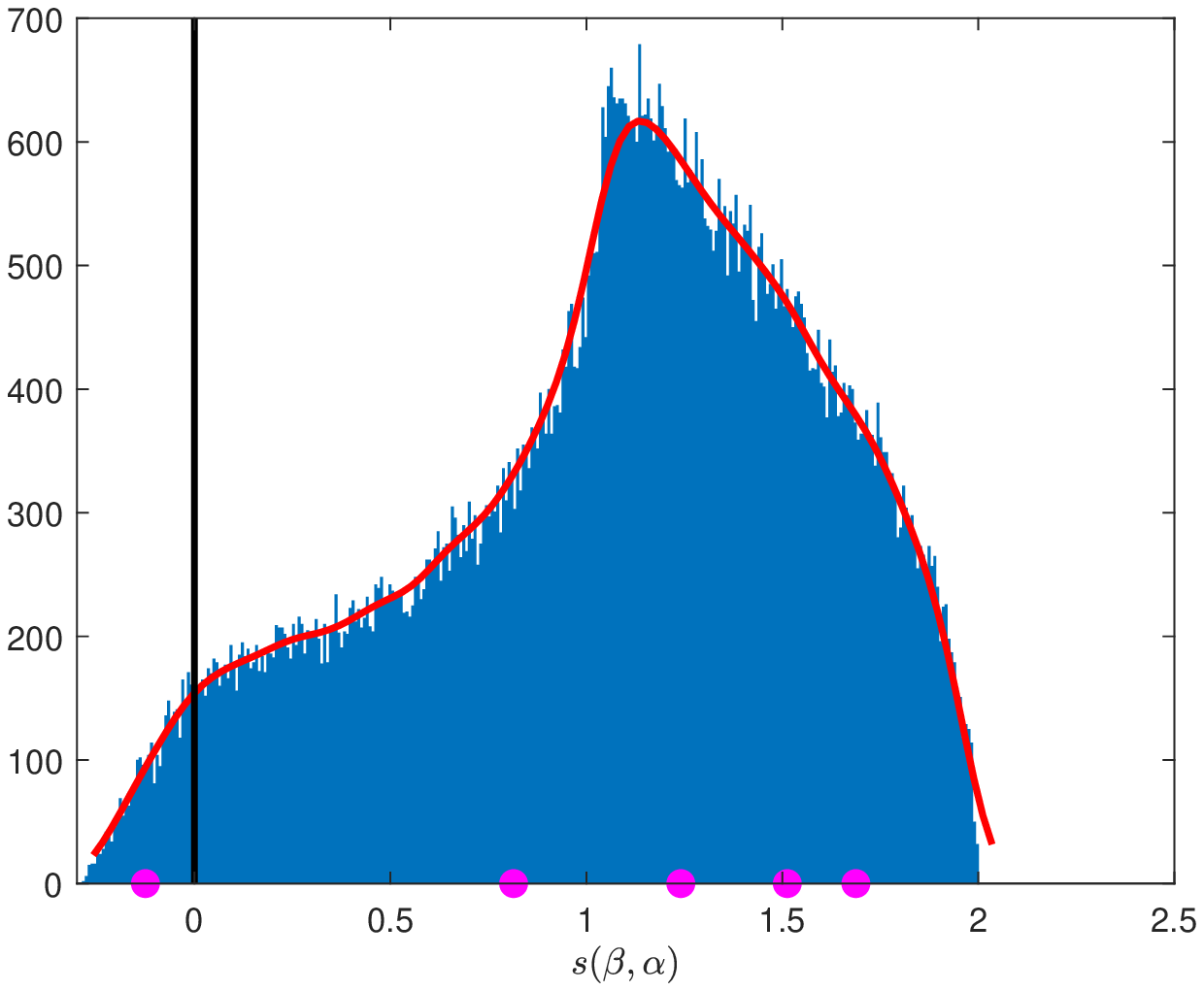}} \hspace{0.2cm}
	\subfloat[\ $a_1=0, b_1=0.2, a_2=0.7, b_2=0.95$]{\includegraphics[width=0.46\textwidth]{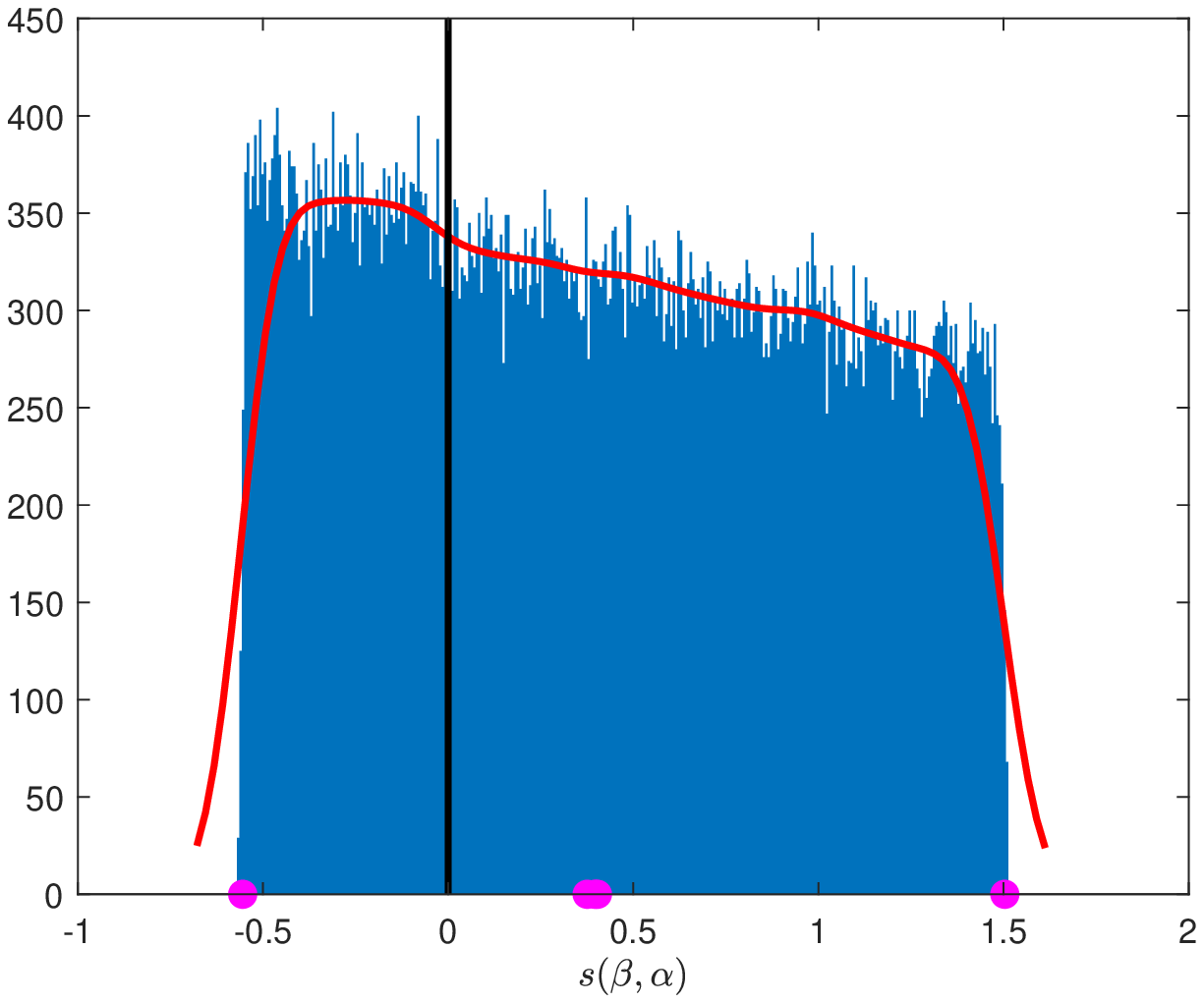}}
	
	\caption{Probability distribution of sign of Lyapunov coefficients determined by $s(\beta,\alpha)$ in \eqref{eq:sign_LyapCoeff_WattGovernor} for different supports of uniform input distribution of $\alpha$ and $\beta$; precision 3 is chosen in sigma point selection of \cite{Julier.2004}; red line is MATLAB kernel density estimate of blue histogram. The sigma points propagated through the nonlinear transformation \eqref{eq:sign_LyapCoeff_WattGovernor} are depicted in magenta; they reflect well the overall distribution of the probability mass for all tested support combinations.}
	
	\label{fig:LyapCoeff_WattGovernor}
\end{figure}

\begin{figure}[t!]
	\centering
	\subfloat[\ $a_1=0, b_1=1, a_2=0, b_2=1$]{\includegraphics[width=0.46\textwidth]{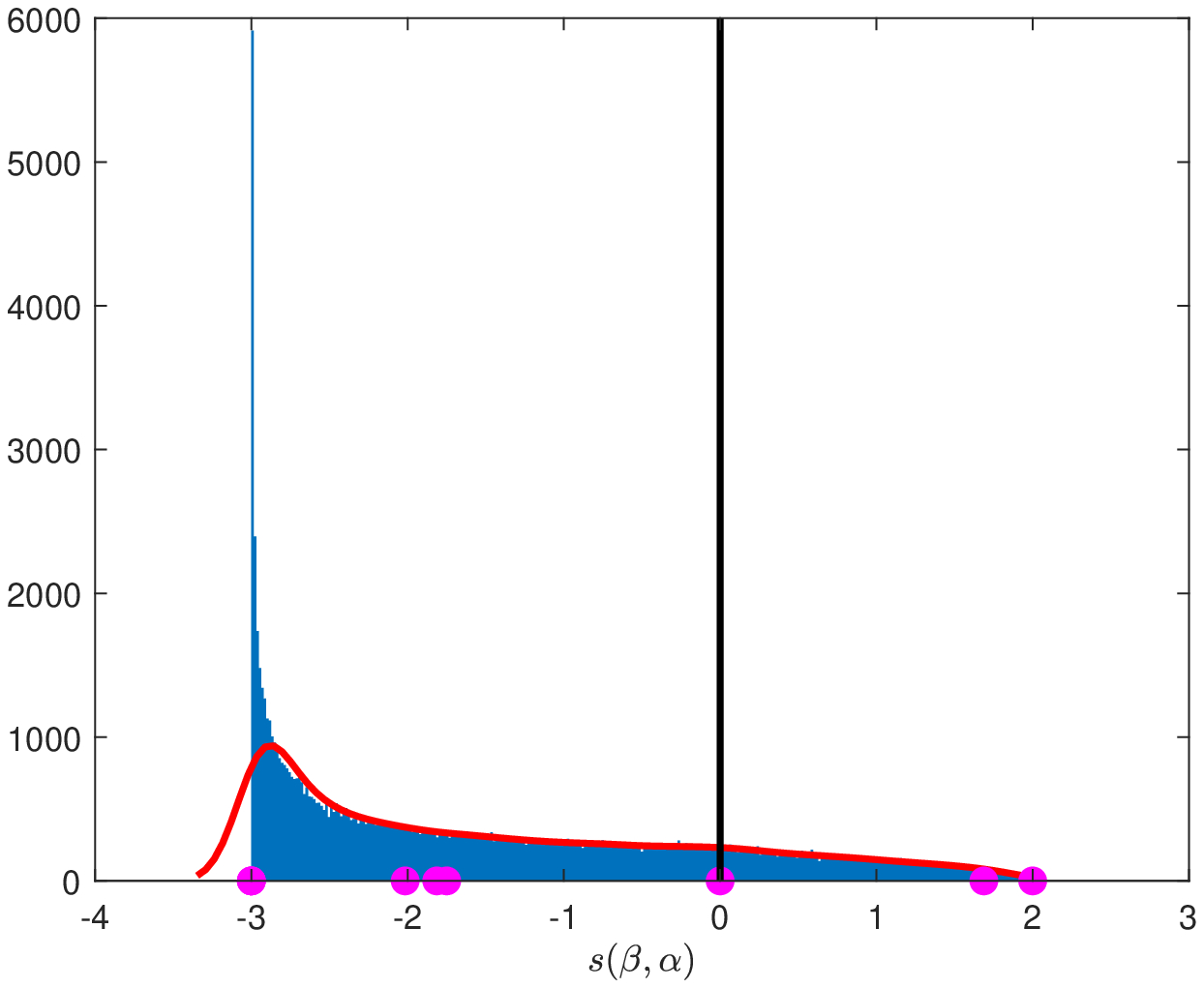}} \hspace{0.2cm}
	\subfloat[\ $a_1=0, b_1=1, a_2=0.6, b_2=1$]{\includegraphics[width=0.46\textwidth]{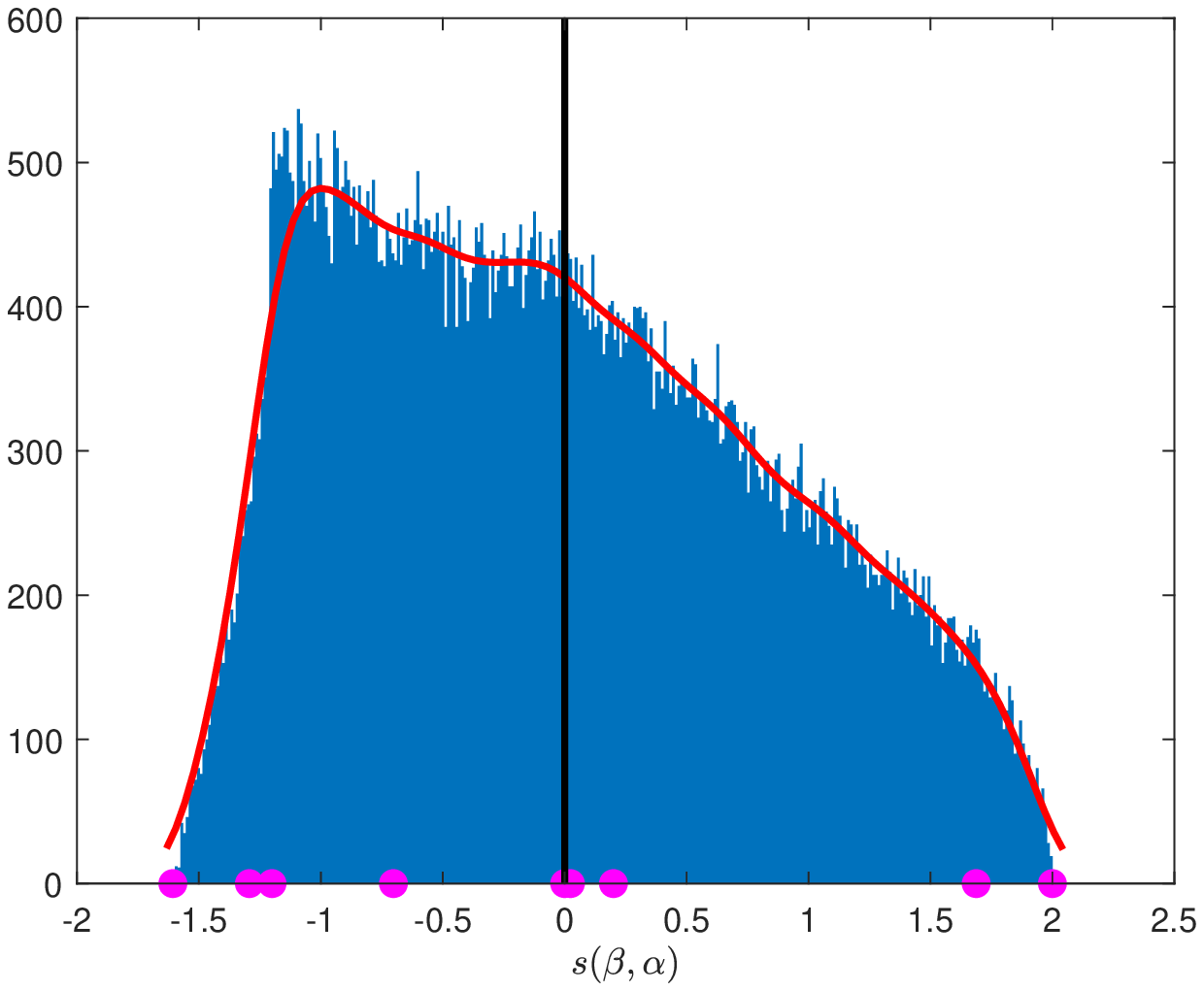}}
	
	\subfloat[\ $a_1=0, b_1=1, a_2=0.9, b_2=1$]{\includegraphics[width=0.46\textwidth]{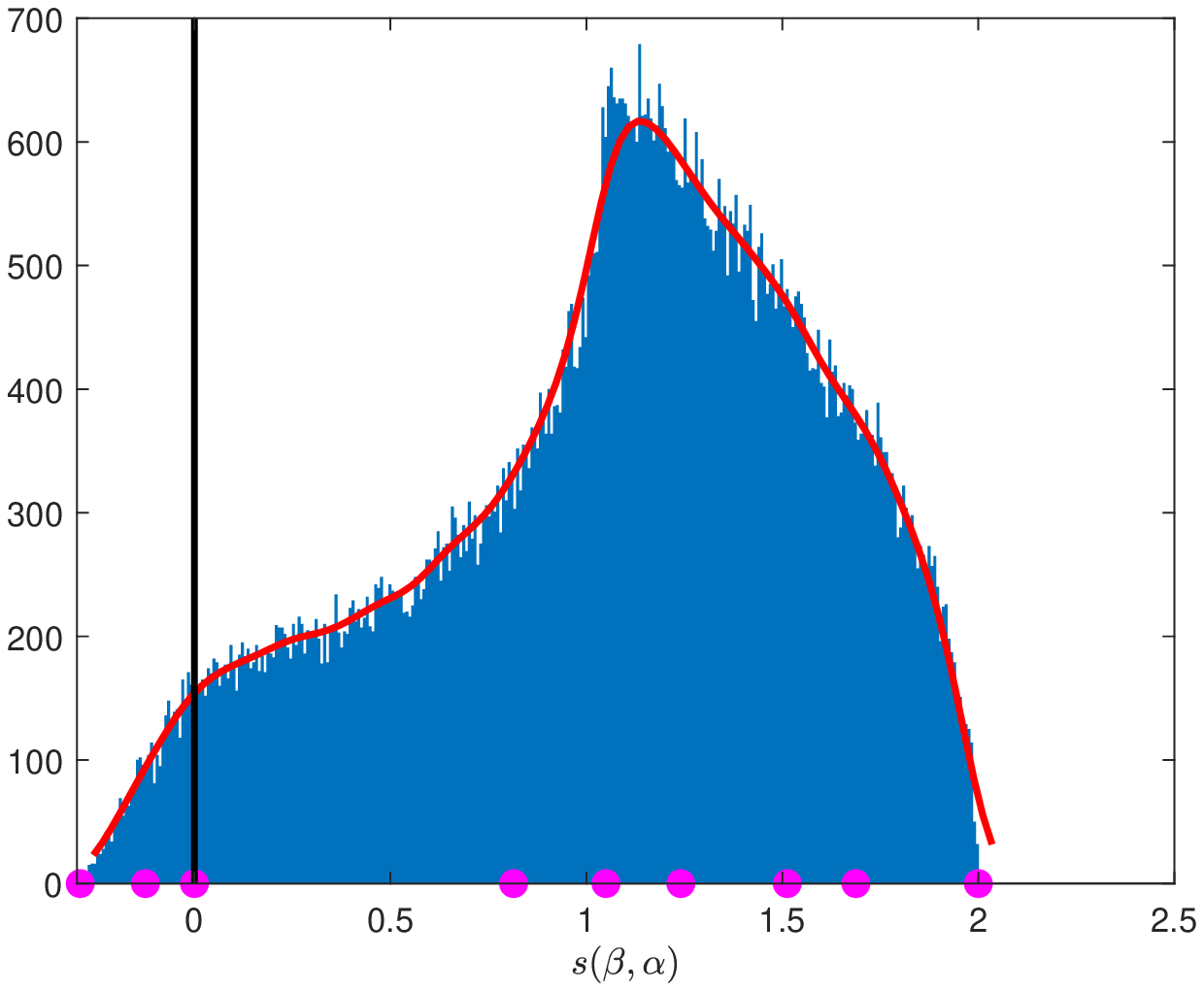}} \hspace{0.2cm}
	\subfloat[\ $a_1=0, b_1=0.2, a_2=0.7, b_2=0.95$]{\includegraphics[width=0.46\textwidth]{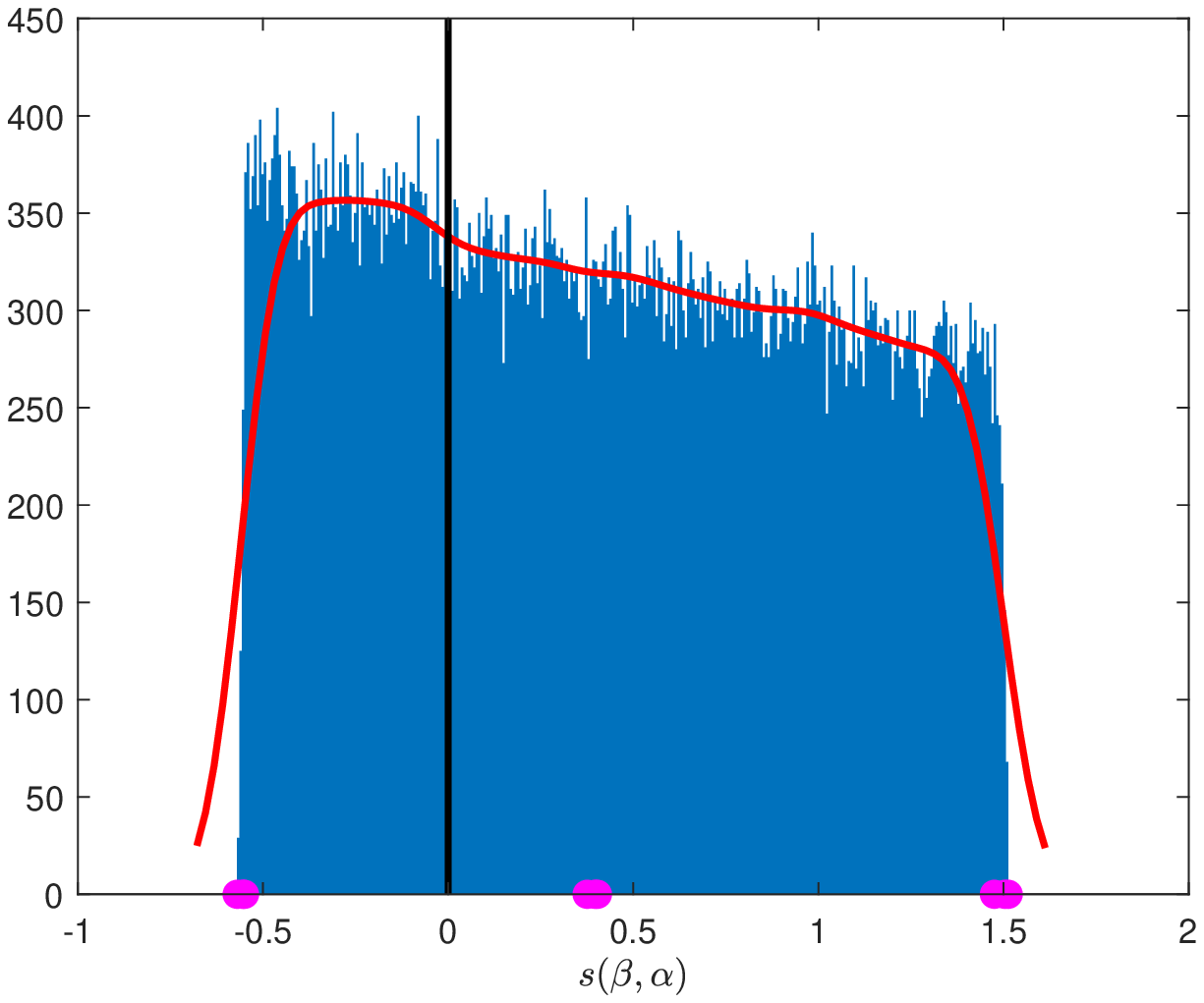}}
	
	\caption{Reproduction of Figure \ref{fig:LyapCoeff_WattGovernor} for precision 5 in sigma point selection algorithm \cite[6.2.3]{Lerner.2002} instead of 3; Probability distribution of sign of Lyapunov coefficients determined by $s(\beta,\alpha)$ in \eqref{eq:sign_LyapCoeff_WattGovernor} for different supports as in Figure \ref{fig:LyapCoeff_WattGovernor}; red line is MATLAB kernel density estimate of blue histogram. The sigma points propagated through the nonlinear transformation \eqref{eq:sign_LyapCoeff_WattGovernor} are depicted in magenta; the reflection of the overall distribution of the probability mass is more detailed for precision 5.}
	
	\label{fig:LyapCoeff_WattGovernor_prec5}
\end{figure}
This set of sigma points is then propagated via the nonlinear transformation. Based on the transformed sigma points the output statistics are calculated. The idea is to avoid an extensive Monte Carlo simulation of the bifurcation coefficient for a large number of realizations of the input parameters and apply the unscented transform approach instead.
Note that in contrast to the extended Kalman filter, no Jacobian or Hesse matrix needs to be calculated/approximated since the unscented transformation does not linearize the transformation (see e.g.\ \cite{Julier.2004}).

\bigskip

Consider the Watt governor system taken from \cite{Sotomayor.2007}, which reads as
\begin{align}
\left\{\begin{array}{ll}
\dx &= y \\
\dy &= z^2\sin{x}\cos{x}-\sin{x}-\varepsilon y \\
\dz &= \alpha(\cos{x}-\beta),
\end{array}\right. . \label{eq:WattGovernorSyst}
\end{align}
where $\alpha>0$, $0<\beta<1$, and $\varepsilon>0$.
According to \cite{Sotomayor.2007}, this system exhibits a Hopf bifurcation at $(x^\ast,y^\ast,z^\ast)=(\arccos{\beta},0,\sqrt{\nicefrac{1}{\beta}})$ for $\varepsilon_c = 2\alpha\beta^{\nicefrac{3}{2}}$. The first Lyapunov coefficient is given in \cite[equ. (60)]{Sotomayor.2007} and reads as
\begin{align}
l(\beta,\alpha,\varepsilon_c) &= -\frac{1}{2}\left(\frac{\alpha\beta^{\nicefrac{3}{2}}(1-\beta^2)(3+(\alpha^2-5)\beta^2+\alpha^4\beta^6)}{(1-\beta^2+\alpha^2\beta^4)(1-\beta^2+4\alpha^2\beta^4)}\right). \label{eq:LyapCoeff_WattGovernor}
\end{align}

Here, we assume that, in \eqref{eq:WattGovernorSyst}, $\alpha$ and $\beta$ are uncertain input parameters and the Lyapunov coefficient \eqref{eq:LyapCoeff_WattGovernor} takes the role of the relevant normal form coefficient. If the sign of \eqref{eq:LyapCoeff_WattGovernor} is positive, the Hopf bifurcation is subcritical. It is supercritical if the sign is negative. The sign of \eqref{eq:LyapCoeff_WattGovernor} is determined by
\begin{align}
s(\beta,\alpha) &= -\left(3+(\alpha^2-5)\beta^2+\alpha^4\beta^6\right). \label{eq:sign_LyapCoeff_WattGovernor}
\end{align}
Hence, we are primarily interested in the probability $P(s(\beta,\alpha)>0)$.

We illustrate the distribution of the Lyapunov coefficient via a concrete example.
Assume $\alpha\sim\UU(a_1,b_1)$ and $\beta\sim\UU(a_2,b_2)$, independent of each other. Then, we obtain the results for the distribution of $s(\beta,\alpha)$ given in  \eqref{eq:sign_LyapCoeff_WattGovernor} (see Figure \ref{fig:LyapCoeff_WattGovernor}). As we deal with two random parameters, we set $n=2$. An application of the unscented transform method with a sigma point selection of precision 3 (see equation (12) of \cite{Julier.2004}) results in $M^{UT} = 2\cdot n + 1 = 5$ deterministically calculated sigma points. The magenta realizations in Figure \ref{fig:LyapCoeff_WattGovernor} are obtained by propagating the sigma points through the nonlinear transformation $s(\beta,\alpha)$ from \eqref{eq:sign_LyapCoeff_WattGovernor}. Despite the tiny number of realizations $M^{UT}$, the magenta points reflect well the distribution of the sign of the Lyapunov coefficient for all tested support combinations $[a_1,b_1]$ for $\alpha$ and $[a_2,b_2]$ for $\beta$. This suggests that an estimation of the probability of the bifurcation type being a subcritical Hopf bifurcation based on the unscented transform might give a first rough approximation of the situation when one wants to avoid an elaborate sampling procedure. Note however that, as soon as the subcritical bifurcation type becomes a rare event, the unscented transform method is no longer suited.

The results of the unscented transformation can be improved by using a sigma point selection of precision 5 from \cite[6.2.3]{Lerner.2002}, which is summarized well in \cite{Menegaz.2015}. Therein, $M^{UT5}=2n^2+1$ points are used. A reproduction of Figure \ref{fig:LyapCoeff_WattGovernor} with sigma point selection of precision 5 instead of 3 is shown in Figure \ref{fig:LyapCoeff_WattGovernor_prec5}. A more detailed representation by the sigma points propagated through \eqref{eq:sign_LyapCoeff_WattGovernor} of the overall distribution of the probability mass is obtained.

\bigskip
As our major focus in this work are estimates of the probability of different bifurcation types under uncertainty, we conclude this work with a calculation of the probability of a subcritical Hopf bifurcation to be present in \eqref{eq:WattGovernorSyst} via a Monte Carlo simulation with $M=10^5$ samples and via the above introduced unscented transform method. Note that $M^{UT} << M$. This low number of points to evaluate might pay off most when we are confronted with a large number $n$ of random parameters in the dynamical system leading to a high-dimensional parameter space. The results for $n=2$ are given in Table \ref{tab:prob_sub_Wiggins_Ex20P6B1c}. Taking the MC estimate as our benchmark, we observe that the rough magnitude of probabilities of a subcritical Hopf bifurcation is matched for all tested support combinations of the input parameters. Furthermore, there is a significant quantitative gain in the precision of the estimate resulting from an increase in the precision within the sigma point selection method for most analyzed support combinations.
\begin{table}[h]
	\begin{adjustbox}{width=\textwidth}
		\begin{tabular}{c|c|c|c}
			& Probability sub UT prec 3 & Probability sub UT prec 5 & Probability sub MC \\
			\hline
			$a_1=0, b_1=1, a_2=0, b_2=1$ & $\nicefrac{1}{5}$ & $\nicefrac{2}{9}$ & 0.1834 \\
			\hline
			$a_1=0, b_1=1, a_2=0.6, b_2=1$ & $\nicefrac{3}{5}$ & $\nicefrac{4}{9}$ & 0.4589 \\
			\hline
			$a_1=0, b_1=1, a_2=0.9, b_2=1$ & $\nicefrac{4}{5}$ & $\nicefrac{6}{9}$ & 0.9670 \\
			\hline
			$a_1=0, b_1=0.2, a_2=0.7, b_2=0.95$ & $\nicefrac{4}{5}$ & $\nicefrac{6}{9}$ & 0.6962 \\
		\end{tabular}
	\end{adjustbox}
	\caption{Estimates of probability of subcritical (sub) Hopf bifurcation for \eqref{eq:WattGovernorSyst}; comparison to MC estimates as benchmark shows that the rough magnitude of probabilities of a subcritical Hopf bifurcation is matched for all support combinations; increasing the precision (prec) in the sigma point selection method for the unscented transformation (UT) improves probability estimates significantly for most analyzed support combinations.}
	
	\label{tab:prob_sub_Wiggins_Ex20P6B1c}
\end{table}

\section{Conclusion and Outlook}

In this work, we have developed several analytical, semi-analytical, and numerical tools to provide estimates for bifurcation types in nonlinear random systems. We have seen that the method of choice depends on the structure and probability distribution of the uncertain input parameters within the RODE \eqref{eq:rODE}. In a first step, we reduce the nonlinear propagation of the parameter uncertainty to the analysis of the normal form coefficient via reduction and transformation steps from bifurcation theory. Then a nonlinear uncertainty propagation problem remains.

The analytical approach from Section \ref{sec:anaAppUncertaintyProp} has the great advantage that the probability distribution of the normal form coefficient is calculated exactly and allows for a rigorous perturbation analysis. The drawback is that it works only for rather simple structures of parameter combinations in the bifurcation coefficient. To address this, we came up with the semi-analytical approach introduced in Section \ref{sec:semiAnaAppUncertaintyProp}. The latter grants more flexibility and the advantages of having exact distributions are carried over as much as possible. The remaining parts of the normal form coefficients are then treated via a PCE and the Mellin transform thereof. We have seen that the semi-parametric estimation procedure provides excellent estimates for the bifurcation-type probabilities. We highlight again that both approaches do not need any samples at all.

By contrast, if a detailed analysis of the nature of the nonlinear transformation is not feasible or too complex, a particular sampling-based approach has been presented in Section \ref{sec:samplingBasedAppUncertaintyProp}. Instead of a huge Monte Carlo sample, only a small number of deterministically specified input parameter combinations is propagated through the nonlinear transformation. We have seen that the estimates of the bifurcation type probability obtained could be improved by enlarging the precision in the sigma point selection algorithm. This gain might be pushed further by also using higher-order moment information in a higher-order unscented transformation framework as presented in \cite{Easley.2020}. The estimates obtained via the unscented transformation procedure can provide first indications for decision-makers to assess the risk exposure to critical transitions.

\bigskip
Our research on bifurcation type probabilities is embedded in a much broader context. Assessing the risk of critical transitions and deriving suitable countermeasures based thereon, comes with a lot of further research challenges and a holistic approach is needed. An ultimate goal for future research would be providing answers to the question: what is the probability that a critical transition in a given real-world system actually occurs?

The problem starts with building a model by combining first principle approaches in combination with data assimilation techniques. The uncertain nature of the model parameters is likely to carry over to important key characteristics of the dynamical system. The invariant sets such as equilibria or stable and unstable manifolds might become random. From a mathematical point of view, this leads to challenging tasks as for example determining the zeros of a random polynomial.

Furthermore, when we face random invariant sets, we may no longer assume knowledge about inherent bifurcation classes. This makes the reduction and transformation procedure to normal form in the first step of our method presented in this work even more involved.

\bigskip
\textbf{Acknowledgment:}
KL and CK would like to thank Jan Sieber for a very valuable discussion of the paper content and the two anonymous referees for their insightful comments that helped to improve the
manuscript. KL and CK acknowledge support of the EU within the TiPES project funded by the European Union's Horizon 2020 research and innovation programme under grant agreement No.\ 820970. This is TiPES publication \#73. CK has also been supported by a Lichtenberg Professorship of the Volks\-wagenStiftung.

\clearpage
\begin{appendix}
\section{Calculation of the Mellin transform of a PCE} \label{app:MellinPCE}
\begin{align}
& \MM\left(\sum_{n=0}^{N} \coeff_n \xi^n\right)(s) \notag \\
=& \int_{0}^{\infty} \left(\sum_{n=0}^{N} \coeff_n u^n\right)^{s-1}\rho_{\xi}(u) du \notag \\
=& \int_{0}^{\infty} \left(\sum_{n=0}^{N-1} \coeff_n u^n + \coeff_Nu^N\right)^{s-1}\rho_{\xi}(u) du \notag \\
=& \int_{0}^{\infty} \left[\sum_{k=0}^{s-1}\binom{s-1}{k}\left(\sum_{n=0}^{N-1} \coeff_n u^n \right)^{s-1-k} \coeff_N^ku^{N\cdot k} \right]\rho_{\xi}(u) du \notag \\
=& \int_{0}^{\infty} \left[\sum_{k=0}^{s-1}\binom{s-1}{k}\left(\sum_{n=0}^{N-2} \coeff_n u^n + \coeff_{N-1}u^{N-1}\right)^{s-1-k} \coeff_N^ku^{N\cdot k} \right]\rho_{\xi}(u) du \notag \\
=& \int_{0}^{\infty} \left[\sum_{k=0}^{s-1}\binom{s-1}{k}\left(\sum_{l=0}^{s-1-k}\binom{s-1-k}{l}\left(\sum_{n=0}^{N-2} \coeff_n u^n\right)^{s-1-k-l}\coeff_{N-1}^l u^{(N-1)\cdot l}\right) \coeff_N^k u^{N\cdot k} \right]\rho_{\xi}(u) du \notag \\
=& \int_{0}^{\infty} \left[\sum_{k=0}^{s-1}\binom{s-1}{k}\left(\sum_{l=0}^{s-1-k}\binom{s-1-k}{l}\left(\sum_{n=0}^{N-3} \coeff_n u^n + \coeff_{N-2}u^{N-2}\right)^{s-1-k-l}\coeff_{N-1}^l u^{(N-1)\cdot l}\right) \coeff_N^k u^{N\cdot k} \right]\rho_{\xi}(u) du \notag \\
=& \int_{0}^{\infty} \left[\sum_{k=0}^{s-1}\binom{s-1}{k}\left(\sum_{l=0}^{s-1-k}\binom{s-1-k}{l}\left(\sum_{m=0}^{s-1-k-l}\binom{s-1-k-l}{m}\left(\sum_{n=0}^{N-3} \coeff_n u^n\right)^{s-1-k-l-m} \coeff_{N-2}^m u^{(N-2)\cdot m}\right)\right.\right. \notag \\
&\left.\left.\cdot\coeff_{N-1}^l u^{(N-1)\cdot l} \vphantom{\sum_{k=0}^{s-1}} \right) \coeff_N^k u^{N\cdot k} \vphantom{\sum_{k=0}^{s-1}}\right]\rho_{\xi}(u) du \notag \\
& \cdots \notag \\
=&  \int_{0}^{\infty} \left[\sum_{k=0}^{s-1}\binom{s-1}{k}\left(\sum_{l=0}^{s-1-k}\binom{s-1-k}{l}\left(\sum_{m=0}^{s-1-k-l} \binom{s-1-k-l}{m}\left(\cdots (\coeff_0 +\coeff_1 u\right)^{s-1-k-l-m-\cdots} \cdots \coeff_{N-2}^m u^{(N-2)\cdot m}\right)\right.\right. \notag \\
&\left.\left.\cdot\coeff_{N-1}^l u^{(N-1)\cdot l} \vphantom{\sum_{k=0}^{s-1}} \right) \coeff_N^k u^{N\cdot k} \vphantom{\sum_{k=0}^{s-1}}\right]\rho_{\xi}(u) du \notag \\
=&  \int_{0}^{\infty} \left[\sum_{k=0}^{s-1}\binom{s-1}{k}\left(\sum_{l=0}^{s-1-k}\binom{s-1-k}{l}\left(\sum_{m=0}^{s-1-k-l} \binom{s-1-k-l}{m}\left( \vphantom{\sum_{k=0}^{s-1}} \cdots \right.\right.\right.\right. \notag \\
&\left.\left.\left.\left. \left(\sum_{j=0}^{s-1-k-l-m-\cdots}\binom{s-1-k-l-m-\cdots}{j}\coeff_0^{s-1-k-l-m-\cdots-j} \coeff_1^j u^{(N-(N-1))\cdot j}\right) \cdots \vphantom{\sum_{k=0}^{s-1}} \right)\right.\right.\right.\notag\\
&\left.\left.\left.\coeff_{N-2}^m u^{(N-2)\cdot m}\vphantom{\sum_{k=0}^{s-1}}\right)\cdot\coeff_{N-1}^l u^{(N-1)\cdot l} \vphantom{\sum_{k=0}^{s-1}} \right) \coeff_N^k u^{N\cdot k} \vphantom{\sum_{k=0}^{s-1}}\right]\rho_{\xi}(u) du \label{eq:MellinSumInNestedSummation}
\end{align}
The nested summation in \eqref{eq:MellinSumInNestedSummation} can be rewritten as
\begin{align*}
\MM\left(\sum_{n=0}^{N} \coeff_n \xi^n\right)(s) &= \int_{0}^{\infty} \left[\sum_{i=0}^{N\cdot(s-1)} \hat{c}_i(s)\cdot u^i \right] \rho_{\xi}(u) du = \sum_{i=0}^{N\cdot(s-1)} \hat{c}_i(s)\cdot  \int_{0}^{\infty} u^i \rho_{\xi}(u) du \\
&= \sum_{i=0}^{N\cdot(s-1)} \hat{c}_i(s) \cdot \MM(\xi)(i+1), 
\end{align*}
where $\hat{c}_i(s)$ denote the cumulated coefficients for the powers of $u$, which depend on the chosen evaluation $s$ of the Mellin transform.

\section{Polynomial approximation of the probability distribution of the bifurcation coefficient} \label{app:probBifGMM_poly}

We will now perform a polynomial approximation the bifurcation coefficient $X$ in the reduced Lorenz system \eqref{eq:redLorenzSyst} and explain why, in our bifurcation context, we prefer the approach via the method of moments to fit a Gaussian mixture distribution, for which we already showed a first result in Figure~\ref{fig:PDF_approx_teaserLorenzSyst}.
Carrying over the estimation techniques from Sections \ref{subsubsec:nonparaEst_PolynomialEstPDF} and \ref{subsubsec:semiParaEst_MethodOfMoments} to non-standard PDFs, such as the one of the bifurcation coefficient $X$ in the reduced Lorenz system \eqref{eq:redLorenzSyst}, turns out to be a tricky task.

\begin{figure}[t]
	\subfloat[\label{fig:histVsPDFapprox_LorenzRV_Uniform46_Gamma81_PolyApprox_chosenSupp01_N2_numMoms7}\ {Support $[0,1]$}]{\begin{overpic}[width=0.5\textwidth]{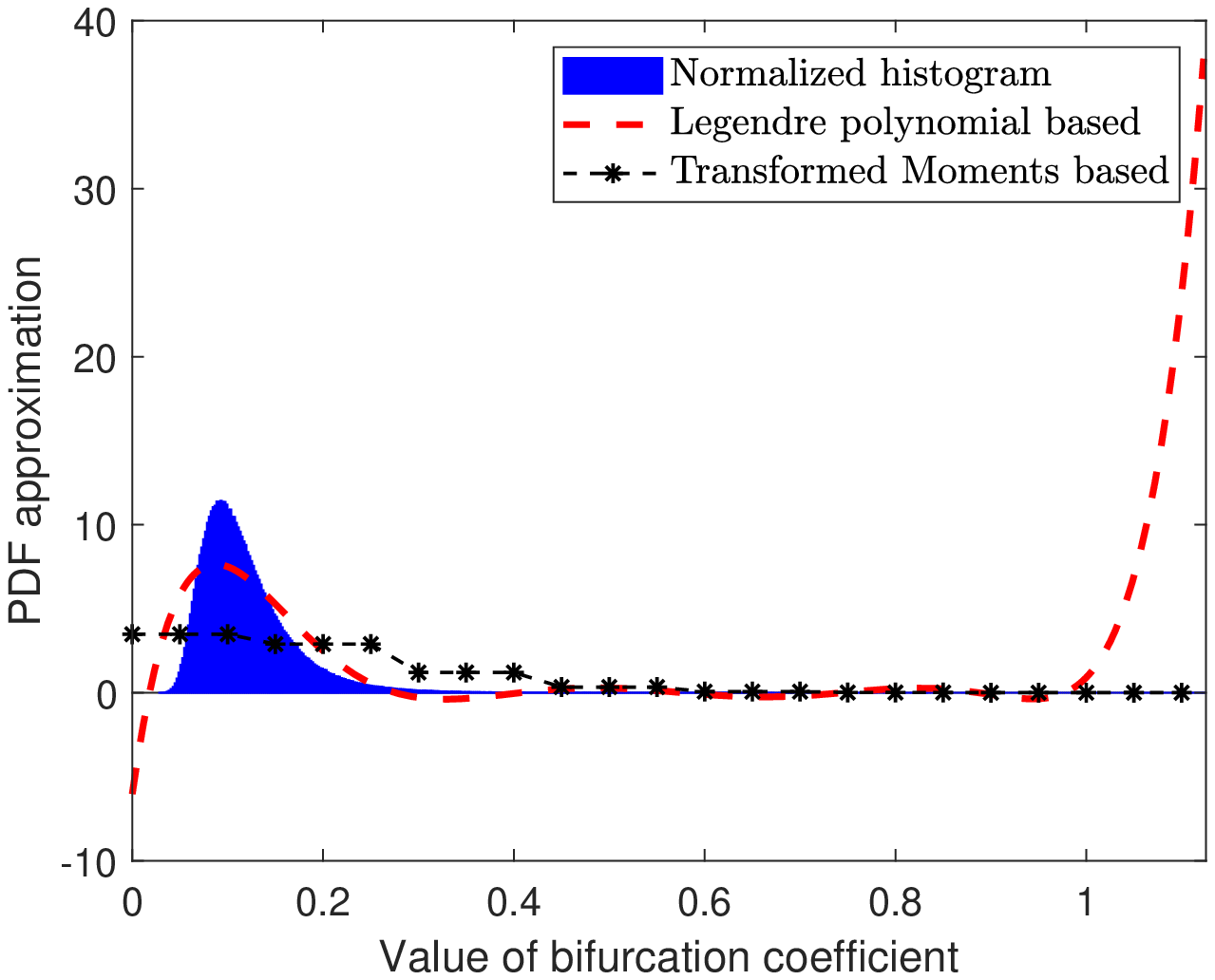}
	\end{overpic}
	}
	\subfloat[\ {Support $[0,0.5]$}]{\begin{overpic}[width=0.5\textwidth]{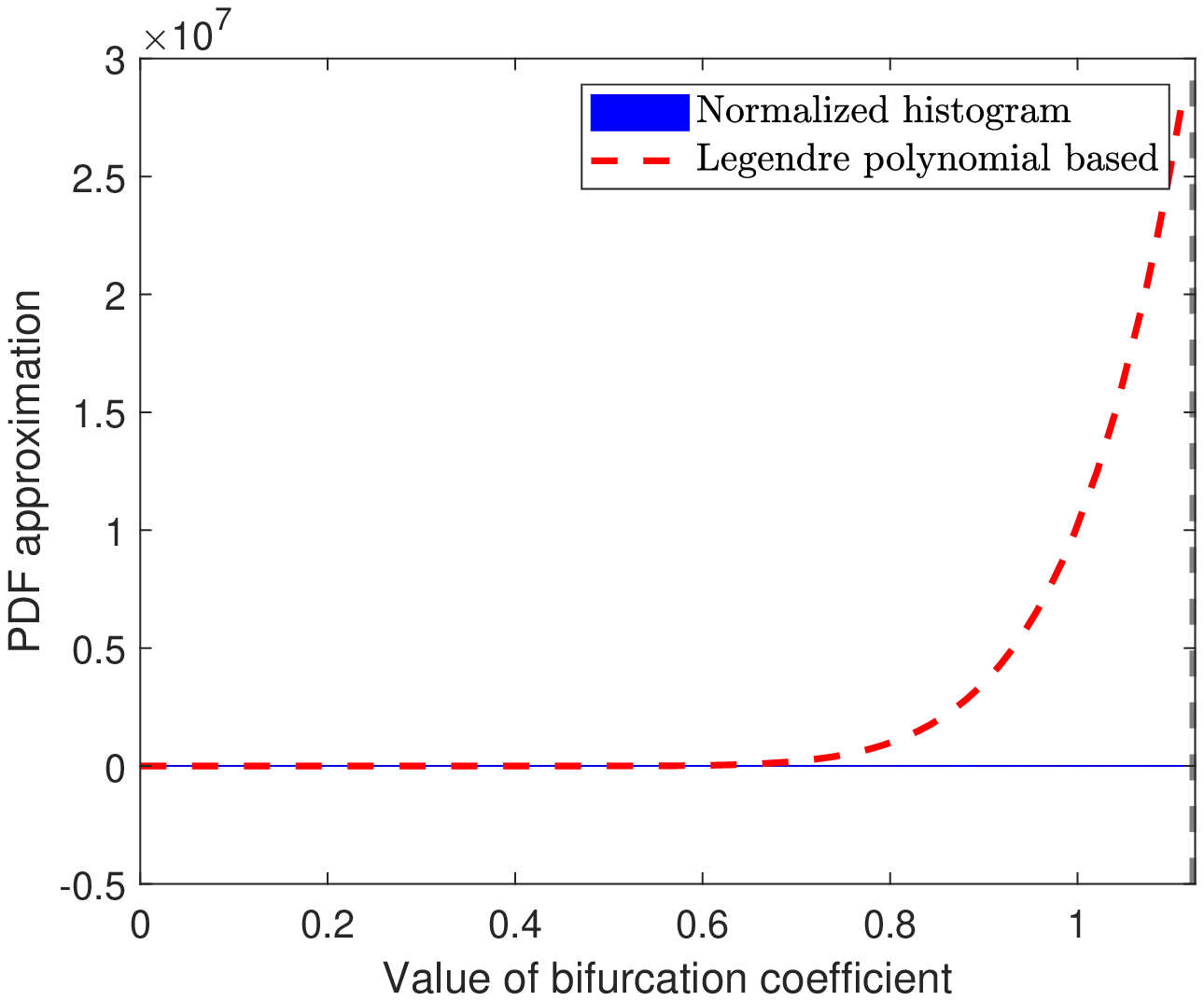}\end{overpic}}
	
	\caption{Reconstruction of the PDF of $X$ in reduced Lorenz system \eqref{eq:redLorenzSyst} via polynomial approximation for PS A; the quality of the approximation heavily depends on the specified support $[a,b]$ of the approximation: whereas in the left figure, PDF approximation \eqref{eq:PDF_approx_legPol} works reasonably for $a=0$, $b=1$, it fails completely in the right figure for $a=0$, $b=0.5$.}
	
	\label{fig:histVsPDFapprox_LorenzRV_Uniform46_Gamma81_PolyApprox_N2_numMoms7}
\end{figure}

Firstly, the support of the bifurcation coefficient is not known exactly. Although the input distributions are known with corresponding supports, a direct analytical transfer to the support of $X$ might not be feasible due to the nonlinear transformation structure. But the support $[a,b]$ crucially affects the polynomial approximations \eqref{eq:PDF_approx_legPol}, \eqref{eq:PDF_approx_monic}, and \eqref{eq:PDF_approx_trafoMom} of the PDF.
This can be seen in Figure \ref{fig:histVsPDFapprox_LorenzRV_Uniform46_Gamma81_PolyApprox_N2_numMoms7} for the PDF approximation based on Legendre polynomials \eqref{eq:PDF_approx_legPol}: whereas for a chosen support of $[0,1]$, the approximation \eqref{eq:PDF_approx_legPol} recovers well the PDF of the Lorenz coefficient, the approximation fails for setting the support to $[0,0.5]$.

Secondly, the moments used are no longer exact but include an approximation in terms of the PCE \eqref{eq:PCEnonlinTrafo_individual}. This might harm the approximation \eqref{eq:PDF_approx_legPol} of the PDF as it is stated in \cite{Provost.2005} that the method is intended to be used with exact moments.

Thirdly, Figure \ref{fig:histVsPDFapprox_LorenzRV_Uniform46_Gamma81_PolyApprox_chosenSupp01_N2_numMoms7} shows that -- as already expected -- the number of moments that can be reasonably well calculated ($\nmoms=7$ here) is not high enough to provide an accurate approximation via the method of transformed moments. Hence, we also have to pay attention to how many moments of the bifurcation coefficient can be reasonably calculated due to a possibly bad scaling in the numerics for the moments in equation \eqref{eq:MellinTransform_PCE_Uniform}.

The method based on monic orthogonal polynomials somehow fails to be in an acceptable range in the reconstruction of the bifurcation coefficient $X$ in the reduced Lorenz system \eqref{eq:redLorenzSyst}.

It is mentioned in \cite{Mnatsakanov.2008a} that an unstable behavior might be observed when the higher-order moments are used in the approximation. In our numerical results for the polynomial-based approximations, we face instability issues as well.
To increase stability, in \cite{Tagliani.2003}, in the context of the maximum entropy method, it is suggested to work with fractional moments.

Instead here, we test the performance of a completely different type of approximation, namely the semi-parametric approach via GMMs that we introduced in Section \ref{subsubsec:semiParaEst_MethodOfMoments}.
\end{appendix}

\bibliographystyle{siam}
\bibliography{myLiterature}

\end{document}